\documentclass{hooml2022}

\usepackage{amsfonts}
\usepackage{microtype}

\usepackage{dsfont}
\newcommand{\E}[1]{\mathds{E} \left[#1\right]}

\usepackage{color}

\title[Effects of momentum scaling for SGD]{Effects of momentum scaling for SGD}

\optauthor{\Name{Dmitry A. Pasechnyuk}\Email{dmivilensky1@gmail.com}\AND\Name{Alexander Gasnikov} \Email{gasnikov@yandex.ru}\AND\Name{Martin Tak{\'a}{\v{c}}}\Email{martin.taki@gmail.com}}

\usepackage{algorithmic}
\usepackage{algorithm}
\usepackage{wrapfig}
\begin{document}

\maketitle

\setlength{\abovedisplayskip}{3pt}
\setlength{\belowdisplayskip}{3pt}

\begin{abstract}
The paper studies the properties of stochastic gradient methods with preconditioning. We focus on momentum updated preconditioners with momentum coefficient $\beta$. Seeking to explain practical efficiency of scaled methods, we provide convergence analysis in a norm associated with preconditioner, and demonstrate that scaling allows one to get rid of gradients Lipschitz constant in convergence rates. Along the way, we emphasize important role of $\beta$, undeservedly set to constant $0.99...9$ at the arbitrariness of various authors. Finally, we propose the explicit constructive formulas for adaptive $\beta$ and step size values.
\end{abstract}
\section{Literature review}
Preconditioning is long and widely known practice in numerical methods and mathematical programming \cite{davidon1959variable,broyden1967quasi,dennis1977quasi}. With the recent surge of interest to statistical learning applications, there were proposed methods applicable to finite-sum function minimization \cite{amari1998natural,duchi2011adaptive,li2018preconditioner}. In combination with momentum technique \cite{polyak1964some}, preconditioning gave rise to adaptive methods, extensively used in applications \cite{zeiler2012adadelta,schaul2013no,kingma2014adam,ene2022adaptive}. Recently, there has been a return to full-matrix methods using momentum \cite{yao2021adahessian,jahani2021doubly,sadiev2022stochastic,beznosikov2022scaled}. This gave significant benefit in practice, but in theory there still was no established with global convergence of SGD \cite{berahas2021quasi,rodomanov2022rates}. Following \cite{scheinberg2016practical}, we analyse convergence of preconditioned gradient descent in preconditioner associated norm to get rid of gradient Lipschitz constant $L$. This allows us to determine step size independent on $L$ and adaptive momentum parameter $\beta$. We also show the starting acceleration of convergence as in \cite{hanzely2021damped}.

\section{One-step effects}
Firstly, we consider optimization problem $\min_{x \in \mathbb{R}^n} f(x)$, where $f: \mathbb{R}^n \rightarrow \mathbb{R}$ is non-convex continuous function. Every additional requirement on $f$ is introduced in appropriate place of the text where it is needed to simplify reasoning. Let us start with considering the simplest preconditioned method \texttt{Scaled SGD}
\[
    x_{t+1} = x_t - \eta_t P_t^{-1} g_t
\]
with variable preconditioner $P_t \in \mathbb{S}^{n}_{++}$. 
Our goal for the nearest narration is to estimate the rate of function decreasing. For this purpose, we need to operate with a majorant of $f$. We assume that $f$ 
has Lipschitz continuous gradient, that is  
$
    f(x_{t+1}) \leqslant f(x_t) + \langle \nabla f(x_t), x_{t+1} - x_t \rangle + \tfrac{L_t}{2} \|x_{t+1} - x_t\|^2,
$
on each of segment $[x_t, x_{t+1}]$, $t=1,2,...$ of methods trajectory with a corresponding $L_t = L({x_t,x_{t+1}})$. This assumption is weaker than uniform $L$-smoothness. This point of view allows to get rid of uniform constant $L$, but still does not tell anything about its value, because properties of the norm $\|\cdot\|$ are not used anyhow. Then, if we go to norm $\|x\|_{P_t} = \langle P_t x, x \rangle^{1/2}$ associated with $P_t$, we have
\[
    f(x_{t+1}) \leqslant f(x_t) - \eta_t \langle \nabla f(x_t), P_t^{-1} g_t \rangle + \tfrac{L_t \eta_t^2}{2} \|g_t\|_{P_t}^{*\textsuperscript{2}},
\]
for some $L_t = L(x_t, x_{t+1}, P_t)$ which is believed to be smaller than previous $L(x_t, x_{t+1})$ if preconditioner is proper. The best case is when $P_t = \nabla^2 f(x_t)$ which implies that condition number is close to $1$ in some vicinity of $x_t$, so it is also common to estimate convergence rate in $\|\cdot\|_{\nabla^2 f(x_t)}$ norm.

Another upper bound for $f$ comes from $M$-Lipschitz continuity of $\nabla^2 f$:
\[
    f(x_{t+1}) \leqslant f(x_t) - \eta_t \langle \nabla f(x_t), P_t^{-1} g_t \rangle + \tfrac{\eta_t^2}{2} \langle \nabla^2 f(x_t) P_t^{-1} g_t, P_t^{-1} g_t \rangle + \tfrac{M \eta_t^3}{6} \|P_t^{-1} g_t\|^3,
\]
which can be equivalently rewritten as
\[
    f(x_{t+1}) \leqslant f(x_t) - \eta_t \langle \nabla f(x_t), P_t^{-1} g_t \rangle + \tfrac{\eta_t^2}{2} {\|g_t\|_{P_t \left[\nabla^2 f(x_t)\right]^{-1} P_t^\top}^{*\textsuperscript{2}}} + \tfrac{M \eta_t^3}{6} \|P_t^{-1} g_t\|^3.
\]
Thus, we have got rid of $L_t$ and replaced it with term without any uniform constant, but in different $P_t \left[\nabla^2 f(x_t)\right]^{-1} P_t^\top$ norm, plus additional cubic term. The new form of quadratic term gives us an opportunity to obtain the explicit replacement for $L_t$. 

Indeed, we can think of $\nabla^2 f(x_t) P_t^{-1}$ as an inexactness of $P_t$, its closeness to the Hessian value, which can be bounded as follows:
\[
    \boxed{\nabla^2 f(x_t) P_t^{-1} \preccurlyeq (1+\Delta_t) I}
\]
where ideally $0 \leqslant \Delta_t \ll 1$.
Hereinafter, $\preccurlyeq$ is used to compare arbitrary matrices with matrices of the form $\text{const}\cdot I$, so we can define it as follows: $A \preccurlyeq bI$ means $\lambda_{\max}(A) \leq b$. Then, we have
\[
    f(x_{t+1}) \leqslant f(x_t) - \eta_t \langle \nabla f(x_t), P_t^{-1} g_t \rangle + \tfrac{(1 + \Delta_t) \eta_t^2}{2} {\|g_t\|_{P_t}^{*\textsuperscript{2}}} + \tfrac{M \eta_t^3}{6} \|P_t^{-1} g_t\|^3.
\]
Thus, we turned 2nd-order term with constant $L$ into 2nd-order term with constant $(1 + \Delta_t)$, which is slightly greater than $1$, and cubic term, so the behaviour of preconditioned gradient descent is as close to that of Newton method as $P_t$ is to $\nabla^2 f(x_{t})$.

From now, let us consider only preconditioners of the form 
\begin{equation}\label{momentum}
    P_{t+1} = \beta_{t+1} P_t + (1-\beta_{t+1}) d_{t+1},
\end{equation}
where $d_{t+1} = \text{diag}\left(\nabla^2 f(x_{t+1})\right)$ and $P_0 = I$. 
One can use any other proper update instead of $d_{t+1}$, which preserves positive definiteness of $P_t$ (if $f$ is non-convex, positive truncation should be applied to $d_{t+1}$, see \cite{paternain2019newton}). We assume that $d_t$ is a good approximation of Hessian, so that $P_t$ is maintained to be close to Hessian. In the case of diagonal $d_t$, we assume that $\nabla^2 f$ is almost diagonal. By introducing two more inexactness relating preconditioner and Hessian to update term $d_t$
\[
    \boxed{\nabla^2 f(x_t) d_t^{-1} \preccurlyeq (1+\sigma) I}\qquad\boxed{(1-\delta_t^-) I \preccurlyeq P_t d_t^{-1} \preccurlyeq (1+\delta_t^+) I,}
\]
we get the opportunity to express each one of $\sigma, \delta$ and $\Delta$ through other ones. Estimating a local Lipschitz constant of $f$ after scaling, we obtain the following proposition that bound $\Delta$.
\begin{proposition}\label{Delta} For preconditioner updated in accordance with \eqref{momentum}, inexactness $\Delta_t$ depends on inexactness $\delta_t^-$ as follows
    $
    \Delta_t \leqslant \tfrac{1+\sigma}{1 - \min\{\delta_t^-, \beta_t\}} - 1.
    $
\end{proposition}
Note, that the dependency on $\beta_t$ is hidden behind $\delta_t^-$. But $\delta_t^-$ grows with $\beta_t$, $\delta_t^-=0$ for $\beta_t=0$ and $\delta_t^- \in [0, 1)$, so our new bound on $\Delta_t$ behaves similarly to the previous one.
\\
It is obvious that $\beta_t=0$ is the best choice for the case $g_t = \nabla f(x_t)$. Otherwise, small $\beta_t$ also leads to additional penalty on the variation of $P_t$. If $g_t$ is unbiased estimator of $\nabla f(x_t)$, this penalty goes with variation $\|s\|_{P_t}^{*\textsuperscript{2}}$, where $s = g_t - \nabla f(x_t)$. Since we consider $\|\cdot\|_{P_t}$ norm, penalty appears when we go from $\|\cdot\|_{P_t}$ to $\|\cdot\|_{P_{t+1}}$ norm. 
\begin{proposition}\label{variance}
For any $s \in \mathbb{R}^n$, it holds that $\|s\|_{P_{t+1}}^{*\textsuperscript{2}} \leqslant  (1 + \tfrac{1 - \beta_{t+1}}{1/\delta_t^+ +\beta_{t+1}} ) \|s\|_{P_t}^{*\textsuperscript{2}}$.
\end{proposition}
Appearing factor is a penalty. For fixed $\delta_t^+$, it decreases inversely proportional to $\beta_t$, have maximum in $\beta_t=0$ with value $1+\delta_t^+$ and minimum in $\beta_t=1$ with value $1$. If $\delta_t^+$ depends on $\beta_t$, penalty may behave in a more complicated way, but still pushes the best value of $\beta_t$ away from zero.
\\
It remains to relate $\delta_t$ and $\beta_t$ to obtain a descent lemma depending only on a choice of $\beta_t$. There are two ways to do this: assuming, that Hessian changes only a little from iteration to iteration, or not.
\begin{proposition}\label{delta_convex}
For any $(d_t \succcurlyeq 0)_{t=0}^\infty$, it holds that $\delta_{t+1}^+ = \beta_{t+1} \kappa_{t+1}, \delta_{t+1}^- = \beta_{t+1} \chi_{t+1}$,
where $\kappa_t = \left[\tfrac{\max_i {[P_{t-1}]_{ii}}}{\min_i {[d_t]_{ii}}} - 1\right]_+, \chi_t = \left[1 - \tfrac{\min_i {[P_{t-1}]_{ii}}}{\max_i {[d_t]_{ii}}}\right]_+$.
\end{proposition}
\begin{proposition}\label{delta_concordance}
If $f$ is strong self-concordant \cite{rodomanov2021greedy}, that is $\forall x,y,z,w \in \mathbb{R}^n$
\begin{equation*}
    \text{diag}\left(\nabla^2 f(y) - \nabla^2 f(x)\right) \preccurlyeq N \|y - x\|_{\text{diag}\left(\nabla^2 f(z)\right)} \text{diag}\left(\nabla^2 f(w)\right)
\end{equation*}
for some $N > 0$, then it holds that $\delta_{t+1}^+ \leqslant \beta_{t+1}  [\delta_t^+ + \delta_t^+ \sqrt{1 + \delta_t^+} N \eta_t \|g_t\|_{P_t}^* - 1 ]_+$.
\end{proposition}
We can finally estimate the descent on the one step, depending only on $\beta$ and not on any of inexactnesses.
\begin{theorem}[Descent Lemma]\label{descent_lemma} Point $x_{t+1}$ generated by \texttt{Scaled SGD} on iteration $t$ satisfies
\begin{align*}
    \E{f(x_{t+1})} &\leqslant f(x_t) - \tfrac{\eta_t}{2} \|\nabla f(x_t)\|_{P_t}^{*\textsuperscript{2}} + \tfrac{\eta_t}{2}  ( (\tfrac{\eta_t (1+\sigma)}{1 - \beta_t \chi_t} )^2 + \tfrac{\eta_t (1+\sigma)}{1 - \beta_t \chi_t} - 1 ) {\|g_t\|_{P_t}^{*\textsuperscript{2}}}
    \\
    &\quad+ \tfrac{\eta_t}{2} 
      (1 + \tfrac{(1 - \beta_{t+1}) \beta_t}{1/\kappa_t+\beta_t \beta_{t+1}} ) \|g_t - \nabla f(x_t)\|_{P_t}^{*\textsuperscript{2}} + \tfrac{M^\prime \eta_t^3}{6}  (\tfrac{1+\sigma}{1 - \beta_t \chi_t} )^{3/2} \|g_t\|_{P_t}^{*\textsuperscript{3}}.
\end{align*}
\end{theorem}
    
\section{Cumulative effects}
Further, we consider finite-sum optimization problem $\min_{x \in \mathbb{R}^n} f(x) := \frac{1}{m} \sum_{i=1}^m f_i(x)$. Such a form of objective function became widely used due to the recent surge of interest in statistical learning applications.

We assume that Lipschitz constants are close to $1$, in view of what we have demonstrated, but further we do not specify their values. This allows us to claim that following results stay valid not only for the preconditioned methods, but for any gradient method if the difference between $L_t$, $t=1,2,...$ is significant.

\begin{theorem}\label{convergence}
If $\eta_t \leqslant \min \{\tfrac{\alpha p}{3}, \frac{3}{4}\tfrac{p}{5p+1} \} \tfrac{1}{L_t}$, sequence of the points generated by \texttt{Scaled L-SVRG} satisfies
$
\E{\|\nabla f(\overline{x}_T)\|^{*\textsuperscript{2}}_{P_t}} \leqslant \tfrac{4}{\alpha p}\tfrac{\overline{L}}{T} 
 [f(x_0) - f(x_*) + 2 \Gamma \textstyle{\sum}_{t=2}^{T+1} L_t \E{(1-\beta_t) \|x_t - y_t\|_2^2}  ],
$
where $\overline{L} = \frac{T}{\frac{1}{L_1} + \dotsi + \frac{1}{L_T}}$ is a harmonic average of $L_t$, $t=1,...,T$.
\end{theorem}
Note, that harmonic average, which appears in convergence rate estimates too rarely by the way, has the wonderful property that it is minimum-dominated. In particular, $\frac{T}{\sum_{t=1}^T 1/L_t} \leqslant T \min_{t=1,...,T} L_t$. 

\begin{wrapfigure}{R}{0.6\textwidth}
\vskip-10pt
    \begin{minipage}{0.6\textwidth}
\begin{algorithm}[H]
 \caption{\texttt{Scaled L-SVRG}}
\begin{algorithmic}
  \STATE {\bf Data:} {$p \in (0, 1), \{\eta_t > 0\}_{t=0}^{\infty}, \{0 < \beta_t < 1\}_{t=0}^{\infty}, x_0$}
  \STATE $P_0 = I$, $y_0 = x_0$\;
  \FOR{$t \geqslant 0$}
   
   \STATE Draw $i_t \in \{1, ..., m\}$ from uniform distribution\;
\STATE    $g_t = \nabla f_{i_t}(x_t) - \nabla f_{i_t}(y_t) + \nabla f(y_t)$\;
    \STATE $y_{t+1} = \left\{\begin{array}{ll}
            y_t & \text{with probability } p\\
            x_t & \text{with probability } 1-p
        \end{array}\right.$\;
    \STATE $x_{t+1} = x_t - \eta_t P_t^{-1} g_t$\;
    \STATE Draw $z_t \in \{-1, 1\}^n$ from Rademacher distribution\;
    \STATE $P_{t+1} = \beta_{t+1} P_t + (1-\beta_{t+1}) \text{diag}\left(z_t \circ \nabla^2 f(x_{t+1}) z_t\right)$\;
 
  \ENDFOR 
\end{algorithmic}
\end{algorithm}
\end{minipage}
\vskip-10pt
\end{wrapfigure}
This implies that it does not matter how big is one of the $L_t$ (or all of them, except one), their harmonic average will be bounded and stuck to minimum of the elements being averaged, even if some of those variable $L_t$ are infinite. Moreover, harmonic average is always less than arithmetic mean. So, we say that local Lipschitz constants are very well-aggregated in final convergence rate. If algorithm maintain $L_t$ close to one for each $t=1,2,...$, final convergence rate is almost independent on global Lipschitz constant!

Note that obtained dependence of the norm of gradient after $T$ iterations of algorithm on $\beta_t$ allows one to choose $\beta_t$ so that error term in Theorem~\ref{convergence} can take a small value. At least, we should guarantee its boundedness for every $T$. This can be achieved by the several ways: by choosing $\beta_t$ depending on pre-known number of iterations $T$, or by making it dependent on some hyperparameter sequence. In the first case, we notice that for $\beta_t = 1 - 1/(T L_t \|x_t - y_t\|_2^2)$, error term is equal to $2 \Gamma$ and does not affect the convergence rate. But with this approach $\beta_t$ values can be chosen too close to $1$, which is not efficient in practice. To prevent this, we can bound the error term with the series with upper limit $+\infty$. Then, it is sufficient to choose $\beta_t = 1 - a_t/(L_t \|x_t - y_t\|_2^2)$, where $a_t > 0$ and $\sum_{t=2}^{+\infty} a_t$ is converging, so that error term is bounded.

\section{Synthesis}
\begin{wrapfigure}{r}{0.37\columnwidth}
    \centering
    \vskip-20pt
     {\includegraphics[width=0.35\textwidth]{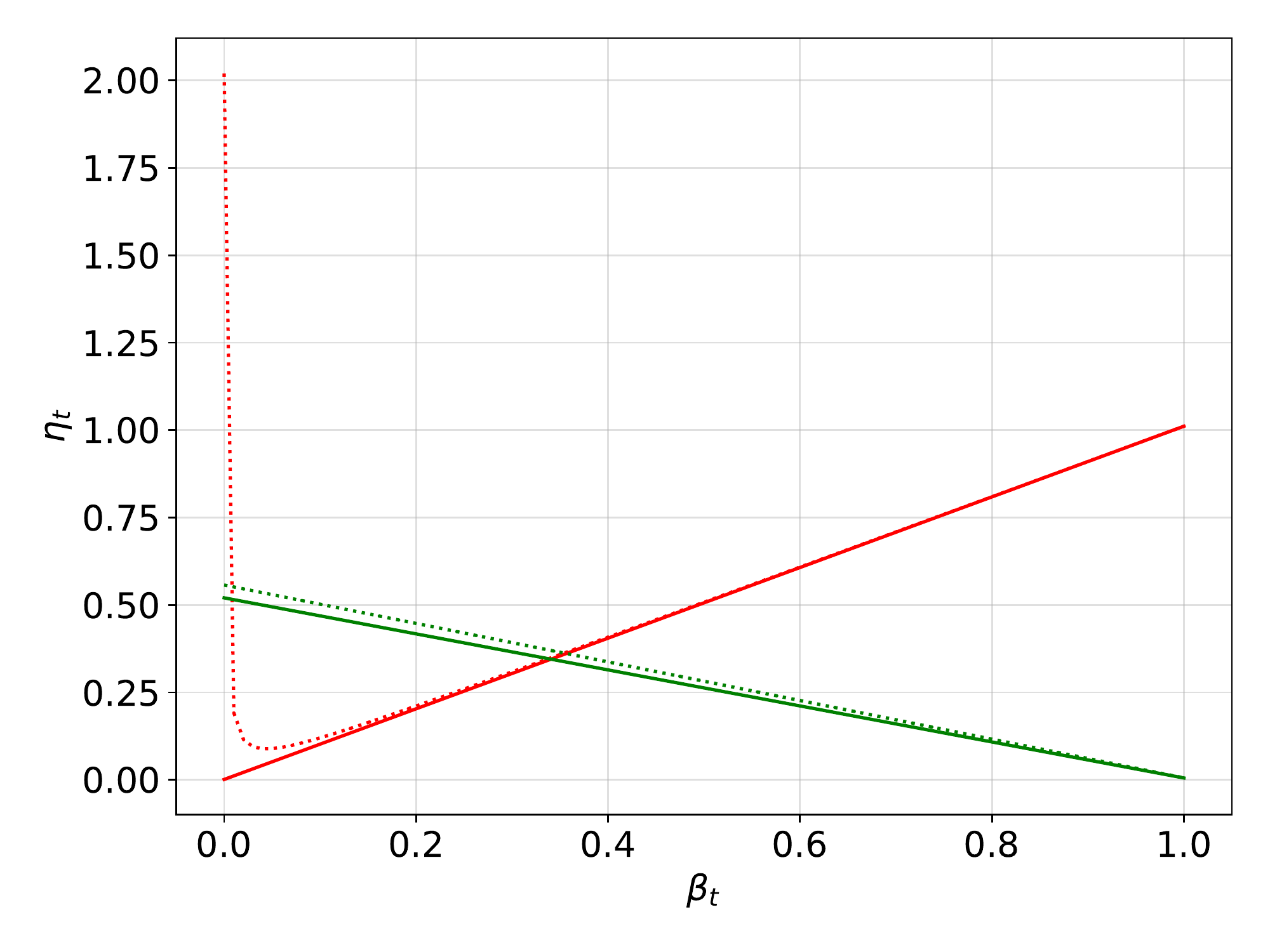}}
     \vskip-18pt
    \caption{Intuition for the $\beta_t$ choice.}
    \vskip-10pt
    \label{beta_intuition}
\end{wrapfigure}

In this section, we continue the convergence analysis started in ``One-step effects'' and based on descent lemma, using the sketch of the proof from ``Cumulative effects''. To simplify the reasoning, we consider \texttt{Scaled SGD} without variance reduction but updating step direction as $g_t = \nabla f(x_t)$ with probability $1-p$.

The following corollary of Theorem~\ref{descent_lemma} shows the convergence of \texttt{Scaled SGD} updating step direction as $g_t = \nabla f(x_t)$ with probability $1-p$.
\begin{theorem}\label{conv_onestep}If $\eta_t \leqslant \min\big\{\tfrac{1 - \beta_t \chi_t}{(1 + \sigma) (\Phi + \sqrt{M^\prime/ 6 \cdot \|g_t\|_{P_t}^*} )}, \tfrac{1/(1 + \kappa_t) + \beta_t}{1-p}\big\}$, sequence of the points generated by \texttt{Scaled SGD} satisfies
$
    \min_{t=1,...,T} \|\nabla f(x_t)\|^{*\textsuperscript{2}}_{P_t} = O (\tfrac{f(x_0) - f(x_*)}{T} ).
$
\end{theorem}
Note that the first term in $\eta_t \leqslant \min\{\cdot,\cdot\}$ is decreasing, and the second one is increasing, so there is $\beta_t$ at which upper bound on $\eta_t$ has a fracture and starting from which it is determined by decreasing second term (see Figure~\ref{beta_intuition}). Thus, upper bound on $\eta_t$ attains its maximum with respect to $\beta_t$ at this point. So, we have equation $\frac{1 - \beta_t \chi_t}{(1 + \sigma) (\Phi + \sqrt{M^\prime/ 6 \cdot \|g_t\|_{P_t}^*} )} = \frac{1/(1+\kappa_t) + \beta_t}{1-p}$,
which leads to
\begin{equation}\label{beta_formula}
    {\small
    \beta_t = \max \Big\{\tfrac{\beta_{t-1}}{2}, \tfrac{(1-p)(1 + \kappa_t + \chi_t)}{(1+\sigma) (\Phi + \sqrt{ M^\prime/6 \cdot\|g_t\|_{P_t}^*} ) + (1-p)\chi_t} - 1 \Big\}.
    }
\end{equation}
\vspace{-20pt}
    
\section{Observed effects}
Let's consider the binary logistic regression task $\min_{x \in \mathbb{R}^n}  \{ f(x) := \frac{1}{m} \sum_{i=1}^m  \log (1+e^{-b_i (a\;\circ\;a_i)^\top x} )  \}$,
where $\{(a_i, b_i)\}_{i=1}^m$ is a dataset containing features $a_i$ and classes $b_i \in \{-1, 1\}$, and $a \in \mathbb{R}^n$ is vector of random i.i.d. scaling factors drawn from $\mathcal{U}[-A, A]$ and corresponding to each feature, $\circ$ denotes Hadamard product. We use \texttt{a9a} dataset with $m = 32561$, $n = 123$. For logistic regression problem in this formulation, $\nabla f$ is Lipschitz continuous with constant $L = \frac{1}{4} \|(a \circ a_1\;...\;a \circ a_m)\|_2 = O(A)$, so that $L$ is proportional to $A$, which is helpful for the design of experiments. In all the experiments we froze the following hyperparameters: batch size $= 100$, $p = 0.9$.
\begin{figure}[ht!]
\centering
    \vskip-6pt
     \includegraphics[width=0.24\columnwidth]{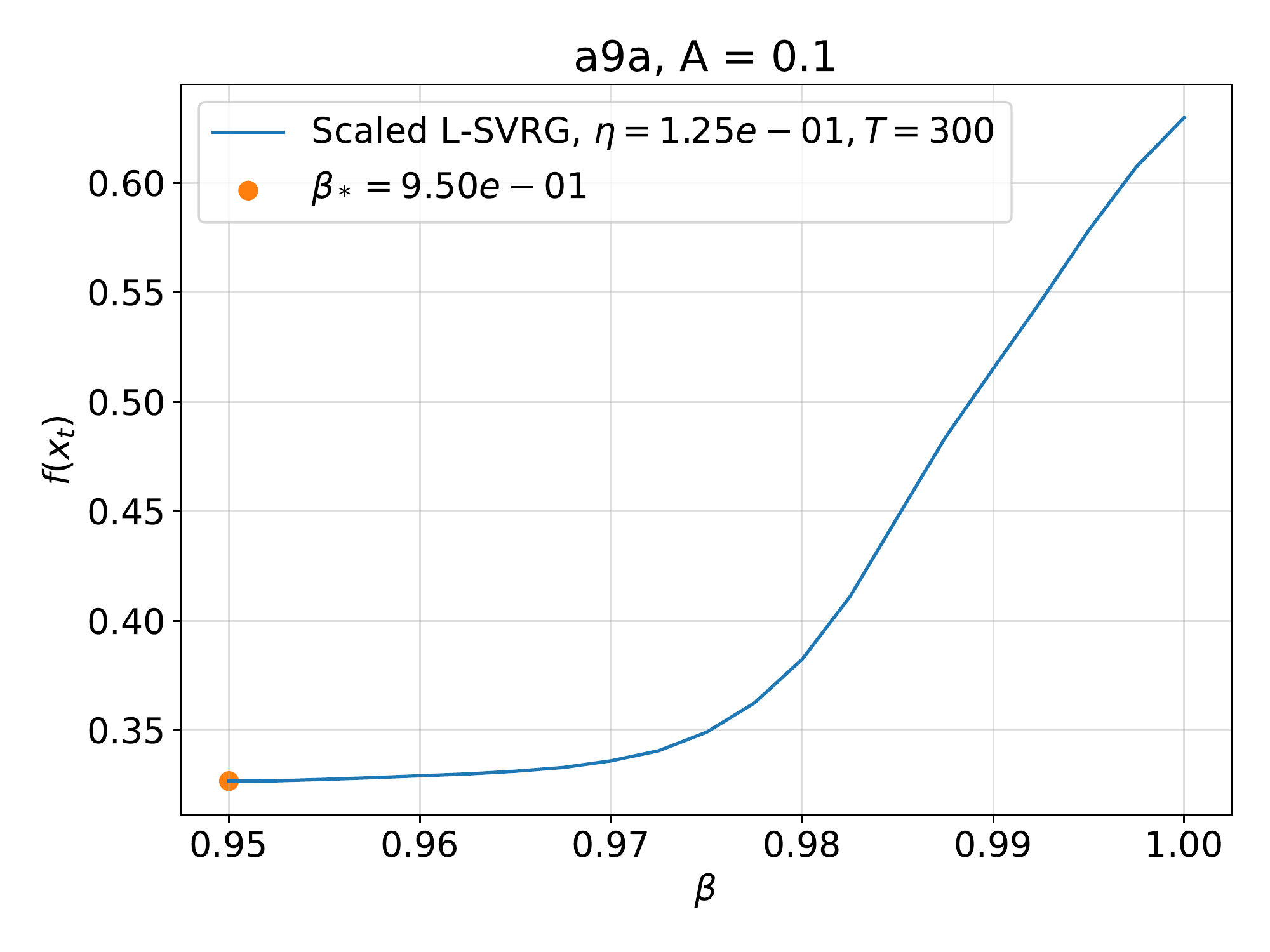} 
     \includegraphics[width=0.24\columnwidth]{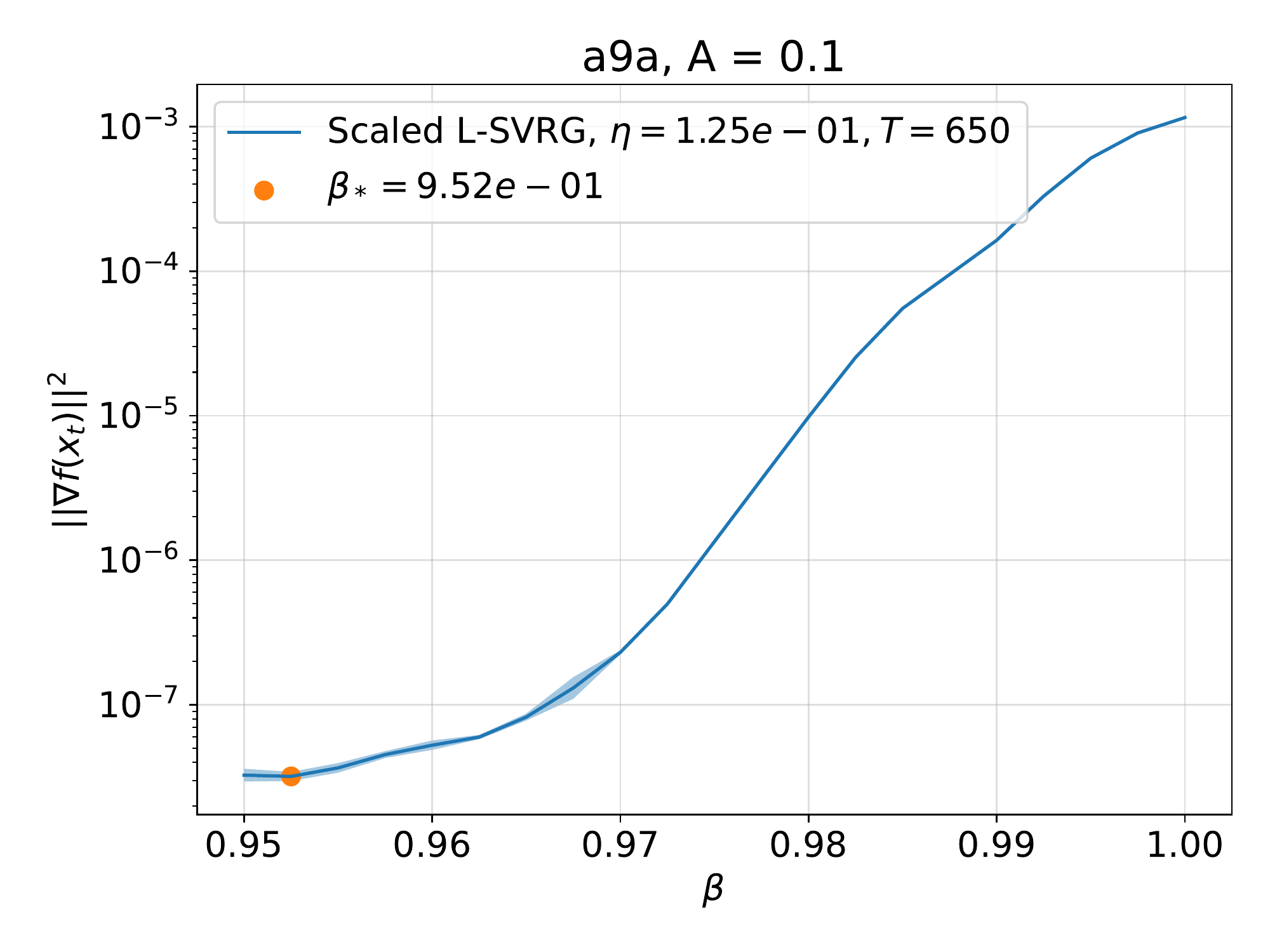} 
 \includegraphics[width=0.24\columnwidth]{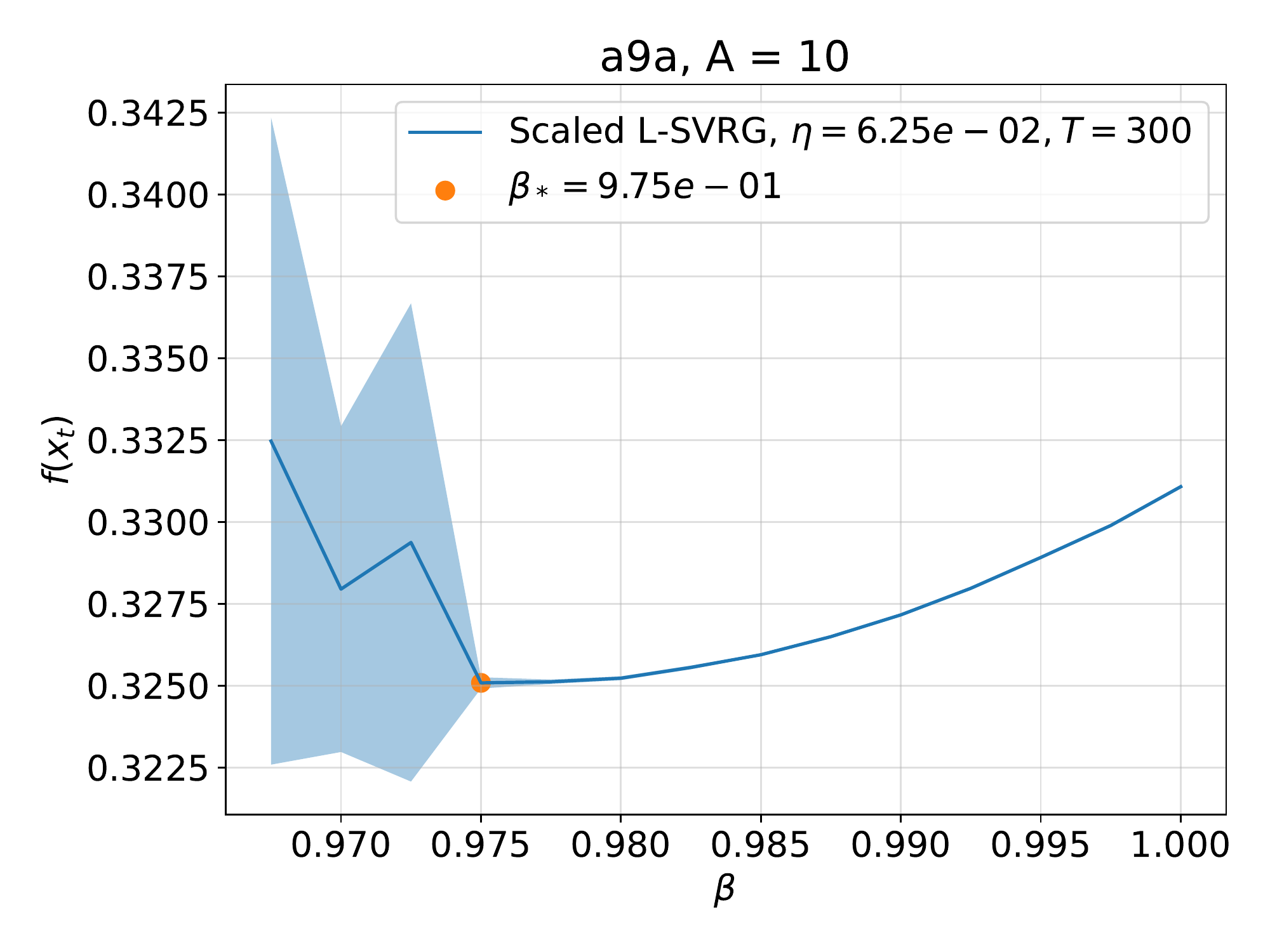}
     \includegraphics[width=0.24\columnwidth]{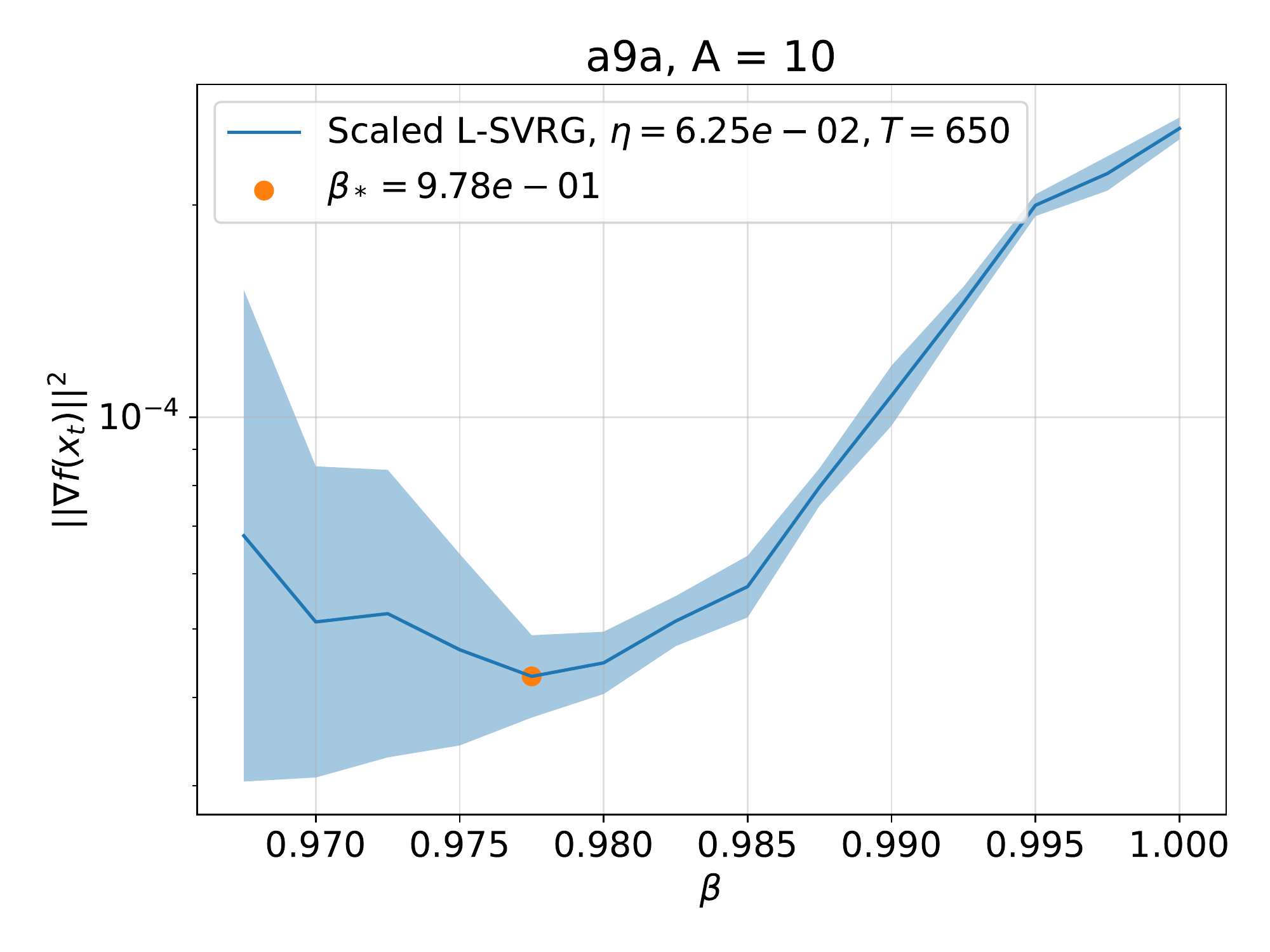}
     \vskip-12pt
    \caption{Dependence of achieved precision on $\beta_t \equiv \beta$.}
    \vskip-6pt
    \label{betas}
\end{figure}

Figure~\ref{betas} summarizes the results of \texttt{Scaled L-SVRG} runs with different $\beta$'s. Horizontal axis measures value of $\beta$, vertical axis measures objective function value $f(x_T)$ (or $\|\nabla f(x_T)\|^2$) after $T = 300$ iterations of the algorithm. Firstly, scaling give significant benefit in comparison with not scaled method, corresponding to $\beta = 1$. (More detailed experiments are presented in Appendices~\ref{sup_exp} and \ref{om_fig}). It can be seen that dependence of achieved precision on $\beta$ changes with increasing of $A$: minimum of the corresponding function is getting closer to $\beta = 1$, its values on the left from minimum are growing and its growth rate near $\beta = 1$ is significantly increasing. This relationship between $\beta$ and $A$ reflects the trade-off between variance compensation and scaling. Variance affects the convergence if $\beta$ is small: increasing of $L$ leads to the increasing of $\delta^+_t$, which increases the accumulating error term in Proposition \ref{variance}. To explain the behaviour near $\beta = 1$, consider the $A = 0.1$ case, where variance error terms are insignificant. Values begin to grow rapidly starting from $\beta \approx 0.97$ and stop on some fixed value at $\beta = 1$. This behaviour is described in Proposition \ref{Delta}, where we have shown the $O(1/(1 - \beta))$ growth of gradients Lipschitz constant. The boundedness at $\beta = 1$ can be explained by the proper choice of $P_0$, such that $\delta^-_t \neq 1$, even if $\beta = 1$. Thus, the main outline of our theory is successfully confirmed on the experiment.

\begin{wrapfigure}{r}{0.53\columnwidth}
    \centering
    \vskip-12pt
     {\includegraphics[width=0.26\columnwidth]{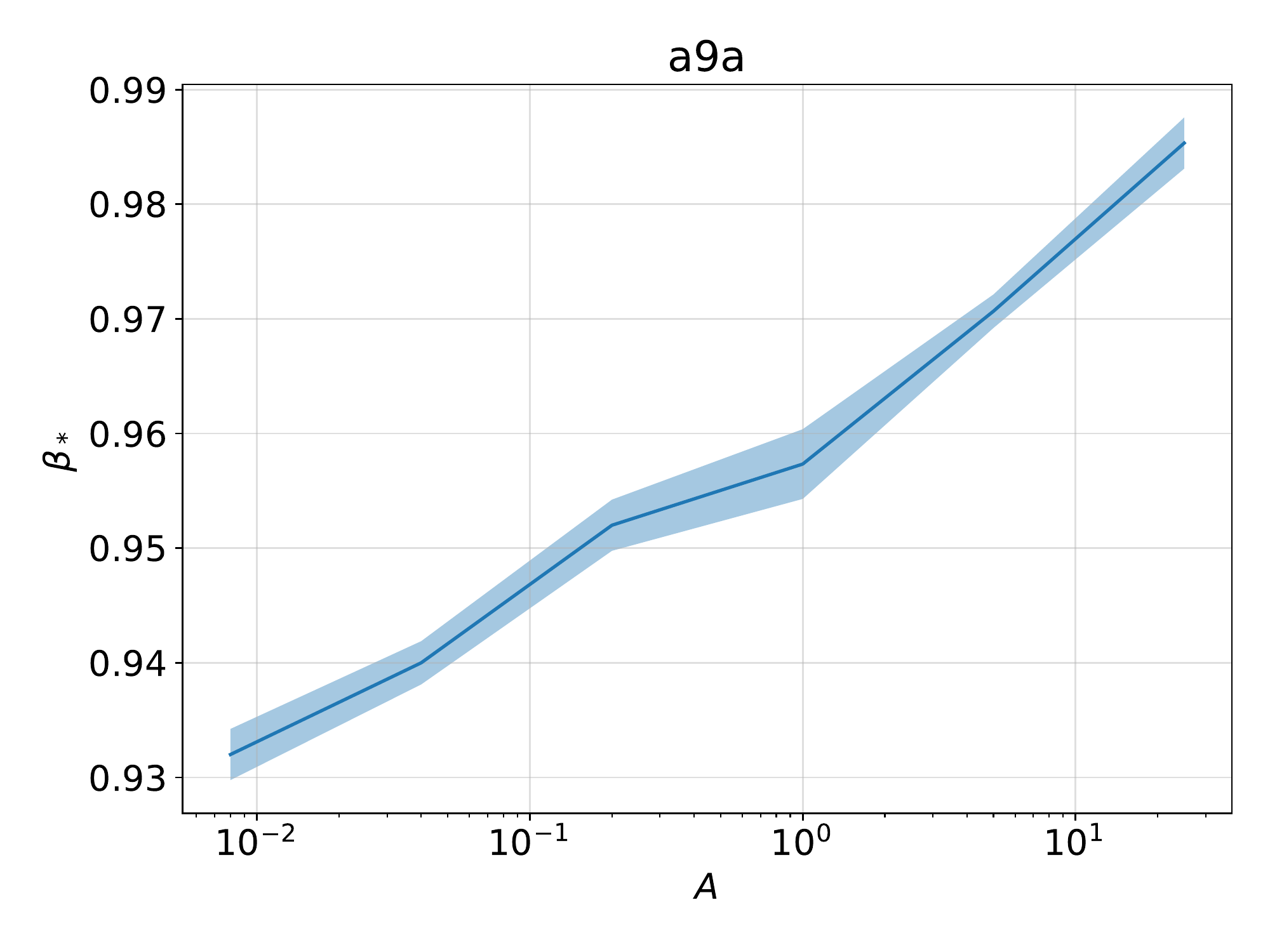}}
     {\includegraphics[width=0.26\columnwidth]{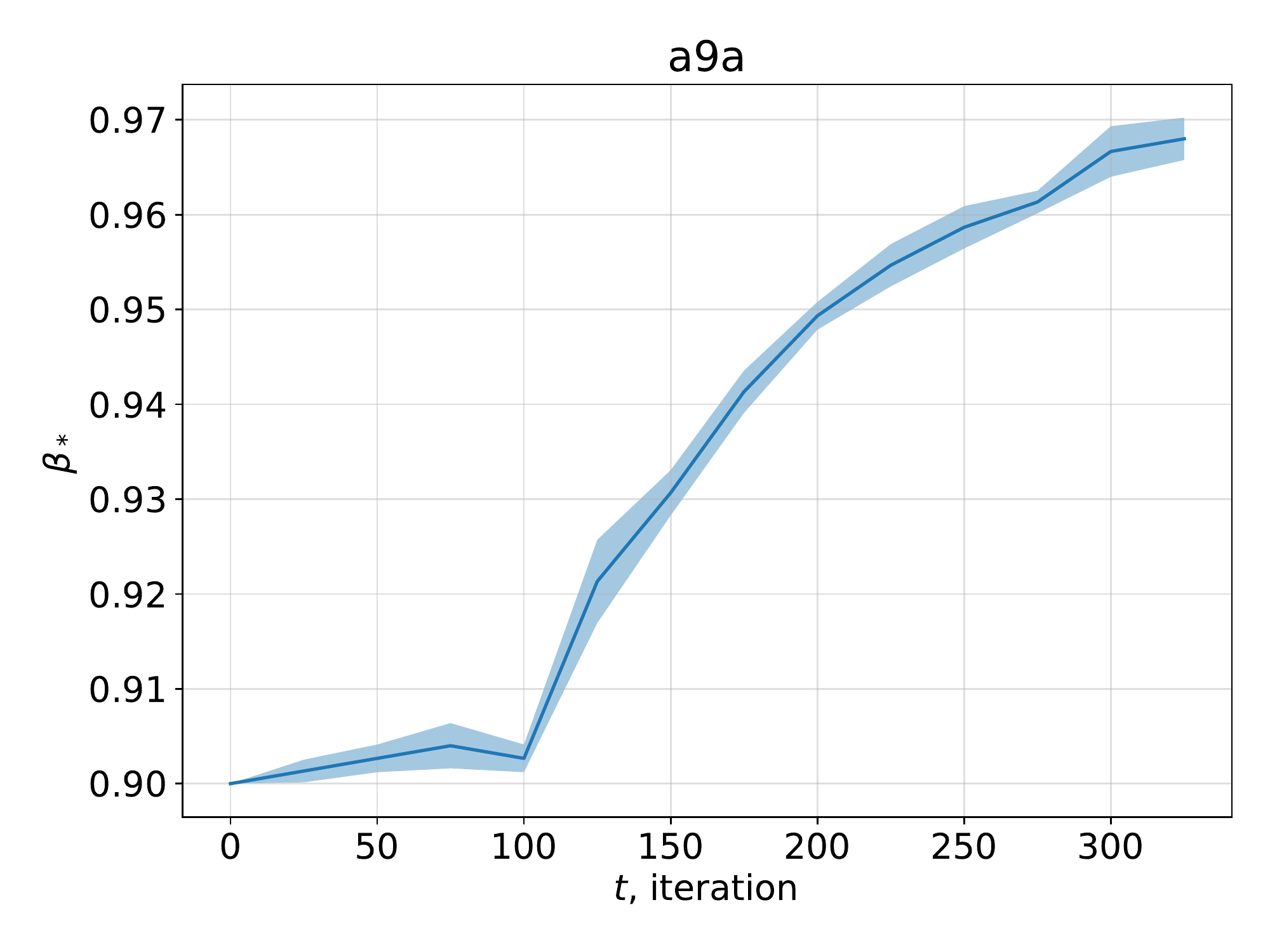}}
     \vskip-12pt
    \caption{Dependencies of optimal $\beta_*$.}
    \vskip-6pt
    \label{optimal_beta}
\end{wrapfigure}
The dependence of optimal $\beta_*$ on smoothness characteristic $A$ is presented on the Figure~\ref{optimal_beta}a. It can be seen that $\beta_*$ grow slowly with increasing of $A$ (plot is close to linear in logarithmic scale for $A$). This means that there is no need for $\beta$ to be in proportional dependence with $\eta$ or $L$. On the other hand, best choice of $\beta$ is close to standard $\beta \approx 0.99$ for large enough $L$ (and it does not matter, $L = 100$ or $L = 200$, because they are of the same order), while small $L$ make choice of $\beta$ very sensitive (that is, one need to make a distinctions between $L = 0.01$ and $L = 0.001$, although they are pretty close to each other).


In addition to dependence of optimal $\beta_*$ on smoothness of the problem, we test the dependence on the number of iterations. This experiment gives a rough estimate for the best choice of $\beta_t$. On the Figure~\ref{optimal_beta}b, one can see the dependence of quasi-optimal (in a sense described above) $\beta_*$ on the number of iterations. Optimal $\beta_*$ is getting closer to $\beta=1$ with increasing number of iterations, while at the beginning of methods operation the best value is significantly lower. The latter can be explained by the need to rapidly adapt $P_t$ to some good estimation of component-wise scaling from its initial value $P_0 = I$. On the other hand, the convergence of optimal $\beta_*$ to $1$ as the iteration number tends to infinity can be explained in view of remark from the end of ``Cumulative effects'' section.

\newpage
\bibliography{main}

\appendix
\section{Additional theory}
\begin{remark}
It turns out that formula in Theorem~\ref{descent_lemma} includes both $\beta_t$ and $\beta_{t+1}$. To find optimal $\beta_t$ we should consider sum of descents on two corresponding steps, but it could make the analysis more complicated. One can assume $\beta_1 = \beta_2 = ... \equiv \beta$ to find optimal $\beta$. But it is easy to see that 
\[
\arg \min_{\beta} \frac{(1 - \beta) \beta}{1/\kappa_t+\beta^2} = 0,\quad \arg \min_{\beta} \frac{1+\sigma}{1 - \beta_t \chi_t} = 0,
\]
so we lose adaptivity. This defect is common for $\delta_t$ linearly dependent on $\beta$.

In practice, $\beta$ is usually chosen to be close to $1$, assuming that dependence of $L_t$ on $\beta$ is ``weaker'' than dependence of penalty for variance on $\beta$. So, it worth estimate the worst case (with respect to the $\beta$) multiplier appearing in variance term. This can be done explicitly
\[
    1 + \max_{\beta} \frac{(1 - \beta) \beta}{1/\kappa_t+\beta^2} = 1 + \frac{1}{2} \left(\sqrt{\kappa_t+1}-1\right) = O\left( \sqrt{\frac{L}{\mu}}\right).
\]
On the other hand, the worst case smoothness multiplier is $L$ and is not being accumulated with iterations of the algorithm. Therefore, it was right to assume that $\beta$ should be close to $1$ giving priority to restrain the growth of variance term.
\end{remark}

\begin{remark}
Note, that second term in \eqref{beta_formula} is less than zero while $\|g_t\|_{P_t}^* \geqslant \displaystyle\frac{6}{M^\prime} \left(\frac{(1-p)(1 + \kappa_t)}{1+\sigma} - \Phi\right)^2$, so $\beta_t$ should be set to zero on the first iterations (the more ill-conditioned is the function, the more iterations are needed). On the other hand, $\|g_t\|_{P_t}^* \rightarrow 0$, so $\beta_t$ on the later iteration is determined by $\displaystyle\frac{(1-p)(1 + \kappa_t) - (1+\sigma)\Phi}{(1-p)\chi_t + (1+\sigma)\Phi}$, which is smaller than one if $\displaystyle1 + \kappa_t - \chi_t < \frac{2 \Phi(1+\sigma)}{1-p}$ (the more well-conditioned is the function, the less $\beta_t$ is). In addition, note that $\beta_t$ is increasing function of $\kappa_t$ and decreasing function of $\chi_t$.
\end{remark}

Let's now establish the starting acceleration of the convergence, achieved by scaled methods. For this purpose we change \eqref{x_bound} as
\[
x \geqslant \frac{-1 + \sqrt{3 + 2A\|g_t\|_{P_t}^*}}{2(1+A\|g_t\|_{P_t}^*)},
\]
which is solution to $(1+A\|g_t\|_{P_t}^*)x^2 + x - 1/2 = 0$. Then, following the analysis proposed in \cite{hanzely2021damped}, we state
\[
\eta_t \geqslant \frac{1}{4\max\left\{1,\sqrt{A \|g_t\|_{P_t}^*}\right\}},
\]
which leads to
\[
\eta_t \|g_t\|_{P_t}^* \geqslant \left\{ \begin{array}{ll}
\displaystyle\frac{\|g_t\|_{P_t}^{*\textsuperscript{2}}}{4},&\text{if }\|g_t\|_{P_t}^* < 1/A\\
\displaystyle\frac{\|g_t\|_{P_t}^{*\textsuperscript{3/2}}}{4 \sqrt{A}},&\text{otherwise}.
\end{array} \right.
\]
Thus, we have two cases of convergence lemma for Lyapunov function $V_t$ introduced above:
\[
\left\{ \begin{array}{ll}
\displaystyle \E{V_{t+1}} \leqslant V_t - \frac{1}{16} (2 \|\nabla f(x_t)\|_{P_t}^{*\textsuperscript{2}} + \|g_t\|_{P_t}^{*\textsuperscript{2}}),&\text{if }\|g_t\|_{P_t}^* < 1/A\\
\displaystyle\E{V_{t+1}} \leqslant V_t - \frac{1}{16\sqrt{A}} \left(2 \frac{\|\nabla f(x_t)\|_{P_t}^{*\textsuperscript{2}}}{\|g_t\|_{P_t}^{*\textsuperscript{1/2}}} + \|g_t\|_{P_t}^{*\textsuperscript{3/2}}\right),&\text{otherwise}.
\end{array} \right.
\]
Thus, at the starting iterations, when $\displaystyle\|g_t\|_{P_t}^* \geqslant \frac{3}{M^\prime} \sqrt{\frac{1+\sigma}{1 - \beta_t \chi_t}}$, convergence rate of scaled method is $\displaystyle O\left(\frac{1}{T^{3/2}}\right)$. Moreover, the greater is $\beta_t$, the shorter this starting acceleration lasts. But $\displaystyle\min_{\beta_t} A = \frac{M^\prime}{3\sqrt{1+\sigma}}$, so there cannot be an acceleration on the latter iterations. Anyway, $\beta_t$ could be chosen in a way to lengthen this starting convergence.
\section{Omitted proofs}
\subsection{Proof of Proposition~\ref{Delta}}

Then, inexactness $\Delta_{t+1}$ can be estimated as follows
\[
    \nabla^2 f(x_{t+1}) P_{t+1}^{-1} = \frac{1}{1-\beta_{t+1}} \nabla^2 f(x_{t+1}) d_{t+1}^{-1} - \frac{1}{1-\beta_{t+1}} \nabla^2 f(x_{t+1}) d_{t+1}^{-1} \left[\frac{1-\beta_{t+1}}{\beta_{t+1}} d_{t+1} P_t^{-1} + I\right]^{-1},
\]
using Woodbury identity
\[
    (A+B)^{-1} = A^{-1} - A^{-1} (A B^{-1} + I)^{-1},
\]
that implies
\[
    \Delta_t \leqslant \frac{1+\sigma}{1-\beta_t} - 1,
\]
which is increasing with $\beta$, always greater than $\sigma$, increases linearly near $\beta = 0$ and inversely proportional near $\beta = 1$.

Estimation for $\Delta$ is rough near $\beta = 1$. Its because we neglect second term in Woodbury identity. To work it out, we use the relation between $P_t$ and $d_{t+1}$ to bound the second term in Woodbury identity, using the formula for $P_{t+1}$:
\begin{align*}
    \frac{1}{1 - \beta_{t+1}} P_{t+1} d_{t+1}^{-1} = \frac{\beta_{t+1}}{1 - \beta_{t+1}} P_t d_{t+1}^{-1} + I,\\
    \frac{\beta_{t+1}}{1 - \beta_{t+1}} P_t d_{t+1}^{-1} = \frac{1}{1 - \beta_{t+1}} \left(P_{t+1} d_{t+1}^{-1} - I + \beta_{t+1} I\right) &\succcurlyeq \frac{\beta_{t+1} - \delta_{t+1}^-}{1 - \beta_{t+1}} I,\\
    \frac{1 - \beta_{t+1}}{\beta_{t+1}} d_{t+1} P_t^{-1} + I &\preccurlyeq \frac{1 - \delta_{t+1}^-}{\beta_{t+1} - \delta_{t+1}^-} I,\\
    \left[\frac{1 - \beta_{t+1}}{\beta_{t+1}} d_{t+1} P_t^{-1} + I\right]^{-1} &\succcurlyeq \frac{\beta_{t+1} - \delta_{t+1}^-}{1 - \delta_{t+1}^-} I.
\end{align*}
Substituting this bound in Woodbury identity finishes the proof.

\subsection{Proof of Proposition~\ref{variance}}

Let's apply Woodbury identity to $P_{t+1}$:
\[
    P_{t+1}^{-1} = \frac{1}{\beta_{t+1}} P_t^{-1} - \frac{1}{\beta_{t+1}} P_t^{-1} \left[\frac{\beta_{t+1}}{1-\beta_{t+1}}P_t d_t^{-1} + I\right]^{-1}.
\]
Then,
\[
    \langle s, P_{t+1}^{-1} s \rangle = \frac{1}{\beta_{t+1}} \langle s, P_t^{-1} s \rangle - \frac{1}{\beta_{t+1}} \langle s, \left(\left[\frac{\beta_{t+1}}{1-\beta_{t+1}}P_t d_t^{-1} + I\right] P_t\right)^{-1} s \rangle.
\]
On the other hand,
\[
    \frac{\beta_{t+1}}{1-\beta_{t+1}}P_t d_t^{-1} + I \preccurlyeq \left(\frac{\beta_{t+1}}{1-\beta_{t+1}} (1+\delta_t^+) + 1\right) I \preccurlyeq \frac{1+\beta_{t+1}\delta_t^+}{1-\beta_{t+1}} I.
\]
Finally, we substitute this bound in equality on variation and obtain that
\begin{equation}
    \|s\|_{P_{t+1}}^{*\textsuperscript{2}} \leqslant \frac{1+\delta_t^+}{1+\beta_{t+1}\delta_t^+} \|s\|_{P_t}^{*\textsuperscript{2}} = \left(1 + \frac{1 - \beta_{t+1}}{1/\delta_t^+ +\beta_{t+1}}\right) \|s\|_{P_t}^{*\textsuperscript{2}}.
\end{equation}

\subsection{Proof of Proposition~\ref{delta_convex}}

It follows from
\begin{gather*}
    P_{t+1} d_{t+1}^{-1} = (1 - \beta_{t+1}) I + \beta_{t+1} P_{t} d_{t+1}^{-1} \preccurlyeq \left[1 + \beta_{t+1} \left(\frac{\max_i {[P_{t}]_{ii}}}{\min_i {[d_{t+1}]_{ii}}} - 1\right)\right] I.
\end{gather*}

\subsection{Proof of Proposition~\ref{delta_concordance}}

We need in the following corollary of strong self-concordance:
\begin{lemma}[Rodomanov--Nesterov \cite{rodomanov2021greedy}]
For all $x, y \in \mathbb{R}^n$, it holds that
\begin{equation*}
    \frac{\text{diag}\left(\nabla^2 f(x)\right)}{1 + N \|y - x\|_{\text{diag}\left(\nabla^2 f(x)\right)}} \preccurlyeq \text{diag}\left(\nabla^2 f(y)\right) \preccurlyeq \left(1 + N \|y - x\|_{\text{diag}\left(\nabla^2 f(x)\right)}\right) \text{diag}\left(\nabla^2 f(x)\right),
\end{equation*}
where $A \preccurlyeq B$ means $\langle (B - A) x, x \rangle \geqslant 0$ for all $x \in \mathbb{E}$.
\end{lemma}
Presented lemma implies the bound on $\delta$ which takes into account that $P_t d_{t+1}^{-1}$ term is as small as the step on iteration $t$. We have
\begin{gather*}
    \left(1 + N \eta_t \|g_t\|_{P_t}^* \sqrt{1 + \delta_t^+}\right)^{-1} d_t \preccurlyeq d_{t+1} \preccurlyeq \left(1 + N\|x_{t+1} - x_t\|_{d_t}\right) d_{t} \preccurlyeq \left(1 + N \eta_t \|g_t\|_{P_t}^* \sqrt{1 + \delta_t^+}\right) d_t,\\
    P_{t+1} d_{t+1}^{-1} = (1 - \beta_{t+1}) I + \beta_{t+1} P_{t} d_{t+1}^{-1} \preccurlyeq \left[1 + \beta_{t+1} (\delta_t^+ + \delta_t^+ \sqrt{1 + \delta_t^+} N \eta_t \|g_t\|_{P_t}^* - 1)\right] I,
\end{gather*}
which implies the statement of the proposition.

\subsection{Proof of Theorem~\ref{descent_lemma}}

\begin{gather*}
    f(x_{t+1}) \leqslant f(x_t) - \eta_t \langle \nabla f(x_t), P_t^{-1} g_t \rangle + \frac{\eta_t^2}{2} \frac{1+\sigma}{1 - \beta_t \chi_t} {\|g_t\|_{P_t}^{*\textsuperscript{2}}} + \frac{M^\prime \eta_t^3}{6} \left(\frac{1+\sigma}{1 - \beta_t \chi_t}\right)^{3/2} \|g_t\|_{P_t}^{*\textsuperscript{3}} \\
    -\eta_t \langle \nabla f(x_t), P_t^{-1} g_t \rangle = - \frac{\eta_t}{2} \|\nabla f(x_t)\|_{P_t}^{*\textsuperscript{2}} - \frac{\eta_t}{2} \|g_t\|_{P_t}^{*\textsuperscript{2}} + \frac{\eta_t}{2} \|g_t - \nabla f(x_t)\|_{P_t}^{*\textsuperscript{2}}\\
    \E{\|g_{t+1} - \nabla f(x_{t+1})\|_{P_{t+1}}^{*\textsuperscript{2}}} \leqslant  \left(1 + \frac{(1 - \beta_{t+1}) \beta_t}{1/\kappa_t+\beta_t \beta_{t+1}}\right) \|g_t - \nabla f(x_t)\|_{P_t}^{*\textsuperscript{2}} + \E{\|\nabla f(x_{t+1}) - \nabla f(x_t)\|_{P_t}^{*\textsuperscript{2}}}\\
    \|\nabla f(x_{t+1}) - \nabla f(x_t)\|_{P_t}^{*\textsuperscript{2}} \leqslant \eta_t^2 \left(\frac{1+\sigma}{1 - \beta_t \chi_t}\right)^2 \|g_t\|_{P_t}^{*\textsuperscript{2}}
\end{gather*}

\subsection{Proof of Theorem~\ref{convergence}}

We will demonstrate the convergence of the method for such a Lyapunov function $V_t = f(x_t) - f(x_*) + a_t \|x_t - y_t\|^2_{P_t}$. Due to $L_t$-Lipschitz smoothness of $f$, we have
\begin{align*}
    \E{V_{t+1}} &\leqslant f(x_t) - f(x_*) - \eta_t \langle \nabla f(x_t), P_t^{-1} \nabla f(x_t) \rangle + \frac{\eta_t^2 L_t}{2} \E{\|g_t\|^{*\textsuperscript{2}}_{P_t}} + a_{t+1} \E{\|x_{t+1} - y_{t+1}\|^2_{P_{t+1}}}.
\end{align*}
It can be easily proven that
\begin{align*}
    \E{\|g_t\|^{*\textsuperscript{2}}_{P_t}} \leqslant 3 \|\nabla f(x_t)\|^{*\textsuperscript{2}}_{P_t} + 6 L_t^2 \|x_t - y_t\|_{P_t}^2.
\end{align*}
We also need in the following lemma describing the properties of Hutchinson diagonal approximation
\begin{lemma}[Jahani et al. \cite{jahani2021doubly}]
For $L$-Lipschitz smooth function $f$, it holds that
\begin{enumerate}
    \item $\left|\left[z_t \circ \nabla^2 f(x_{t+1}) z_t\right]_i\right| \leqslant \Gamma \leqslant \sqrt{n} L$.
    \item $\exists \delta \leqslant 2(1 - \beta)\Gamma$ such that $\forall t: \|P_{t+1} - P_t\|_\infty \leqslant \delta$.
\end{enumerate}
\end{lemma}
\noindent Then, we can bound last term of $V_{t+1}$ as follows
\begin{align*}
    \E{\|x_{t+1} - y_{t+1}\|_{P_{t+1}}^2} &\leqslant p \eta_t^2 \E{\|g_t\|^{*\textsuperscript{2}}_{P_t}} + (1-p)(1+\eta_t b_t) \|x_t - y_t\|_{P_t}^2\\
    &\quad+ (1-p)\frac{\eta_t}{b_t} \|\nabla f(x_t)\|^{*\textsuperscript{2}}_{P_t} + 2 \Gamma \E{(1-\beta_{t+1}) \|x_{t+1} - y_{t+1}\|_2^2}
\end{align*}
due to Fenchel--Young inequality $\langle \nabla f(x_t), x_t - y_t\rangle \leq \frac{1}{b_t} \|\nabla f(x_t)\|^{*\textsuperscript{2}} + b_t \|x_t - y_t\|^2$ for some sequence $b_t > 0$, $t=1,2,...$ we specify later. Now,
\begin{align*}
    \E{V_{t+1}}&\leqslant f(x_t) - f(x_*) - \eta_t \left(1 - (1-p)\frac{a_{t+1}}{b_t}\right) \|\nabla f(x_t)\|^{*\textsuperscript{2}}_{P_t} + a_{t+1} (1-p)(1+\eta_t b_t)\|x_t - y_t\|^2_{P_t}\\
    &\quad+\eta_t^2 \left(\frac{L_t}{2} + a_{t+1}\right) \left(3 \|\nabla f(x_t)\|^{*\textsuperscript{2}}_{P_t} + 6 L_t^2 \|x_t-y_t\|^2_{P_t}\right) + 2 a_{t+1} \Gamma \E{(1-\beta_{t+1}) \|x_{t+1} - y_{t+1}\|_2^2}\\
    &\leqslant f(x_t) - f(x_*) - \eta_t \left(1 - (1-p)\frac{a_{t+1}}{b_t} - 3 a_{t+1} \eta_t - 3L_t \eta_t\right) \|\nabla f(x_t)\|^{*\textsuperscript{2}}_{P_t}\\
    &\quad + a_{t+1}\left((1-p)(1+\eta_t b_t) + 3 \eta^2 \left(\frac{L_t}{a_{t+1}} + 2\right) L_t^2\right) \|x_t - y_t\|^2_{P_t}\\
    &\quad+ 2 a_{t+1} \Gamma \E{(1-\beta_{t+1}) \|x_{t+1} - y_{t+1}\|_2^2}.
\end{align*}

Finally, we determine the step size $\eta_t$ satisfying
\begin{equation*}
\left\{ \begin{array}{ll}
\displaystyle 1 - (1-p)\frac{a_{t+1}}{b_t} - 3a_{t+1}\eta_t - 3 L_t \eta_t \geqslant \frac{1}{4},\\
\displaystyle (1-p)(1+\eta_t b_t) + 3\eta_t^2 \left(\frac{L_t}{a_{t+1}} + 2\right) L_t^2 \leqslant \frac{a_t}{a_{t+1}}.
\end{array} \right.
\end{equation*}
We can set $a_{t+1} = L_t, b_t = p/\eta_t, \eta_t = c_t/L_t$ (this is only a comfortable option, but it is not unique; one can try to find an optimal one) for every $t=1,2,...$ and variable sequence $c_t$, and after substitution we solve it with respect to $c_t$ to obtain $\eta_t \leqslant \min\left\{\frac{\alpha p}{3}, \frac{3}{4}\frac{p}{5p+1}\right\} \frac{1}{L_t}$ for some $\alpha>0$ as a sufficient condition\footnote{To be pedantic, $\eta_t \leqslant \min\left\{\frac{\sqrt{p^2 - 1 + L_{t-1} / L_t}}{3}, \frac{3}{4}\frac{p}{5p+1}\right\} \frac{1}{L_t}$, but it is necessary that $L_{t-1} / L_t \rightarrow 1$ on the one hand and $L_t$ is close to $1$ on the other hand, so we can find a proper $\alpha>0$; note, that if we replace $a_{t+1} = L_t$ with $a_{t+1} = \max_{k=1,...,t} L_k$ we do not have any problems with $a_t / a_{t+1} \leqslant 1$, but need to set $\eta \propto 1/L$ that we do not want to do~--- this is why $L_t$ is close to $1$ is important condition.}. Since the system of inequalities above holds now, we have
\[
\E{V_{t+1}} \leqslant V_t - \frac{\eta_t}{4} \|\nabla f(x_t)\|^{*\textsuperscript{2}}_{P_t} + 2 L_{t-1} \Gamma (1-\beta_{t+1}) \|x_{t+1} - y_{t+1}\|_2^2,
\]
that is summed up for $t=1,...,T$ to
\[
\E{\|\nabla f(\overline{x}_T)\|^{*\textsuperscript{2}}_{P_t}} \leqslant \frac{1}{\sum_{t=1}^T \eta_t} \left[V_0 - V_* + 2 \Gamma \sum_{t=2}^{T+1} L_{t-1} \E{(1-\beta_t) \|x_t - y_t\|_2^2} \right],
\]
where $\overline{x}_T$ is such that
$\overline{x}_T = x_t$ with probability $\eta_t / \sum_{k=1}^T \eta_k$ or, in the best-iteration manner,
\[
\min_{t=1,...,T} \|\nabla f(x_t)\|^{*\textsuperscript{2}}_{P_t} \leqslant \frac{1}{\sum_{t=1}^T \eta_t} \left[V_0 - V_* + 2 \Gamma \sum_{t=2}^{T+1} L_t \E{(1-\beta_t) \|x_t - y_t\|_2^2} \right].
\]

It is time to remember that, in opposite to standard analysis, $\eta_t$ depends on $L_t$ in our case, so we have factor of the form $1/ \left(\sum_{t=1}^T 1/L_t\right)$ in convergence rate. Writing it as harmonic average finishes the reasoning.

\subsection{Proof of Theorem~\ref{conv_onestep}}

Firstly, let us rewrite the descent lemma in a way to get rid of cubic term, as follows
\begin{align*}
    \E{f(x_{t+1})} &\leqslant f(x_t) - \frac{\eta_t}{2} \|\nabla f(x_t)\|_{P_t}^{*\textsuperscript{2}} + \frac{\eta_t}{2} (1-p) \left(1 + \frac{(1 - \beta_{t+1}) \beta_t}{1/\kappa_t+\beta_t \beta_{t+1}}\right) \|g_t - \nabla f(x_t)\|_{P_t}^{*\textsuperscript{2}}\\
    &\quad+ \frac{\eta_t}{2} \left[\left(1 + \frac{M^\prime}{3} \sqrt{\frac{1 - \beta_t \chi_t}{1+\sigma}} \|g_t\|_{P_t}^*\right) \left(\frac{\eta_t (1+\sigma)}{1 - \beta_t \chi_t}\right)^2 + \frac{\eta_t (1+\sigma)}{1 - \beta_t \chi_t} - 1\right] \|g_t\|_{P_t}^{*\textsuperscript{2}}.
\end{align*}
Denoting $\displaystyle\frac{\eta_t (1+\sigma)}{1 - \beta_t \chi_t}$ as $x$ and $\displaystyle\frac{M^\prime}{3} \sqrt{\frac{1 - \beta_t \chi_t}{1+\sigma}}$ as $A$, we get the $\displaystyle\frac{\eta_t}{2} \left((1 + A \|g_t\|_{P_t}^*) x^2 + x - 1\right)$ factor for the term $\|g_t\|_{P_t}^{*\textsuperscript{2}}$, which is less than zero and can be neglected if 
\begin{equation}\label{x_bound}
    0 < x \leqslant \frac{1}{\Phi + B (1 - \beta_t \chi_t)^{1/4} \sqrt{\|g_t\|_{P_t}^*}} = \frac{2}{1+\sqrt{5} + 2\sqrt{A \|g_t\|_{P_t}^*}} \leqslant \frac{-1 + \sqrt{5 + 4 A \|g_t\|_{P_t}^*}}{2 (1 + A \|g_t\|_{P_t}^*)},
\end{equation}
where $B$ denotes the $\displaystyle\sqrt{\frac{M^\prime}{6}} (1+\sigma)^{-1/4} \leqslant \sqrt{\frac{M^\prime}{6}}$ and $\displaystyle\Phi = \frac{1 + \sqrt{5}}{2}$ is a golden ratio. This leads to the first upper bound on a step size:
\begin{equation*}
    \eta_t \leqslant \frac{1 - \beta_t \chi_t}{(1 + \sigma) \left(\Phi + B \sqrt{\|g_t\|_{P_t}^*}\right)} \leqslant \frac{1 - \beta_t \chi_t}{(1 + \sigma) \left(\Phi + B (1 - \beta_t \chi_t)^{1/4} \sqrt{\|g_t\|_{P_t}^*}\right)}.
\end{equation*}

Secondly, to establish the decrease of variance term, we consider Lyapunov function $V_t := f(x_t) - f(x_*) + \|g_t - \nabla f(x_t)\|_{P_t}^{*\textsuperscript{2}}$. Then, we need to upper bound step size once again to obtain
\begin{equation*}
    \frac{\eta_t}{2} (1-p) \left(1 + \frac{(1 - \beta_{t+1}) \beta_t}{1/\kappa_t+\beta_t \beta_{t+1}}\right) \leqslant 1.
\end{equation*}
Further, we need to distance $\beta_{t+1}$ from zero. It would be to rough to lower bound all $\beta$ with some fixed value, so we add the relation limiting the decreasing of $\beta$ instead: $\beta_{t+1} \geqslant \beta_t / 2$. This relation forces algorithm to be conservative and do not change preconditioner too much on a later iterations, that seems to be natural, because for a wide class of functions (self-concordant, for example) Hessian does not change much if step is small enough which holds for small gradients. So, upper bound on step size in our case looks like
\begin{equation*}
    \eta_t \leqslant \frac{1/(1 + \kappa_t) + \beta_t}{1-p} \leqslant \frac{2/\kappa_t+\beta_t^2}{(1-p) \left(1/\kappa_t + \beta_t\right)} \leqslant \frac{2}{(1-p) \left(1 + \displaystyle\frac{(1 - \beta_{t+1}) \beta_t}{1/\kappa_t+\beta_t \beta_{t+1}}\right)}
\end{equation*}

Thus, if $\eta_t \leqslant \min\Big\{\displaystyle\frac{1 - \beta_t \chi_t}{(1 + \sigma)\left(\Phi + \sqrt{M^\prime/ 6 \cdot \|g_t\|_{P_t}^*}\right)}, \frac{1/(1 + \kappa_t) + \beta_t}{1-p}\Big\}$, we have (now without accumulating errors!)
\begin{equation*}
    \E{V_{t+1}} \leqslant V_t - \frac{\eta_t}{4} \|\nabla f(x_t)\|^{*\textsuperscript{2}}_{P_t}\;\Longrightarrow\;\min_{t=1,...,T} \|\nabla f(x_t)\|^{*\textsuperscript{2}}_{P_t} = O\left(\frac{f(x_0) - f(x_*)}{\sum_{t=1}^T \eta_t}\right).
\end{equation*}

\newpage
\section{Supplementary numerical experiments}\label{sup_exp}
Firstly, we compare the performance of \texttt{Scaled L-SVRG} and ordinary \texttt{L-SVRG} on the problems with different smoothness characteristics $A$. For each value of $A$, we determine the best values of $\beta_t \equiv \beta$ and $\eta_t \equiv \eta$ for \texttt{Scaled L-SVRG}, and $\eta_t \equiv \eta$ for ordinary \texttt{L-SVRG} by logarithmically spaced grid search: $\eta \in \{2^{-2}, ..., 2^{-10}\}$, $\beta \in \{1-2^{-5}, ..., 1-2^{-10}\}$, while the precision achieved by tuned algorithm is estimated on average of 3 runs with random sequences of batches. 

\begin{figure}[ht!]
\centering
     \includegraphics[width=0.3\textwidth]{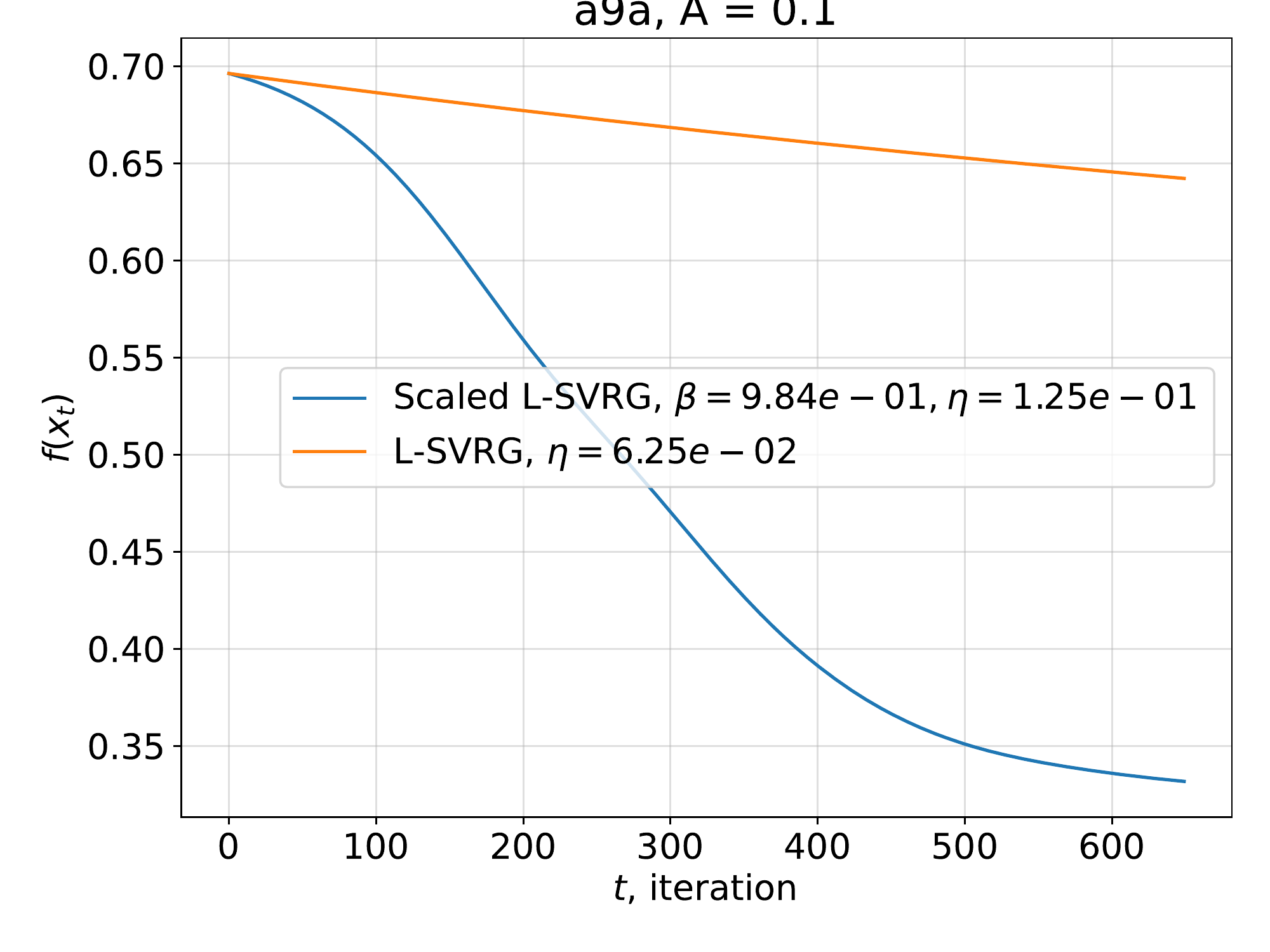}
     \includegraphics[width=0.3\textwidth]{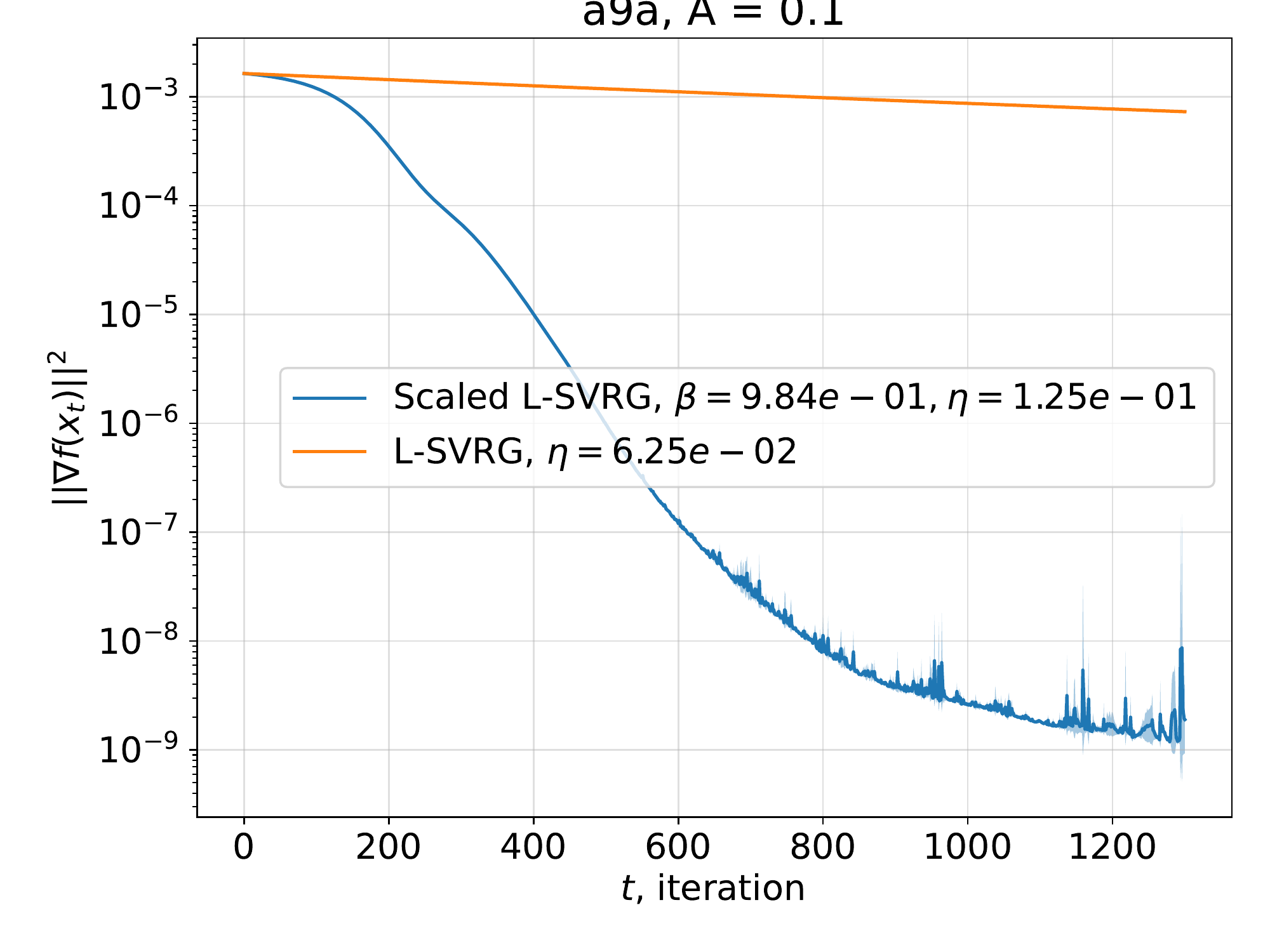}
     \includegraphics[width=0.31\textwidth]{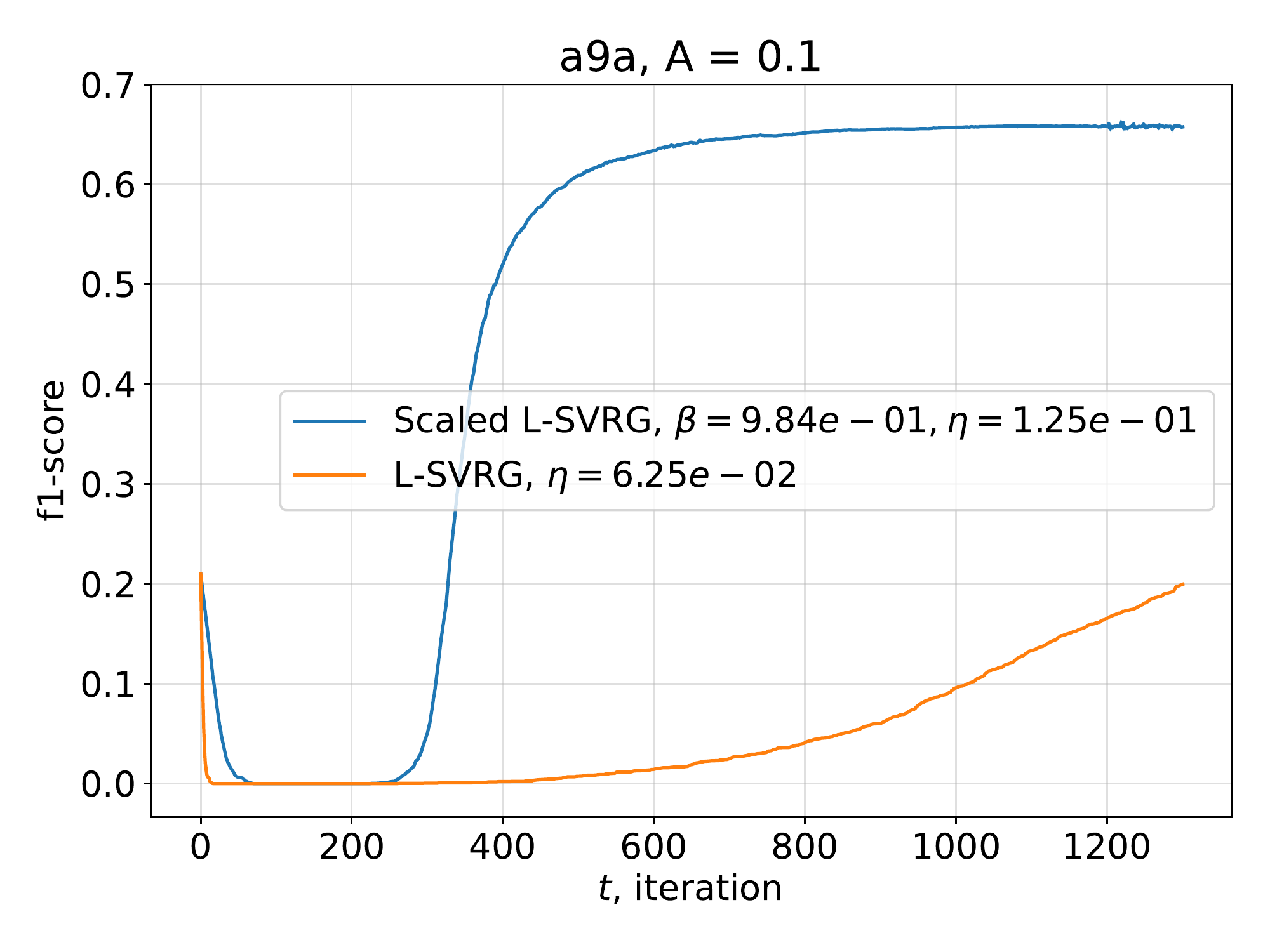}
     \includegraphics[width=0.3\textwidth]{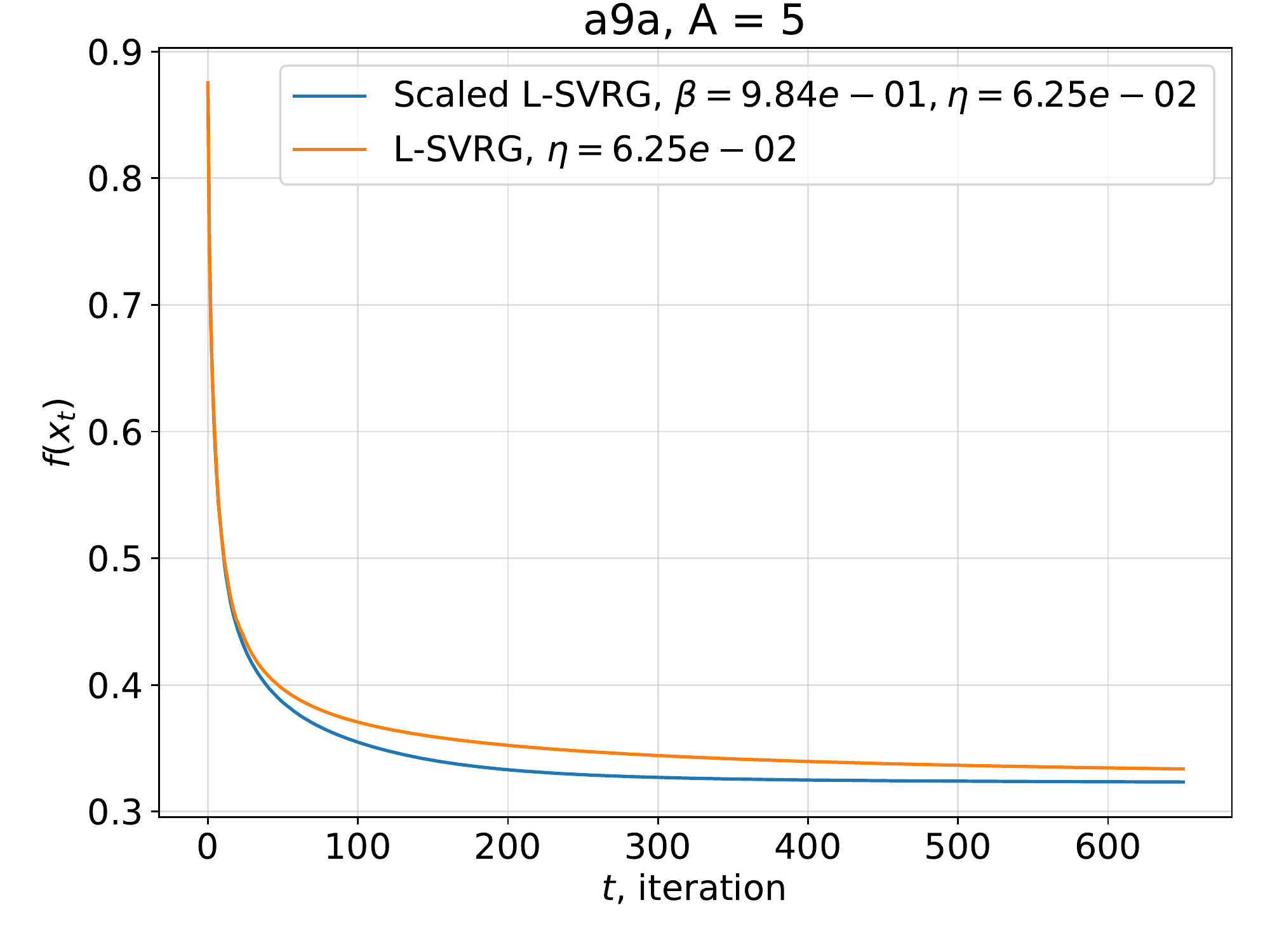}
     \includegraphics[width=0.3\textwidth]{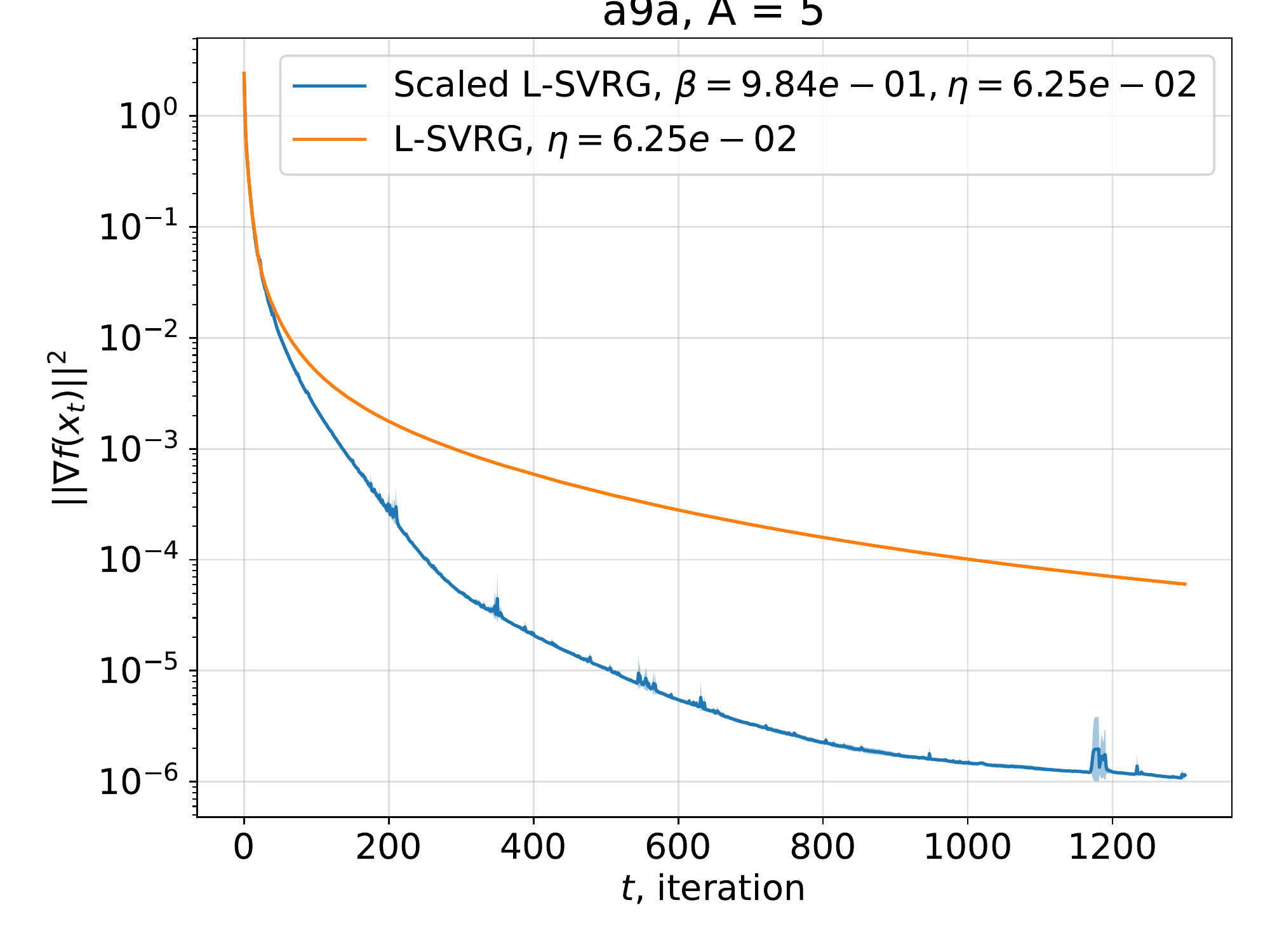}
     \includegraphics[width=0.31\textwidth]{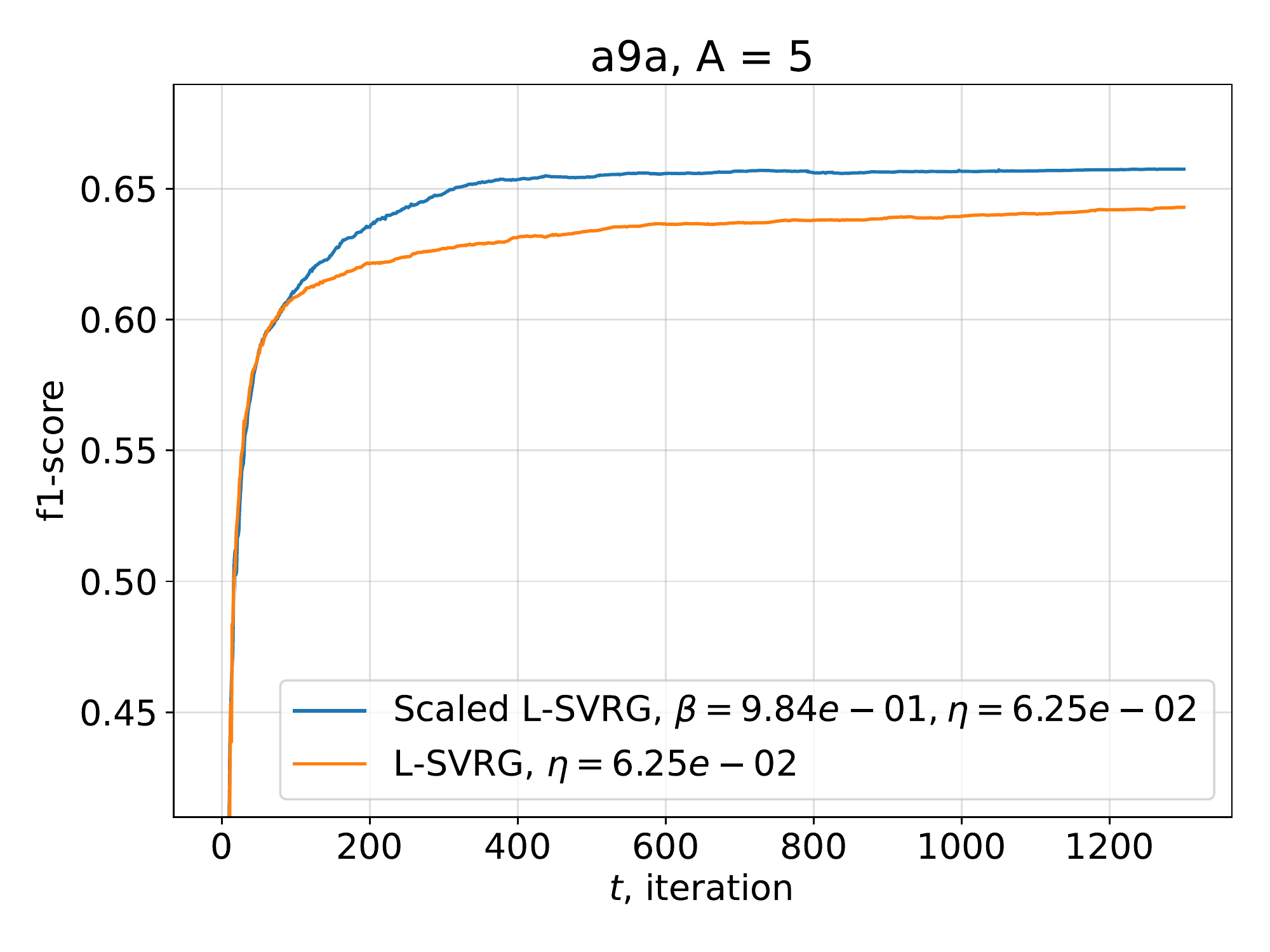}
     \includegraphics[width=0.3\textwidth]{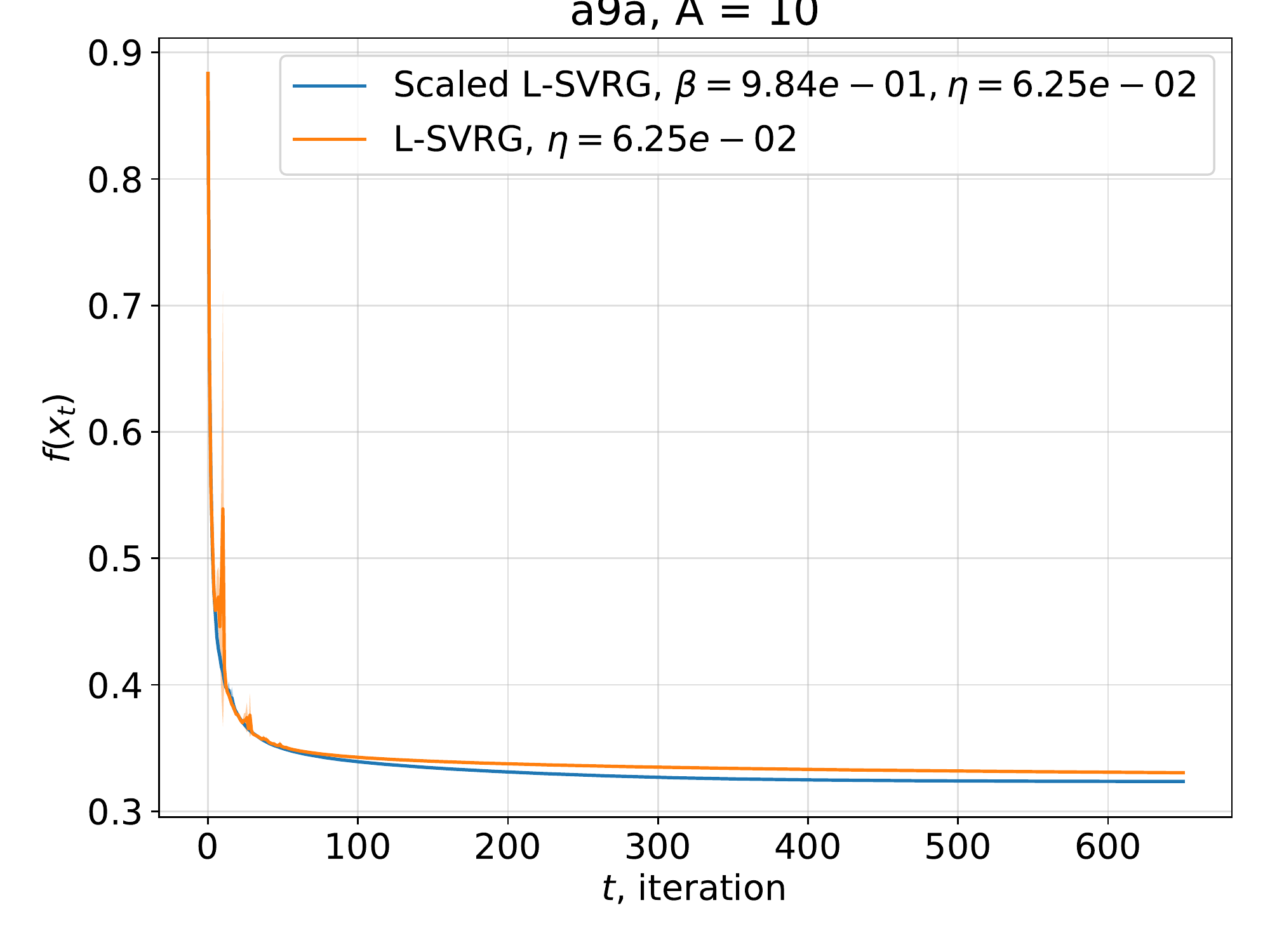}
     \includegraphics[width=0.3\textwidth]{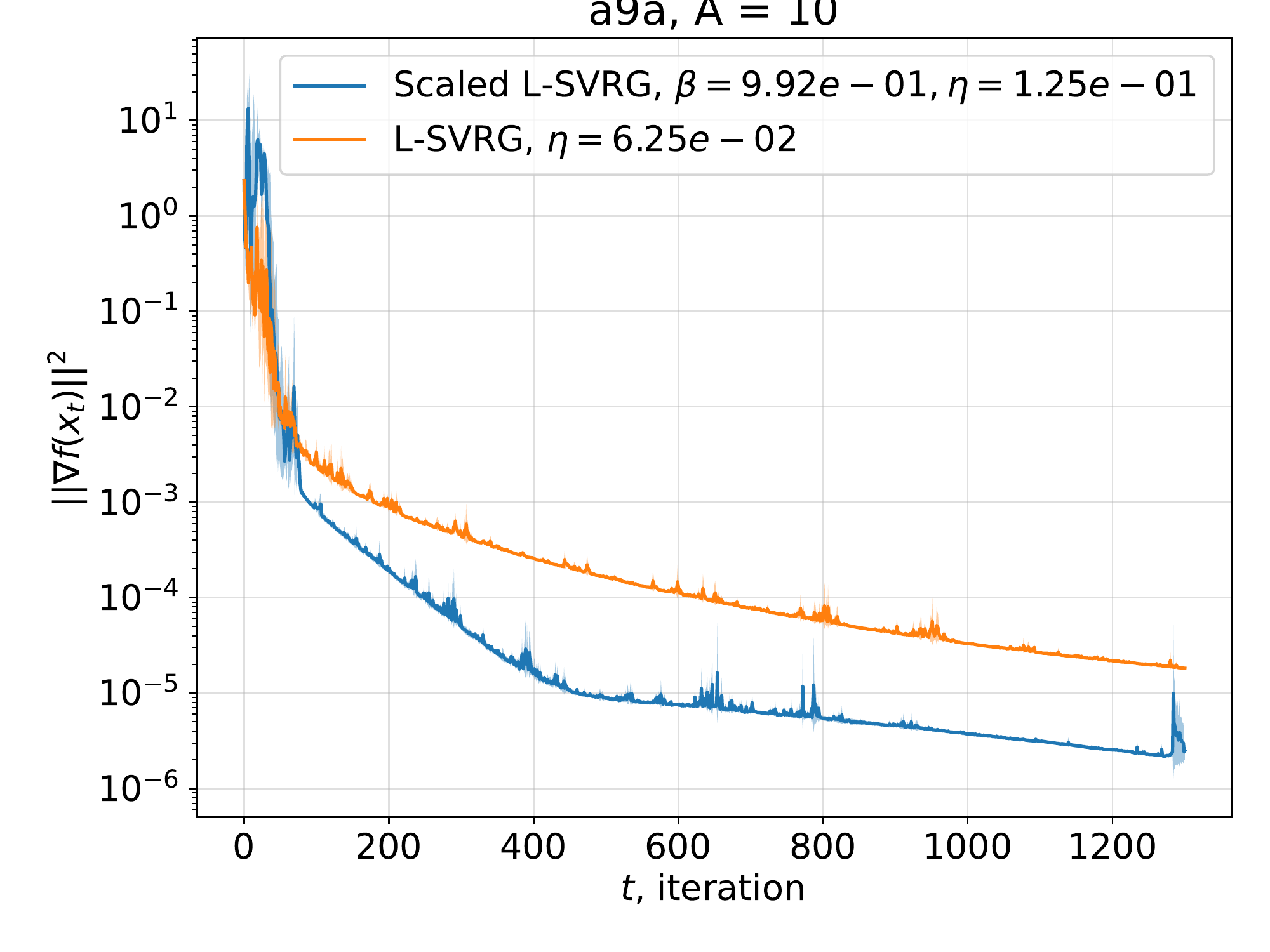}
     \includegraphics[width=0.31\textwidth]{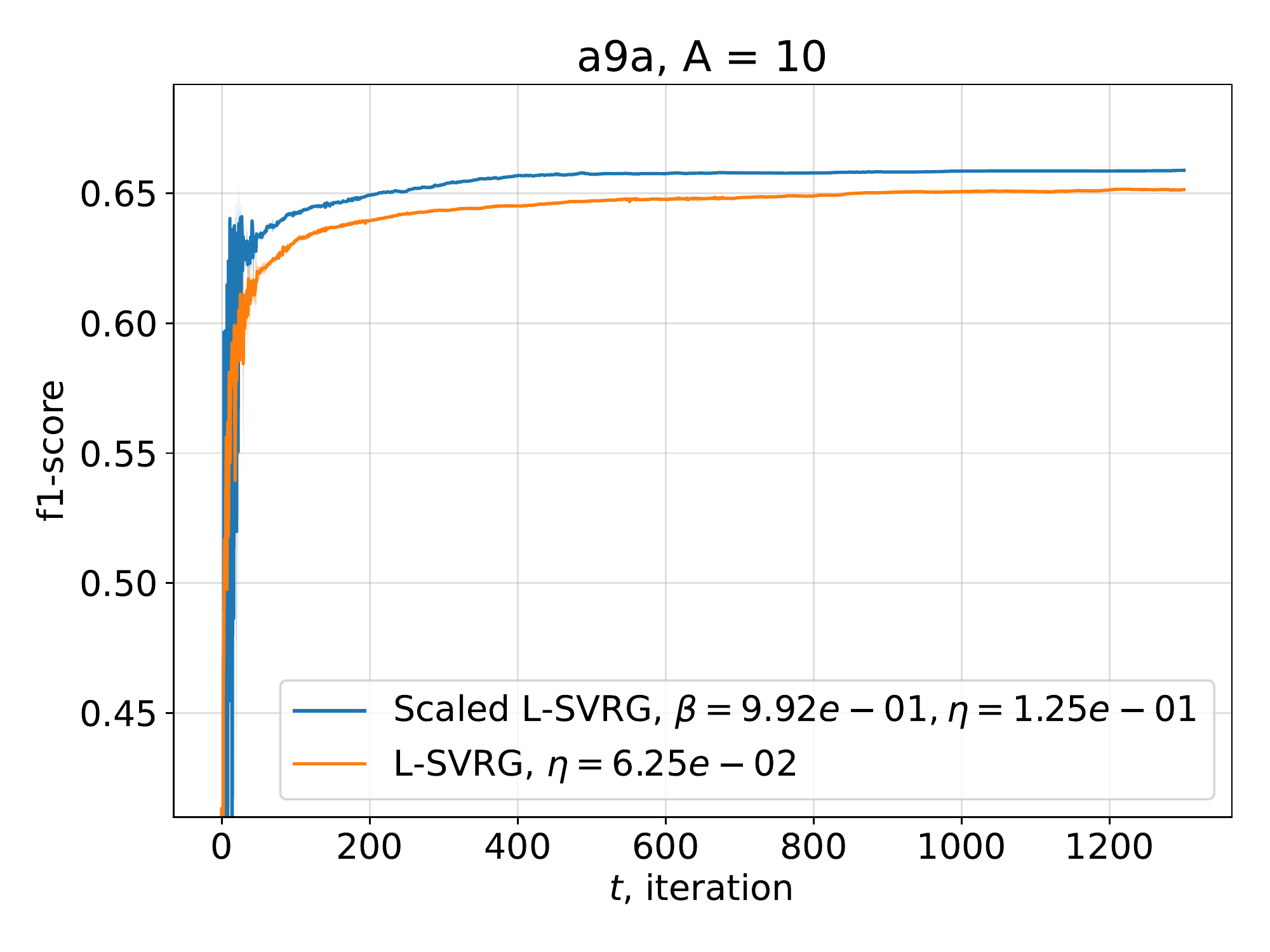}
     \includegraphics[width=0.3\textwidth]{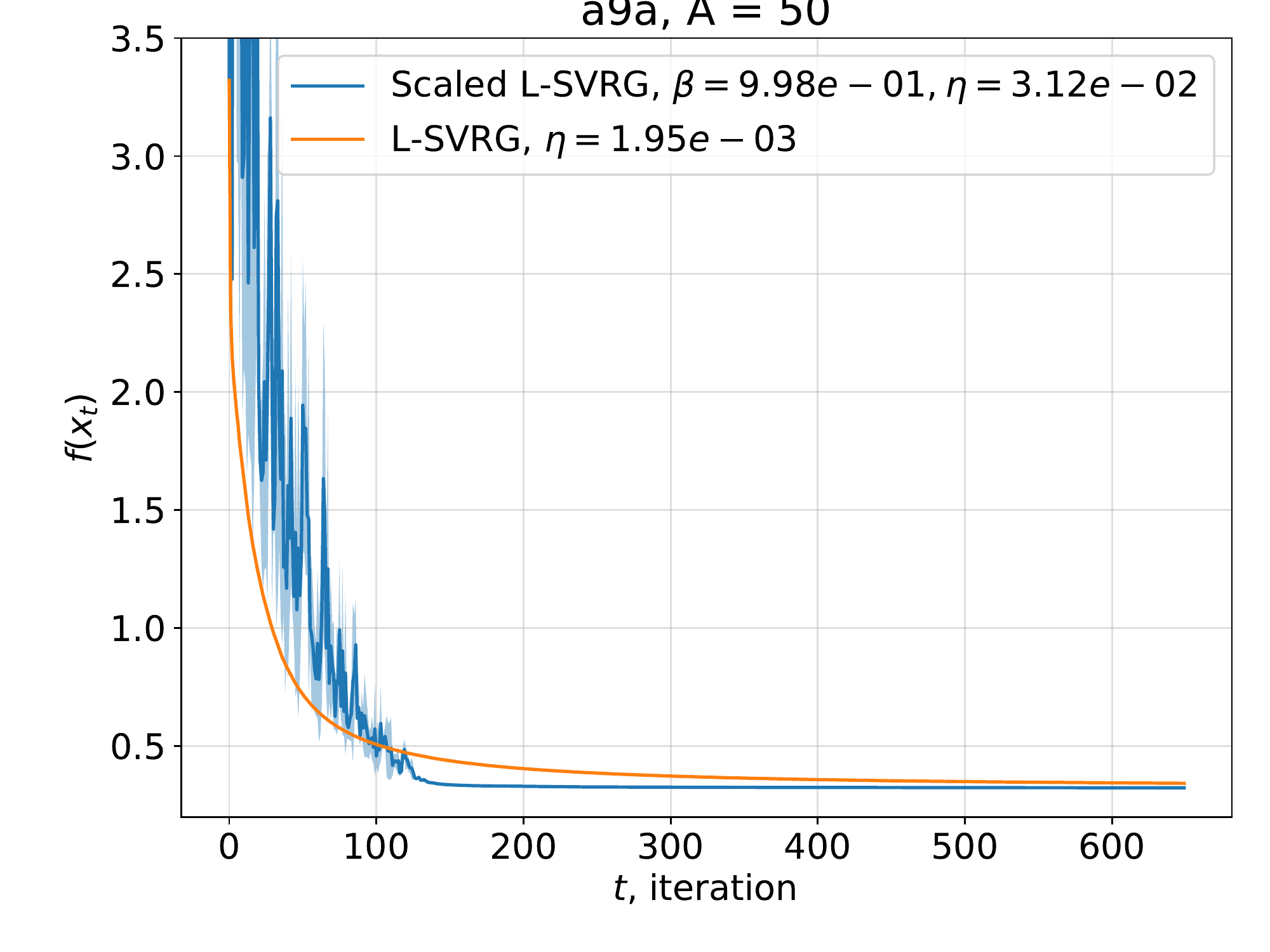}
     \includegraphics[width=0.3\textwidth]{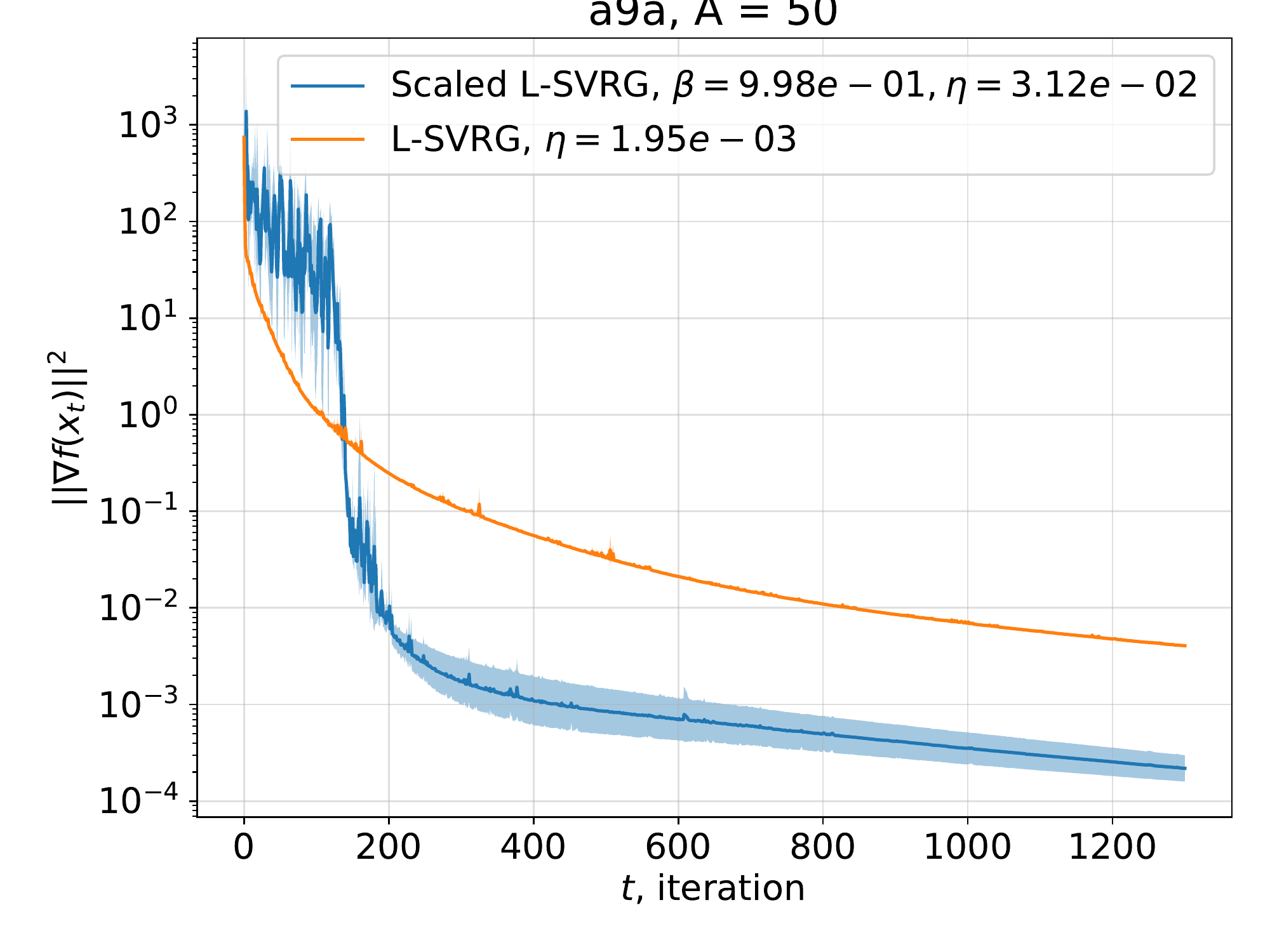}
     \includegraphics[width=0.31\textwidth]{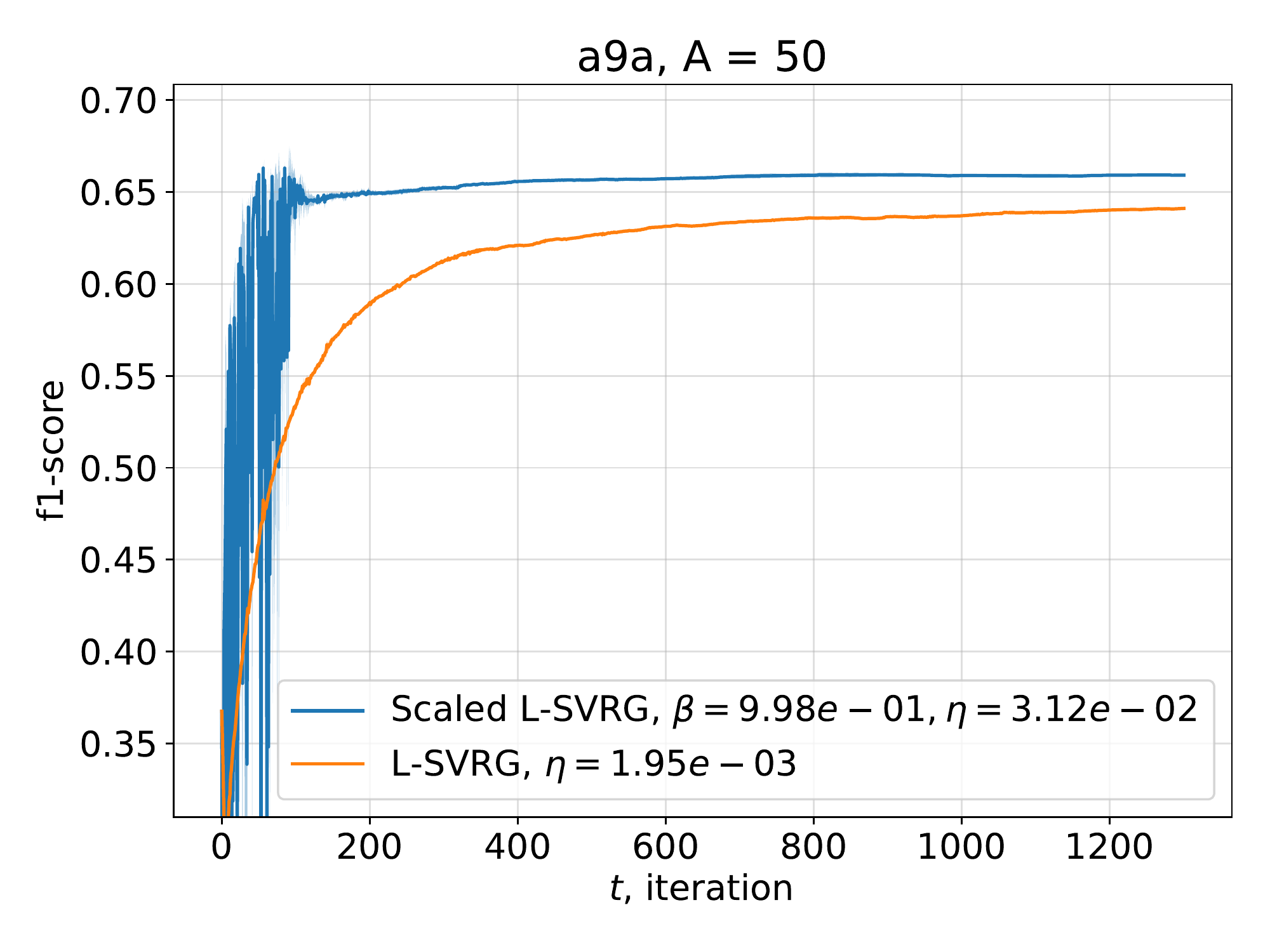}
    \caption{Convergence curves of \texttt{L-SVRG} and \texttt{Scaled L-SVRG} with optimal choice of $\beta_t \equiv \beta$ and $\eta_t \equiv \eta$.}
    \label{conv}
\end{figure}

Figure~\ref{conv} shows the results of comparison of \texttt{Scaled L-SVRG} and ordinary \texttt{L-SVRG} for $A \in \{0.1, 5, 10, 50\}$. Horizontal axis measures the number of iterations (stochastic gradient evaluations), vertical axis measures the objective function value $f(x_t)$, which is more practically interesting, or squared norm of the gradient $\|\nabla f(x_t)\|^2$, which is main for the theory in non-convex case. Convergence curves show the average value of quantity measured for 3 runs with random sequences of batches and are equipped with transparent shades of the size of standard deviation of the measurements. One can see that \texttt{Scaled L-SVRG} converge significantly faster than \texttt{L-SVRG} in all the cases. \texttt{Scaled L-SVRG} allows one to choose bigger step size even if its value is the same for all the iterations. Such a significant superiority of \texttt{Scaled L-SVRG} in the case of $A = 0.1$ might seem to be unexpected, because scaling with $A < 1$ leads to decreasing of Lipschitz constant and increasing of effective step size $\propto 1/L$ whilst scaling introduced by \texttt{Scaled L-SVRG} seeks to eliminate this effect. Nevertheless, scaling in algorithm turns out to be efficient through component-wise adaptivity~--- we encourage this effect by scaling features with random factors $a$ parametrized by $A$.

Next experiment is devoted to the choice of step size $\eta_t \equiv \eta$ for different values of $A$ ($A \in \{0.1, 5, 10, 50\}$). For each $A$, we set $\beta_t \equiv \beta$ to the best value determined for the \texttt{Scaled L-SVRG} in the previous experiment and consider $\eta \in \{2^{-4}, ..., 2^{-10}\}$. We also do not equip corresponding convergence curves with standard deviation shades in this experiment: it is not so significant here, and for most of runs one can estimate the scale of variance with the unaided eye. 

\begin{figure}[ht!]
\centering
    \vskip-12pt
     \includegraphics[width=0.24\textwidth]{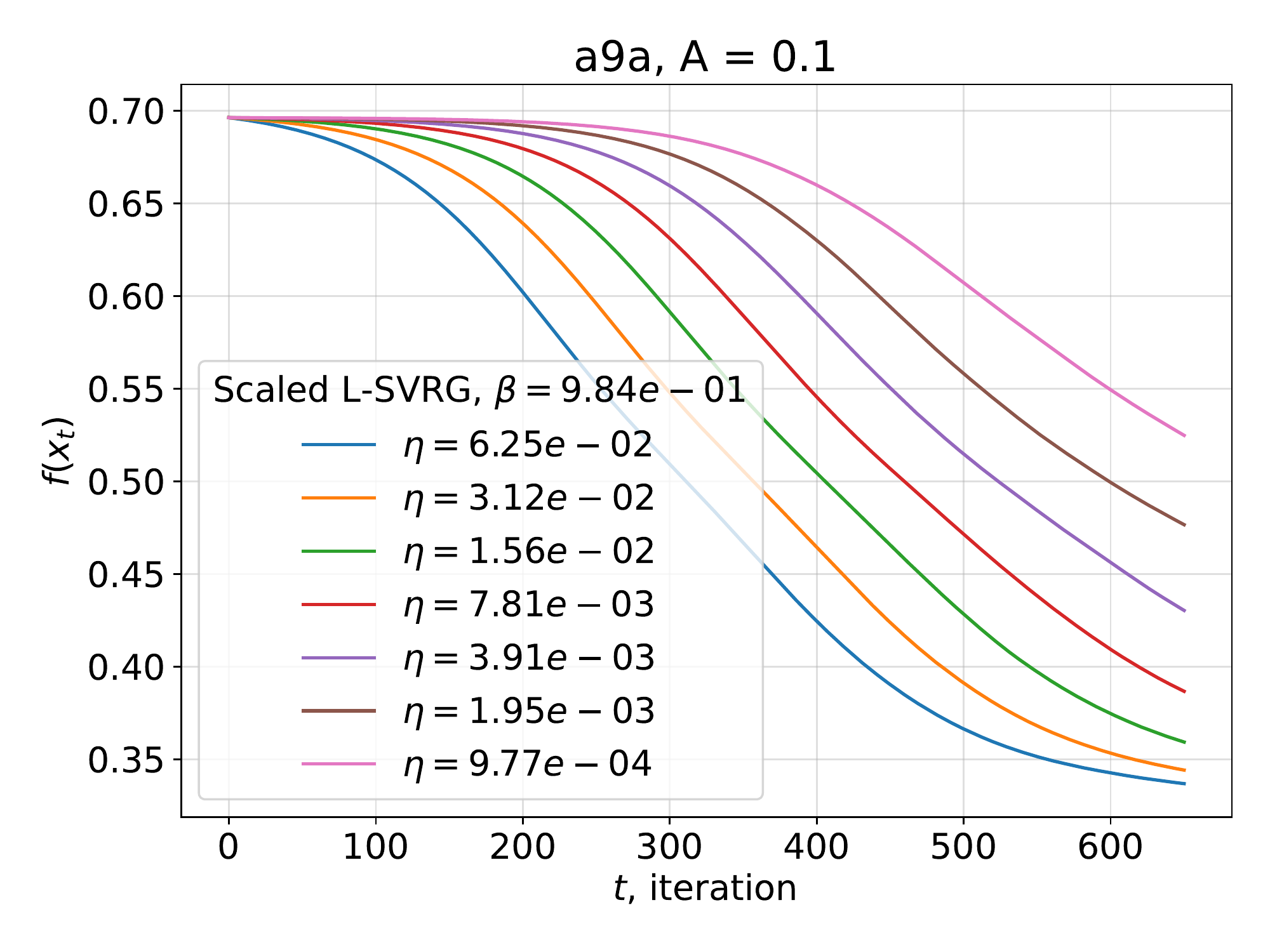}
     \includegraphics[width=0.24\textwidth]{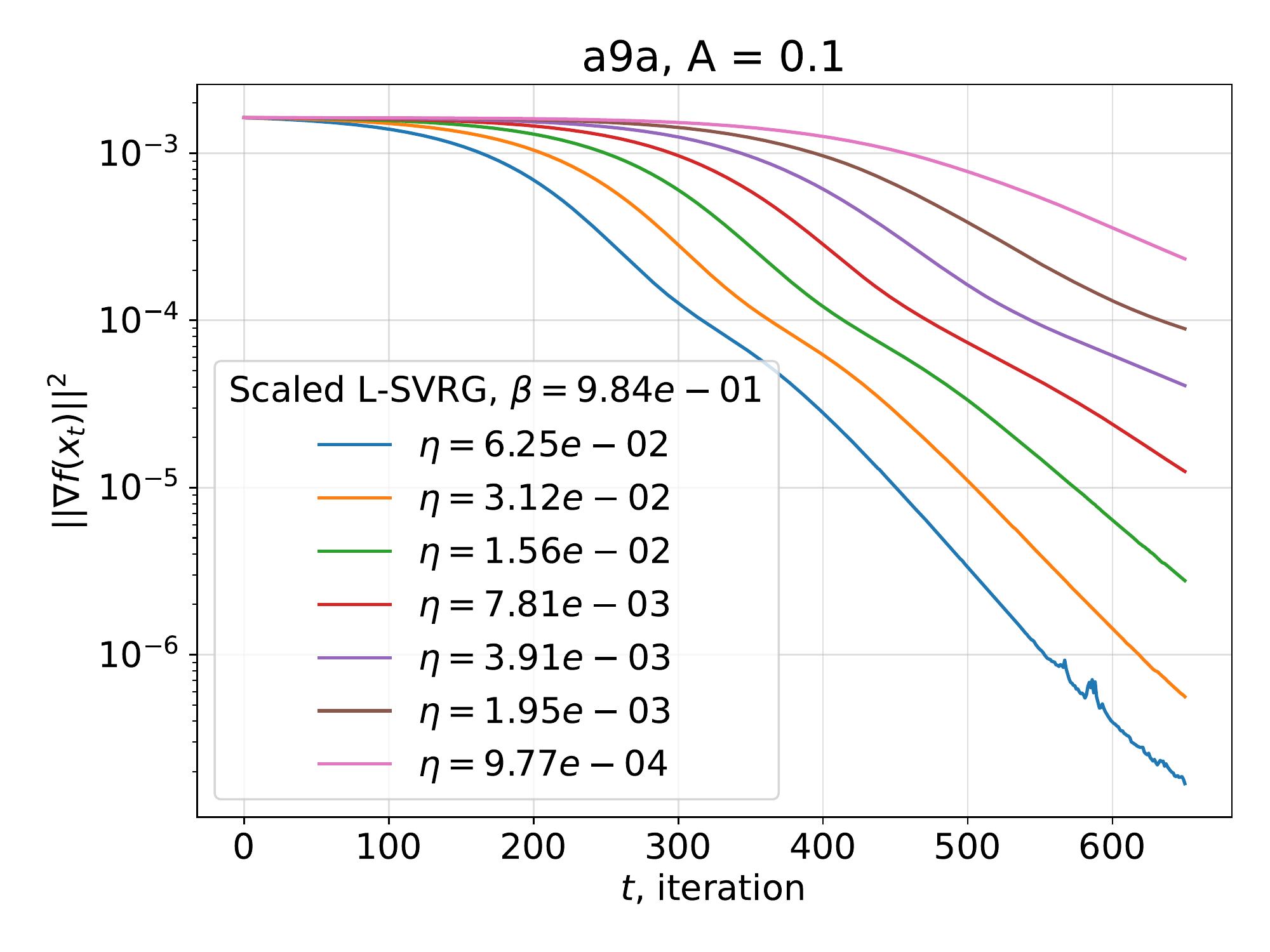}
     \includegraphics[width=0.24\textwidth]{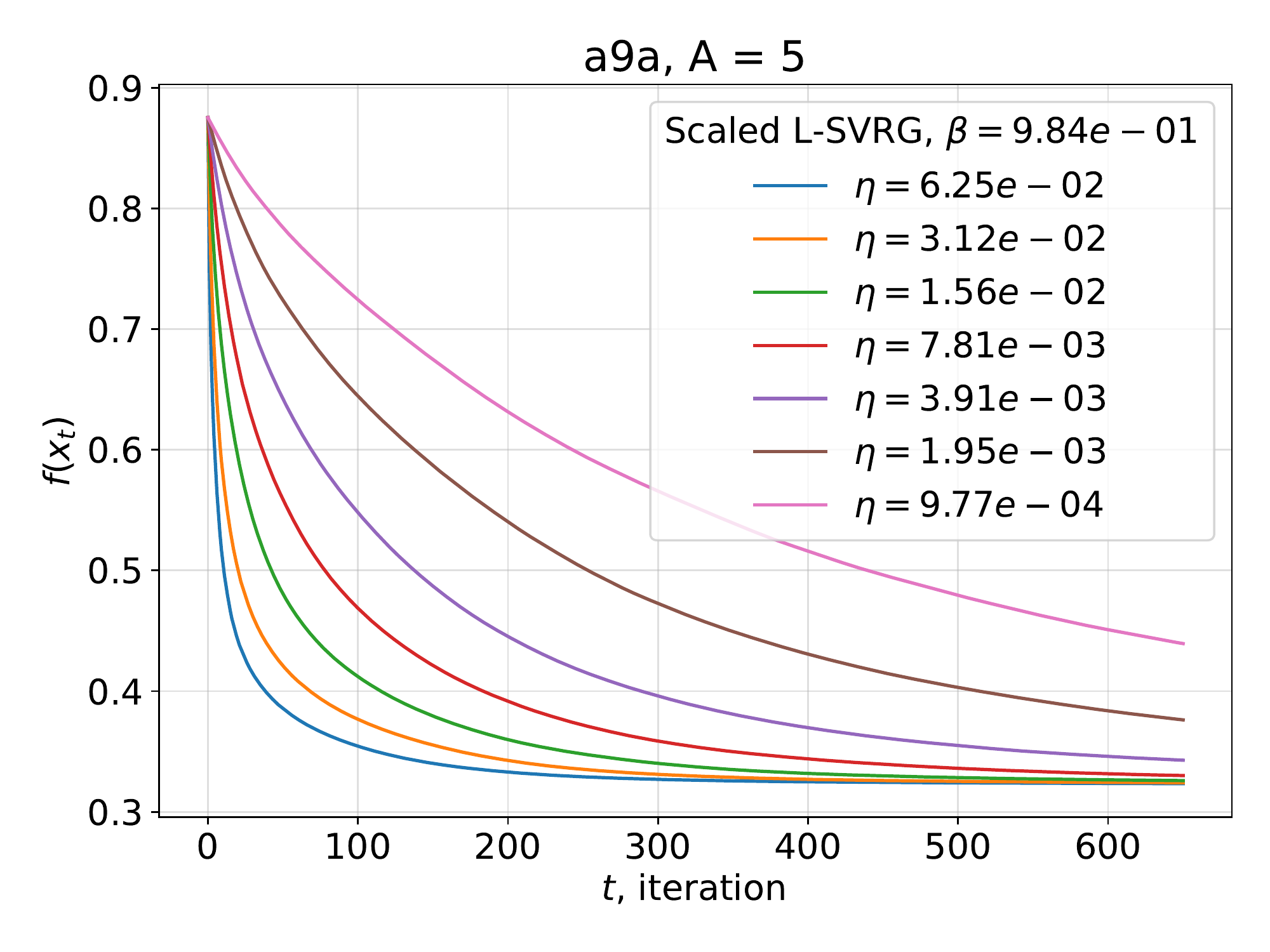}
     \includegraphics[width=0.24\textwidth]{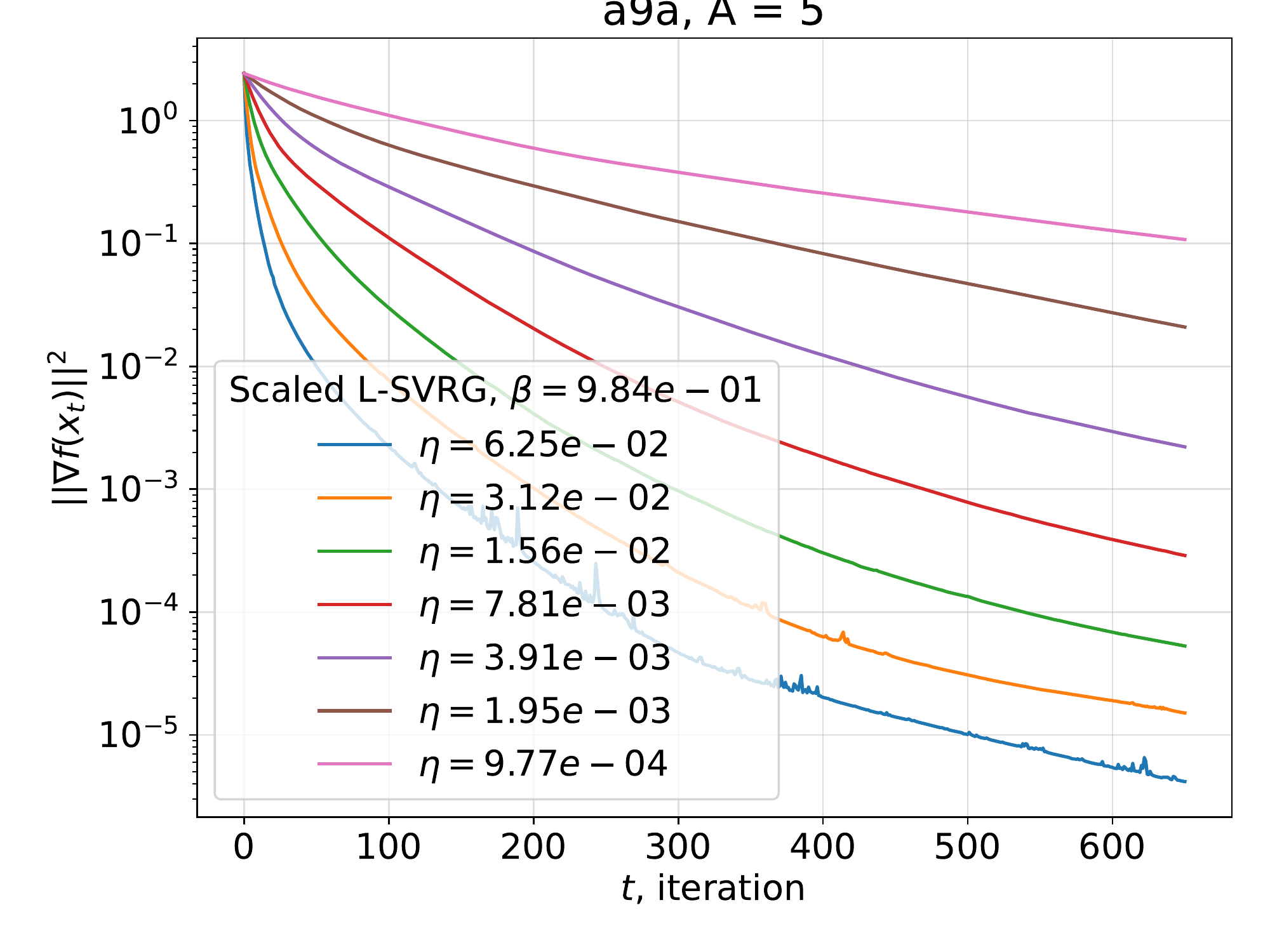}
     \includegraphics[width=0.24\textwidth]{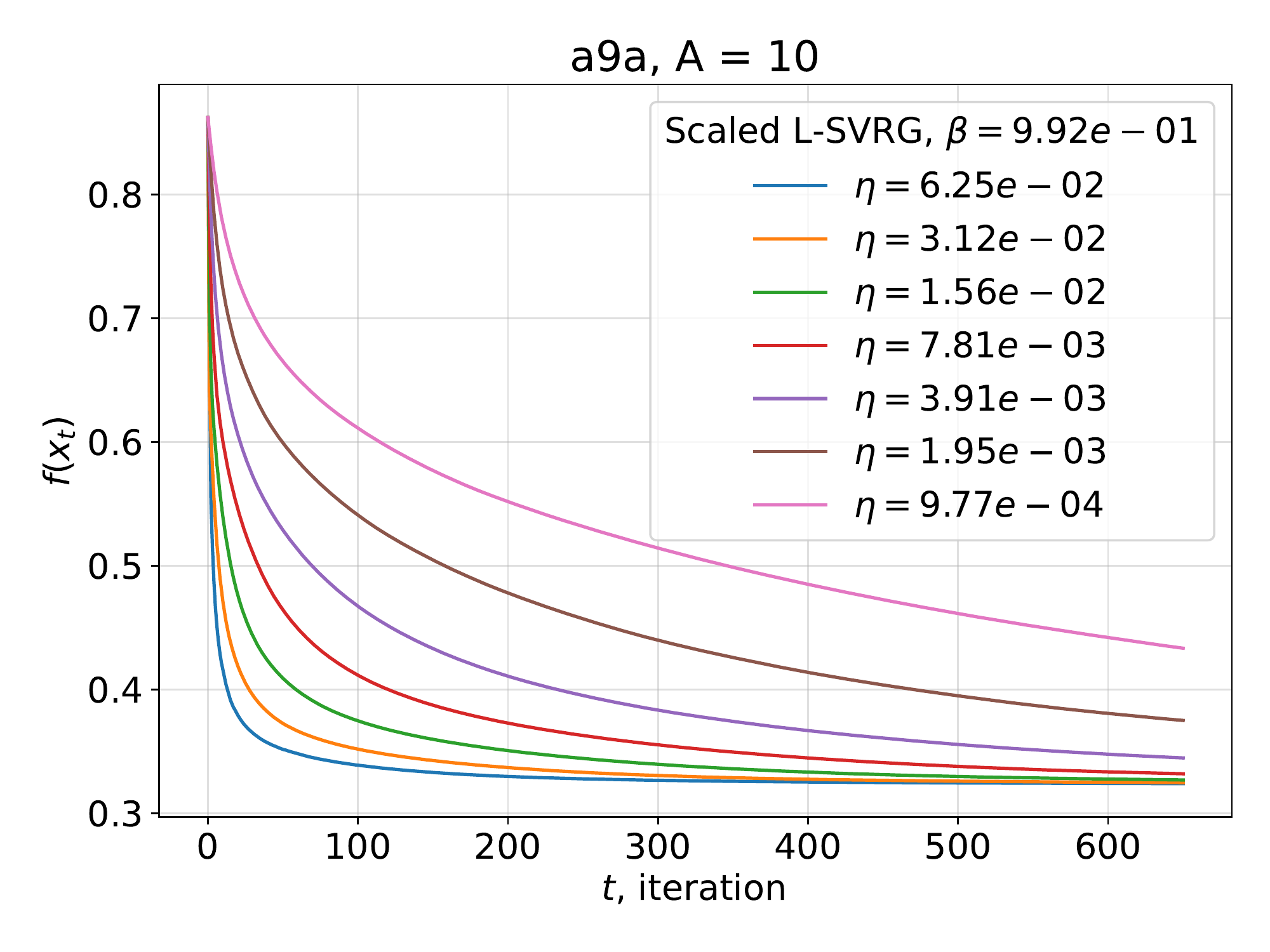}
     \includegraphics[width=0.24\textwidth]{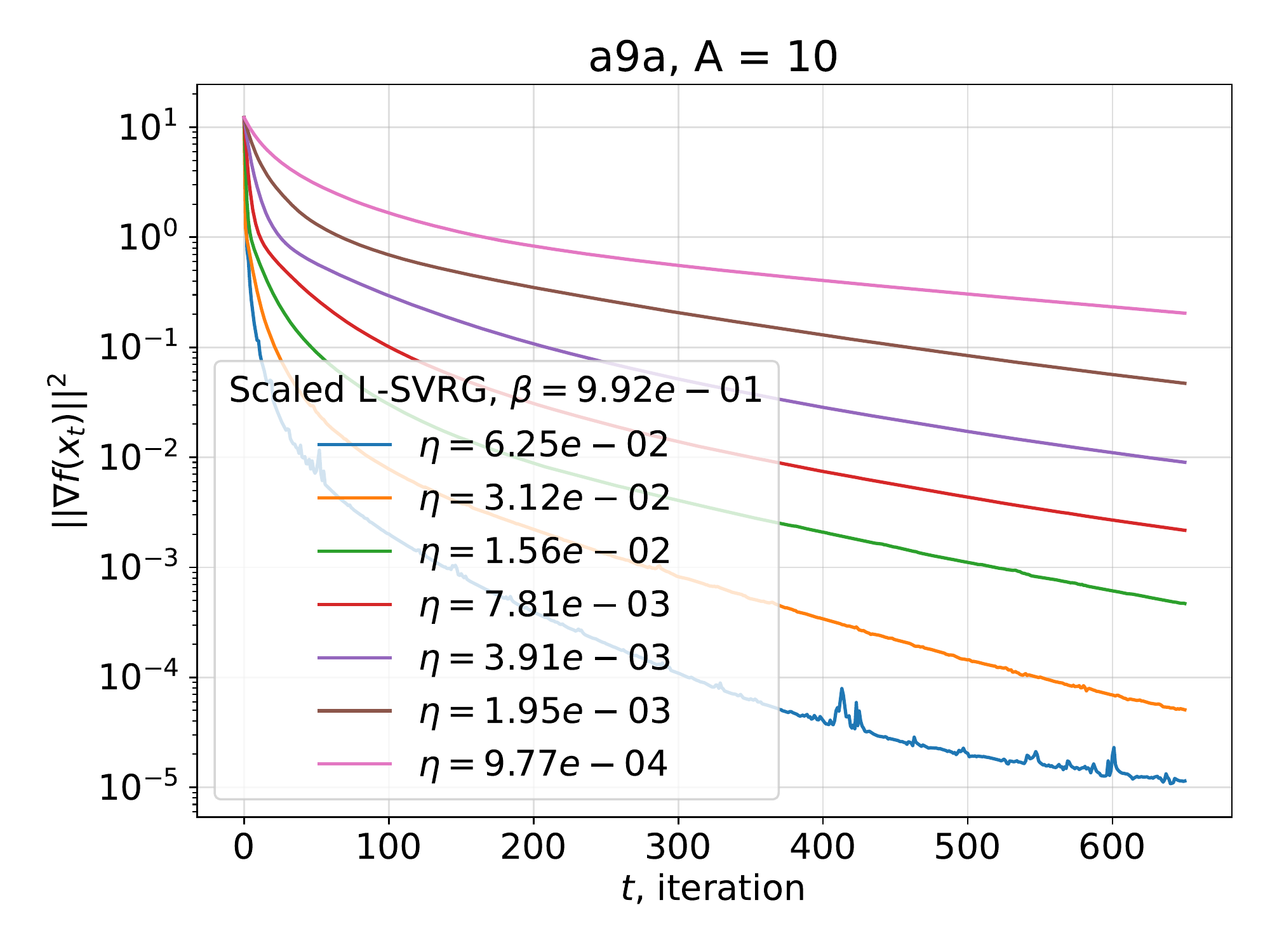}
     \includegraphics[width=0.24\textwidth]{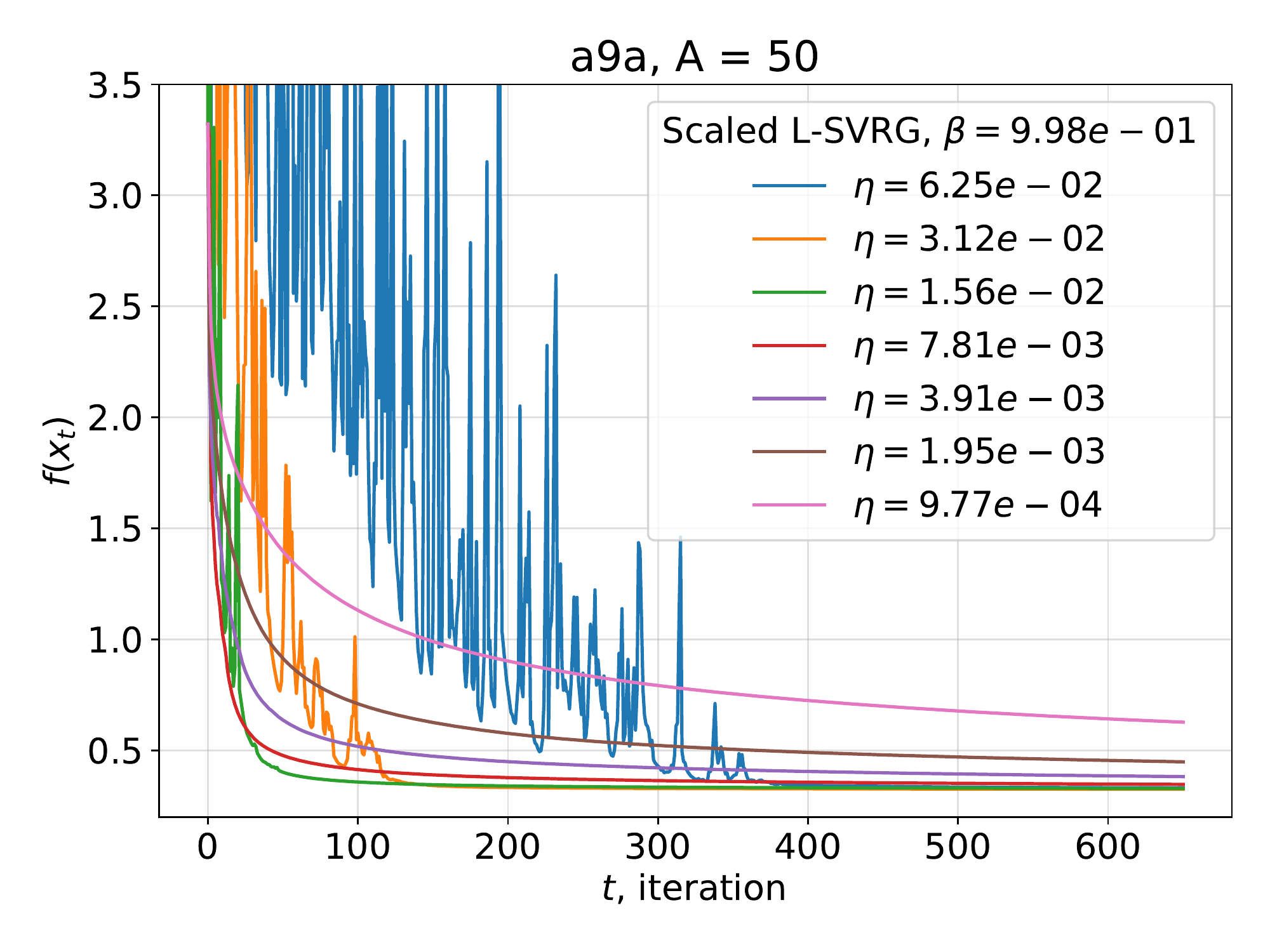}
     \includegraphics[width=0.24\textwidth]{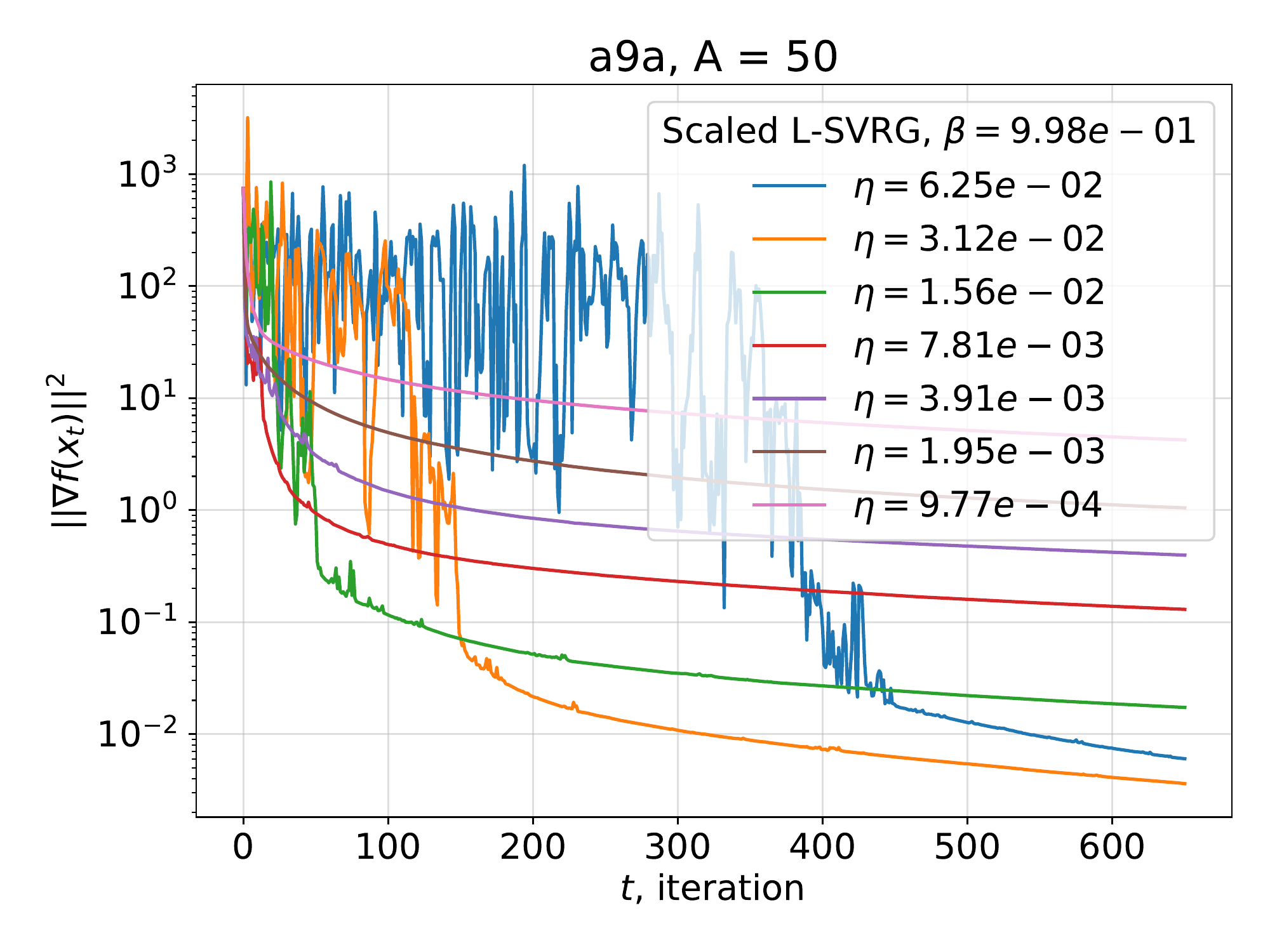}
    \caption{Convergence curves of \texttt{Scaled L-SVRG} with different step sizes on logistic regression problems with different Lipschitz constants.}
    \label{steps}
\end{figure}

Figure~\ref{steps} shows the difference in convergence rate of \texttt{Scaled L-SVRG} in dependence on choice of step size. With the increasing of $A$ (and hence $L$) convergence curves are pressed against the horizontal axis. Its natural, because effective step size is $\propto 1/L$, so efficiency of particular step size $\eta$, getting closer to effective step size, is improving as well, if $\eta$ is small enough. Starting from $A=50$, big step sizes, getting closer to the bound on a step size $\propto 1/L$ guaranteeing the compensation of variance, become inefficient. 

Further, we focus on the behaviour of \texttt{Scaled L-SVRG} algorithm in dependence on the choice of $\beta_t \equiv \beta$, for varying $A$. We consider the case of constant $\beta$ to validate the results obtained in ``One-step effects'' section. 

\begin{figure}[ht!]
\centering
    \vskip-12pt
     \includegraphics[width=0.24\textwidth]{figures/a9a/betas/a9a_A=0.1_different-beta_losses.pdf}
     \includegraphics[width=0.24\textwidth]{figures/a9a/betas/a9a_A=0.1_different-beta_gradients.pdf}
     \includegraphics[width=0.24\textwidth]{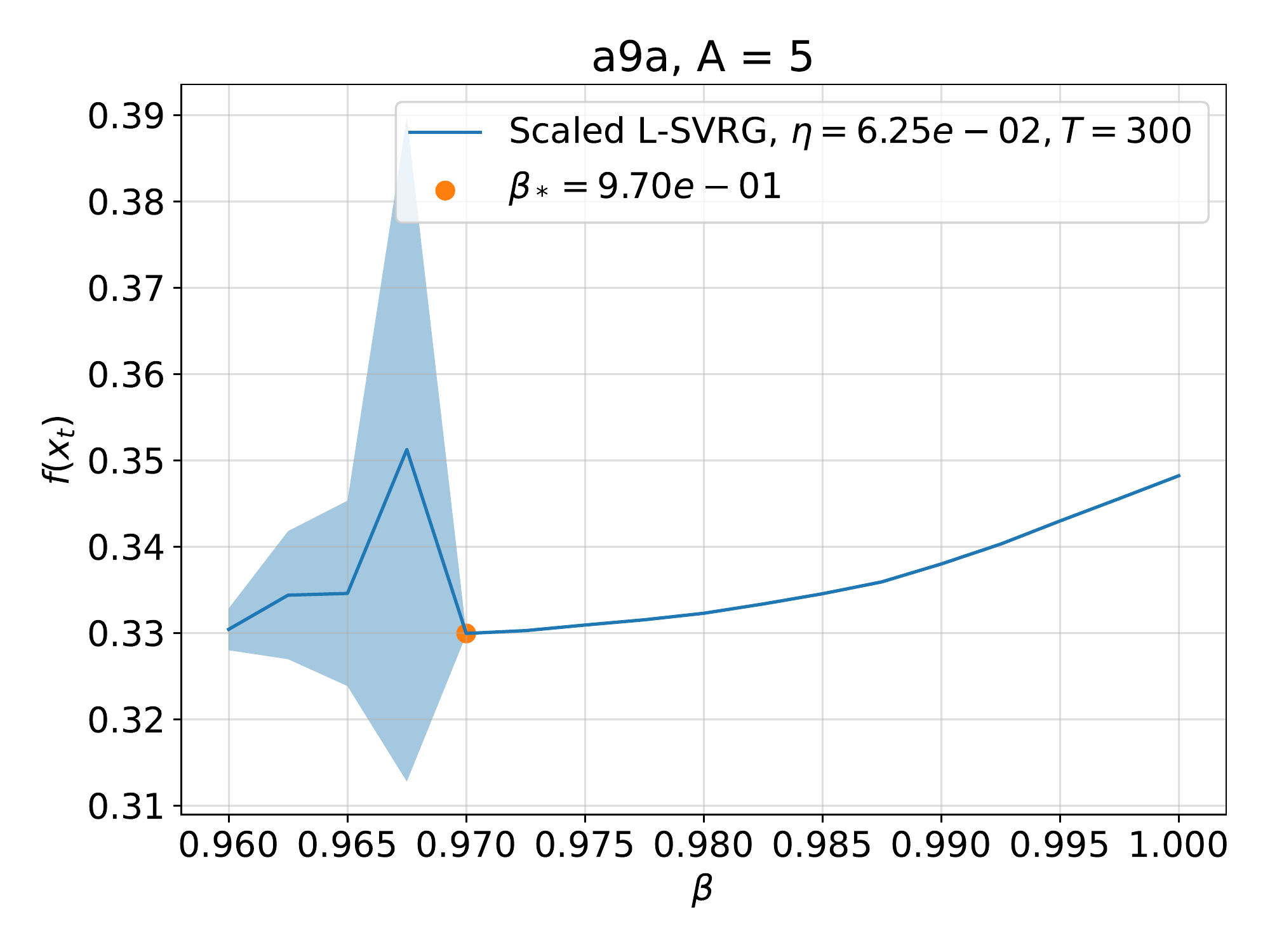}
     \includegraphics[width=0.24\textwidth]{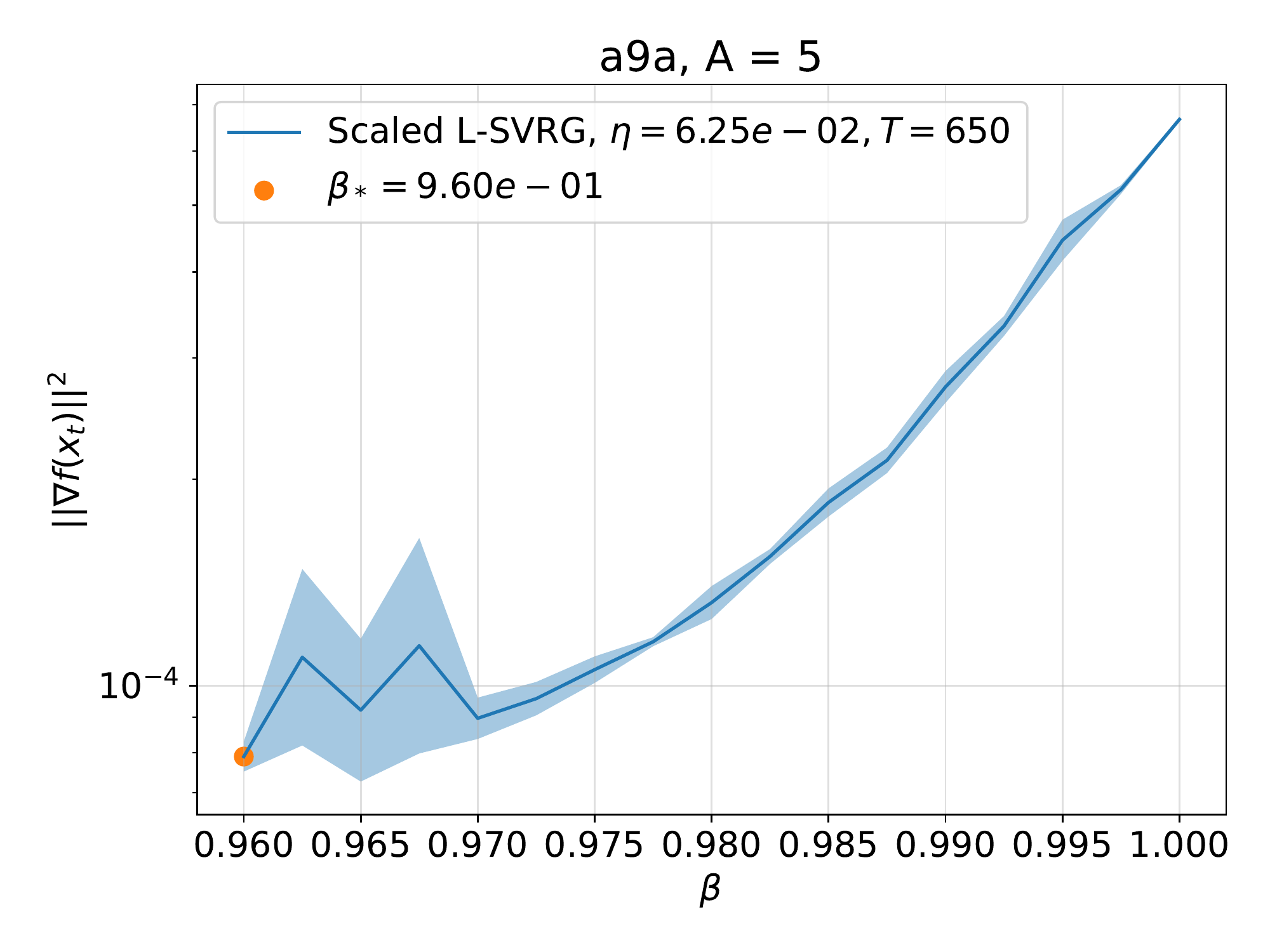}
     \includegraphics[width=0.24\textwidth]{figures/a9a/betas/a9a_A=10_different-beta_losses.pdf}
     \includegraphics[width=0.24\textwidth]{figures/a9a/betas/a9a_A=10_different-beta_gradients.pdf}
     \includegraphics[width=0.24\textwidth]{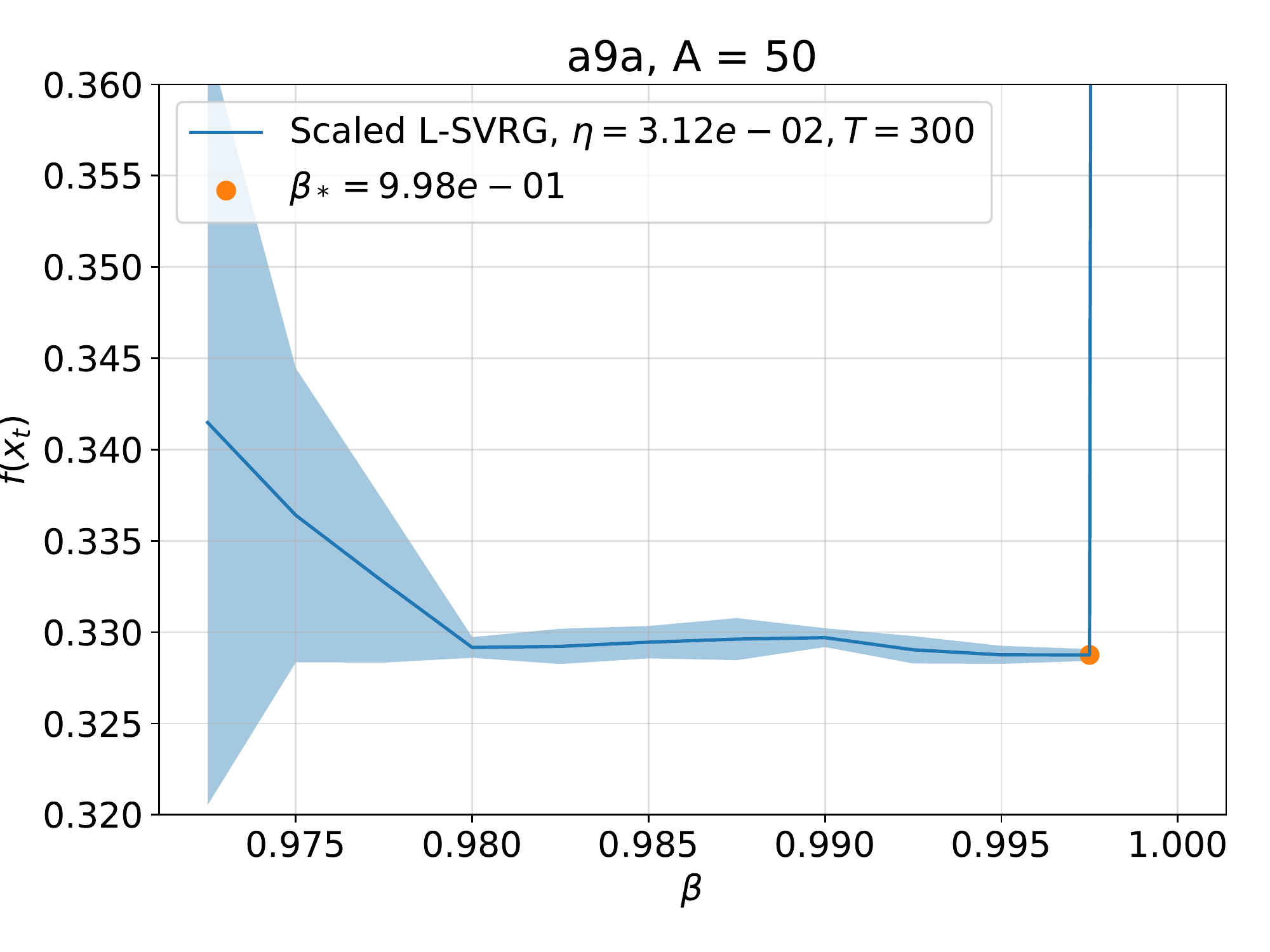}
     \includegraphics[width=0.24\textwidth]{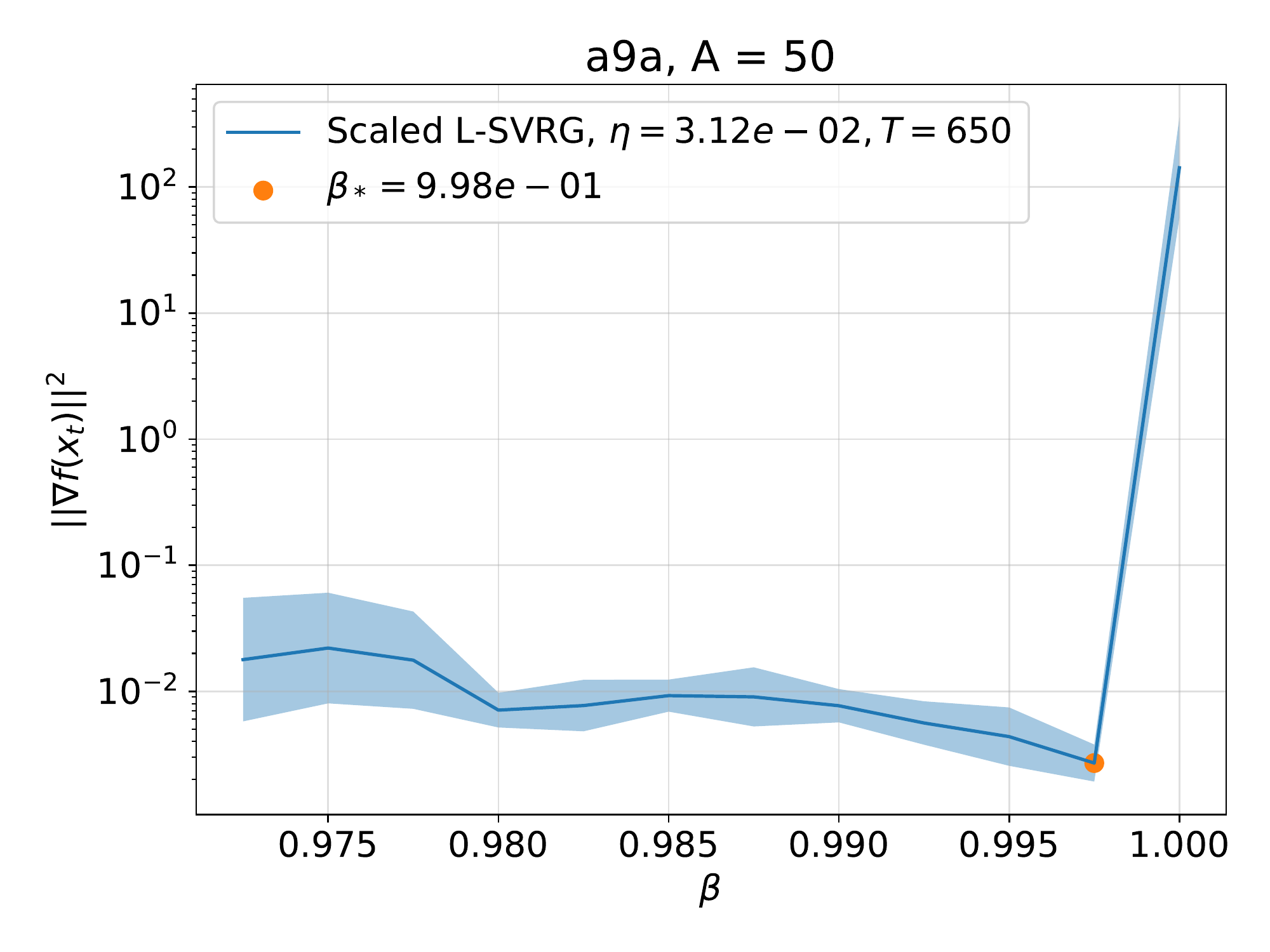}
    \caption{Dependence of achieved precision on $\beta_t \equiv \beta$. (Note: on the lower left plot, function value at $\beta = 1$ is too big, so we cropped the picture, that is why the curve is almost vertical there)}
    \label{betas_full}
\end{figure}

Figure~\ref{betas_full} summarizes the results of \texttt{Scaled L-SVRG} runs with $\beta \in \{0.95, 0.95 + \frac{1-0.95}{20}, ..., 1\}$. Horizontal axis measures value of $\beta$, vertical axis measures objective function value $f(x_T)$ (or squared norm of the gradient $\|\nabla f(x_T)\|^2$) after $T = 300$ iterations of the algorithm. Curves show the average value of quantity measures in 5 runs with random sequences of batches and are equipped with shades of the size of standard deviation of measurements. It can be seen that dependence of achieved precision on $\beta$ changes with increasing of smoothness characteristic $A$: minimum of the corresponding function is getting closer to $\beta = 1$, its values on the left from minimum are growing and its growth rate near $\beta = 1$ is significantly increasing (which is especially noticeable for $A = 50$). This relationship between $\beta$ and $L$ (through $A$) reflects the trade-off between variance compensation and scaling gradients Lipschitz constant. Variance affects the convergence if $\beta$ is small (this fact also leads to divergence for too small $\beta$'s); increasing of $L$ leads to the increasing of $\delta^+_t$, which increases the accumulating error term in \eqref{variance}, so, values for the small $\beta$'s grow. To explain the behaviour near $\beta = 1$, it is reasonable to go to $A = 0.1$ case, where variance error terms are insignificant. Values begin to grow rapidly starting from $\beta \approx 0.97$ and stop on some fixed value at $\beta = 1$. This behaviour is described in \eqref{Delta}, where we have shown the $O(1/(1 - \beta))$ growth of gradients Lipschitz constant. The boundedness at $\beta = 1$ can be explained by the proper choice of $P_0$, such that $\delta^-_t \neq 1$, even if $\beta = 1$. Thus, the main outline of our theory is successfully confirmed on the experiment.

\begin{figure}[ht!]
\centering
    \vskip-6pt
     \includegraphics[width=0.4\textwidth]{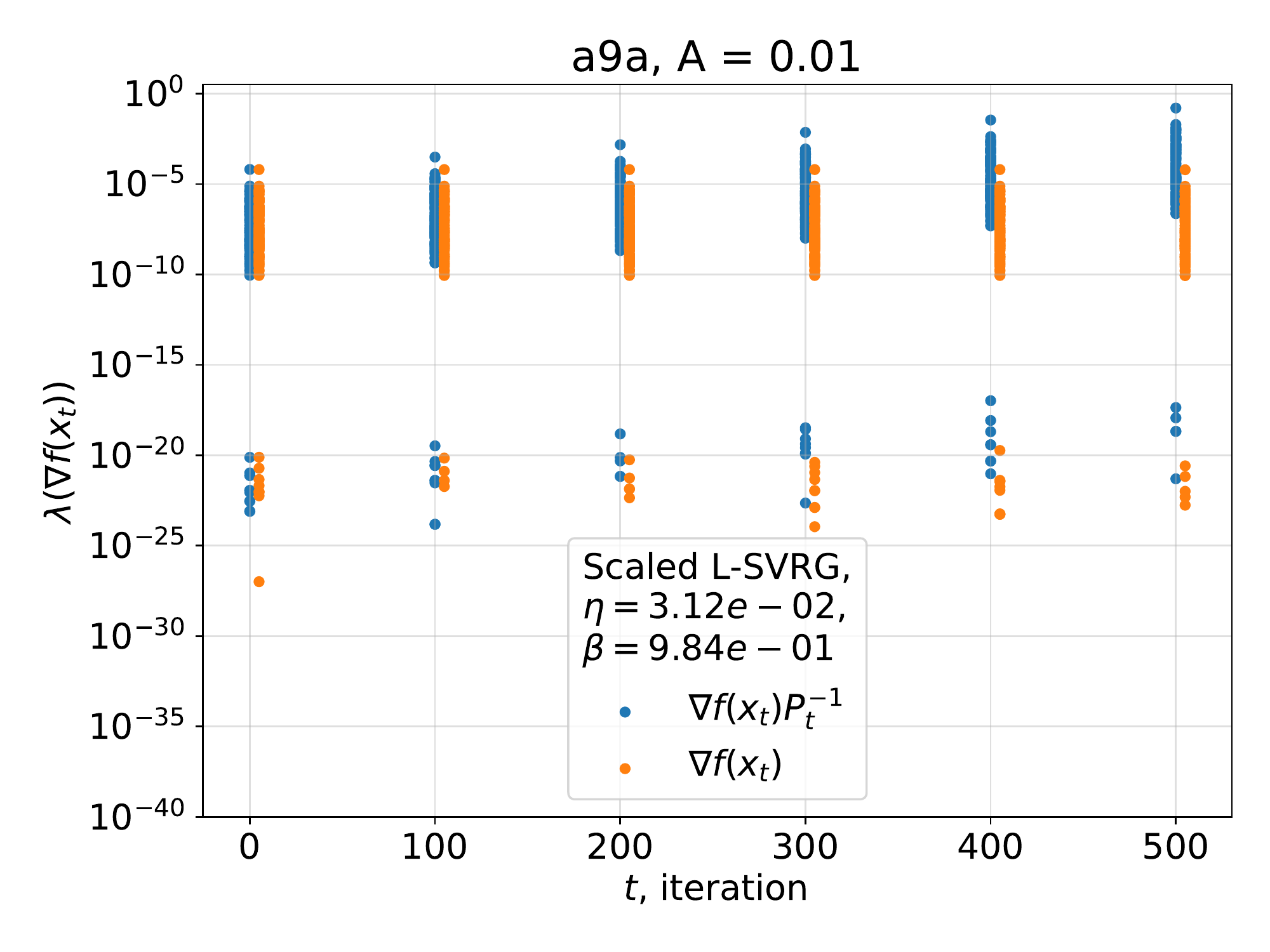}
     \includegraphics[width=0.4\textwidth]{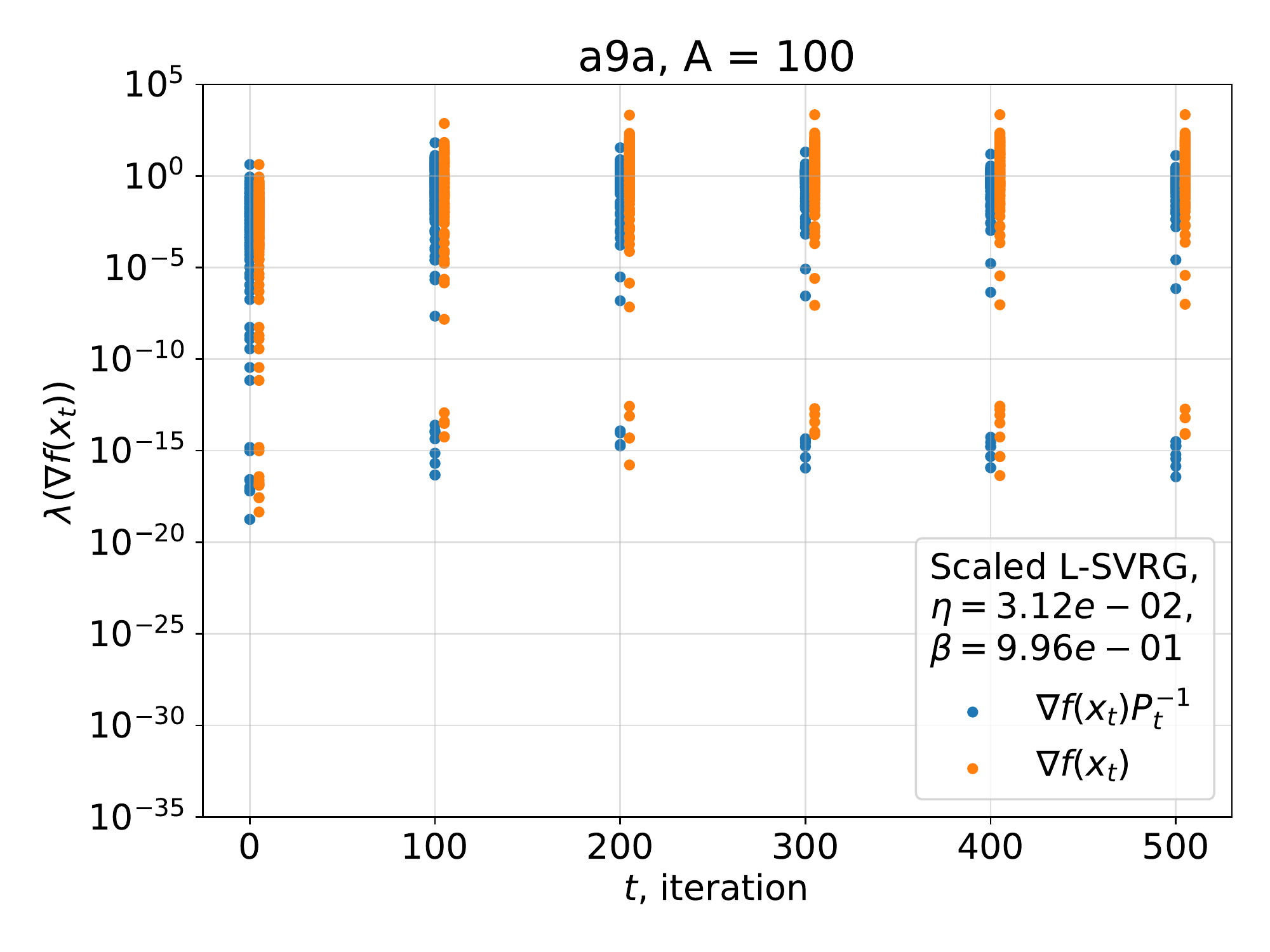}
     \includegraphics[width=0.4\textwidth]{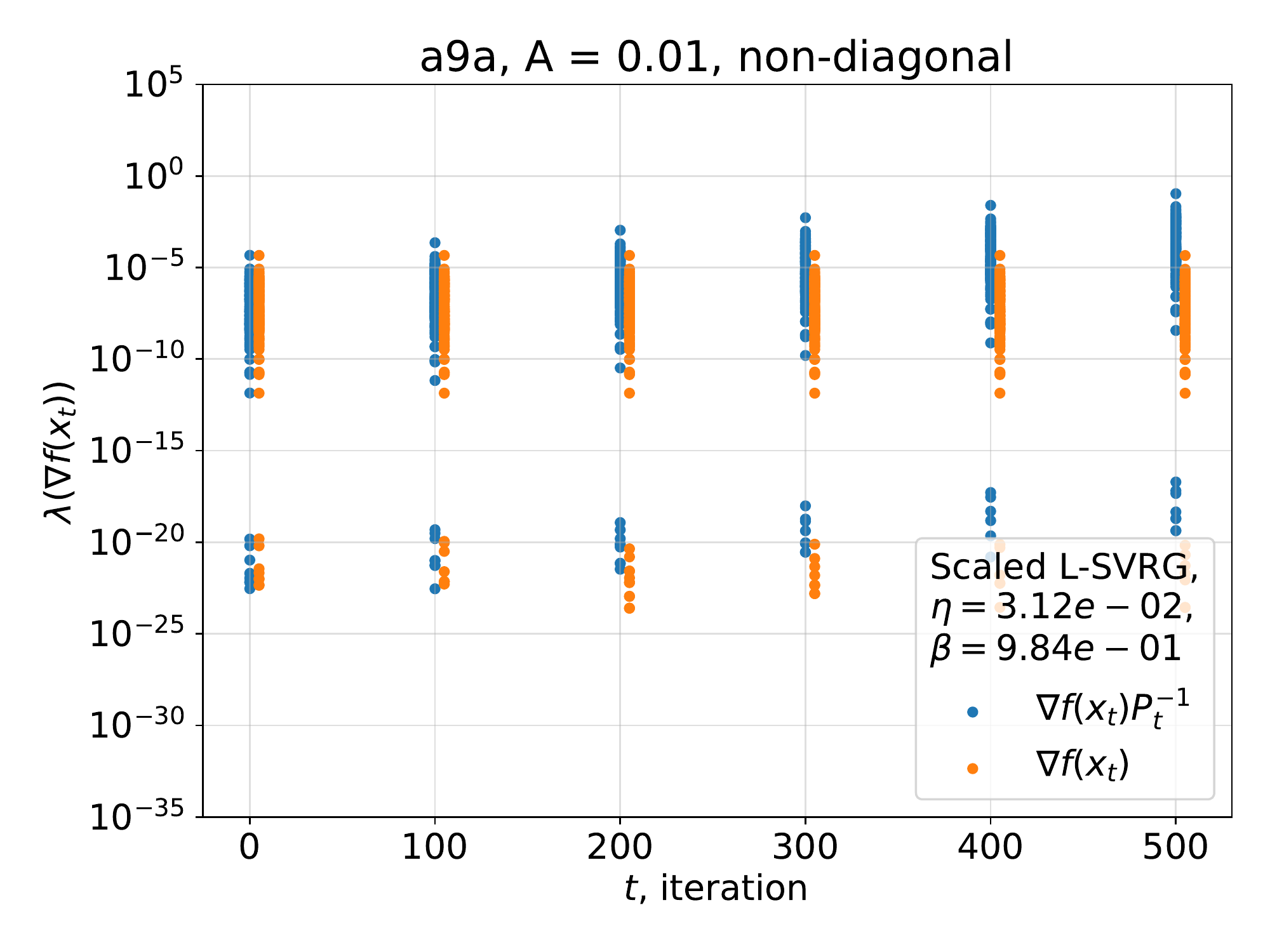}
     \includegraphics[width=0.4\textwidth]{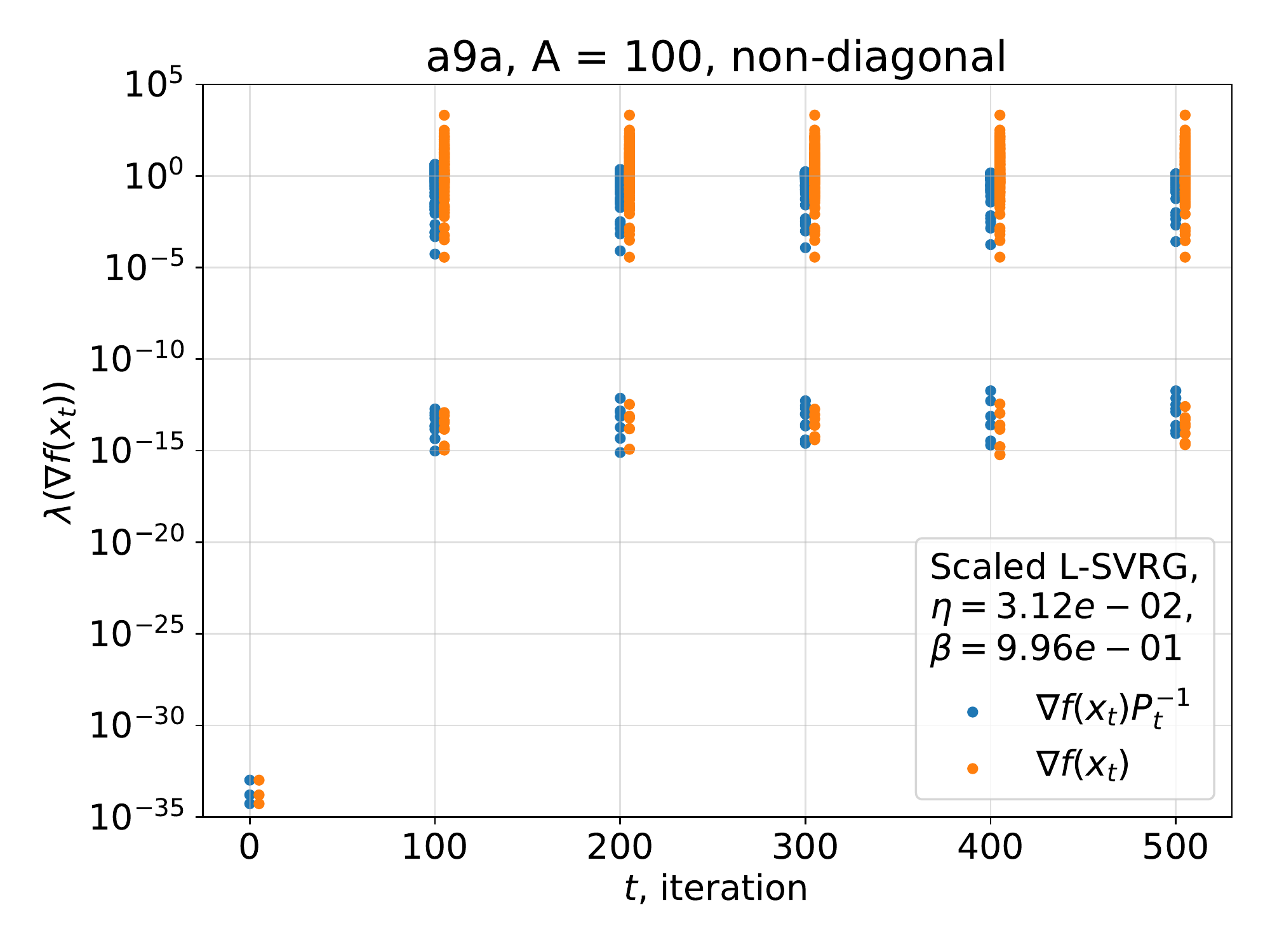}
     \vskip-20pt
    \caption{Dependence of spectrum of Hessian and $\nabla^2 f(x_t) P_t^{-1}$ (characterized by $\Delta_t$) on number of iterations, for diagonal and non-diagonal preconditioning.}
    \vskip-10pt
    \label{spectrum}
\end{figure}

Next part of the experiments is about the patterns in preconditioners and related inexactnesses changing, and corresponding effect on the convergence of \texttt{Scaled L-SVRG} algorithm. Firstly, we consider the dynamic of the spectrum of Hessian and scaled Hessian, that is $\nabla^2 f(x_t) P_t^{-1}$, with increasing number of iterations. The main inexactness $\Delta_t$ upper bounds the largest eigenvalue of the scaled Hessian (minus one), and we do not present curve for $\Delta_t$, because one can estimate it up to the order with the unaided eye. This is a first time we consider non-diagonal preconditioning in our experiments; such an update is defined by $d_t = |\nabla^2 f_{B_t}(x_t)|_\epsilon$, with the same $B_t$ and $\epsilon$, which now requires singular value decomposition of $\nabla^2 f_{B_t}(x_t)$ at the every iteration. Strictly speaking, our theory is not well-suited to this case, but this practical consideration will give us an additional information about behaviour of the algorithm when the smallest and other eigenvalues are scaled in a proper way. 

\begin{figure}[ht!]
\centering
    \vskip-12pt
     \includegraphics[width=0.4\textwidth]{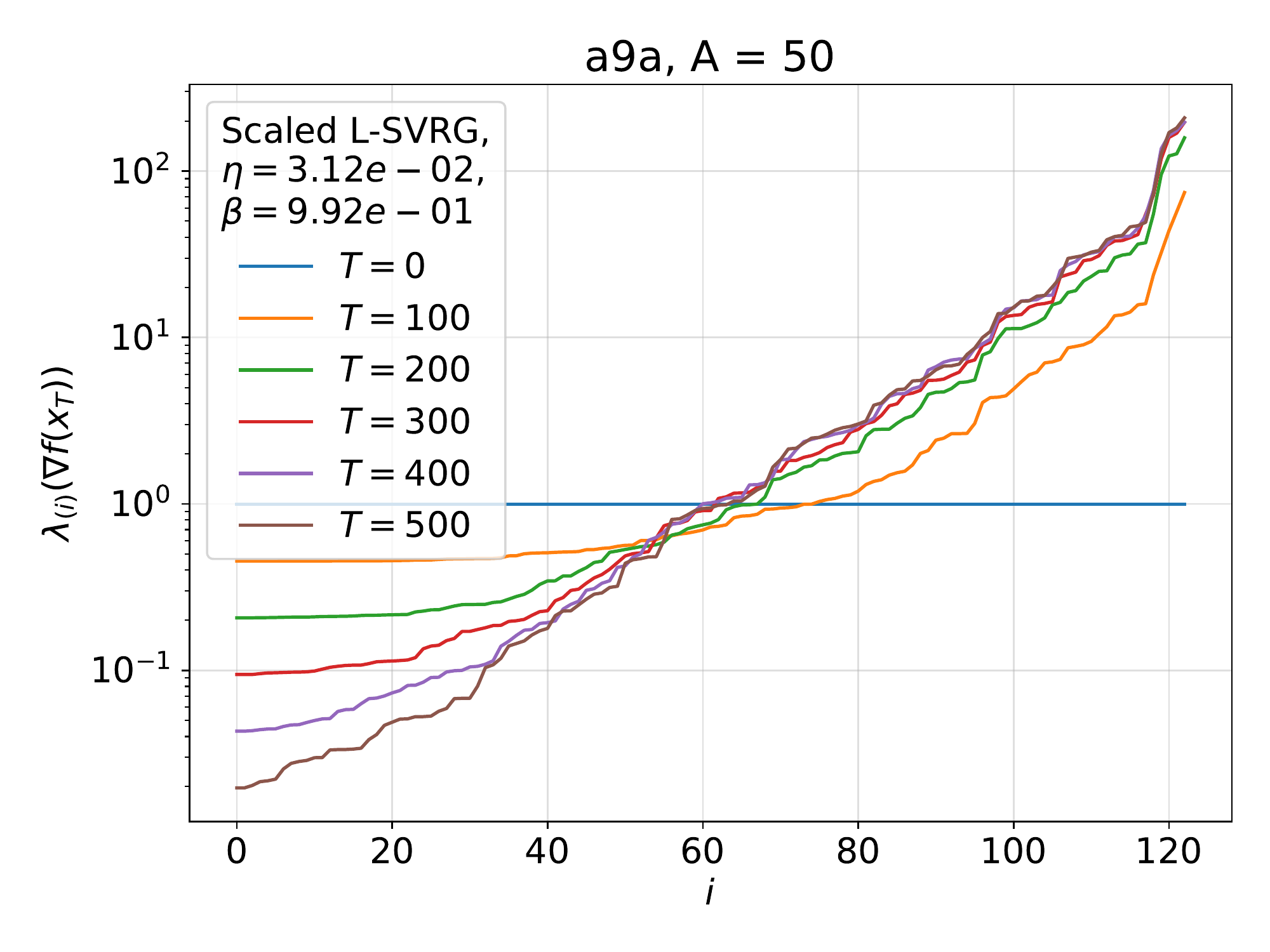}
     \includegraphics[width=0.4\textwidth]{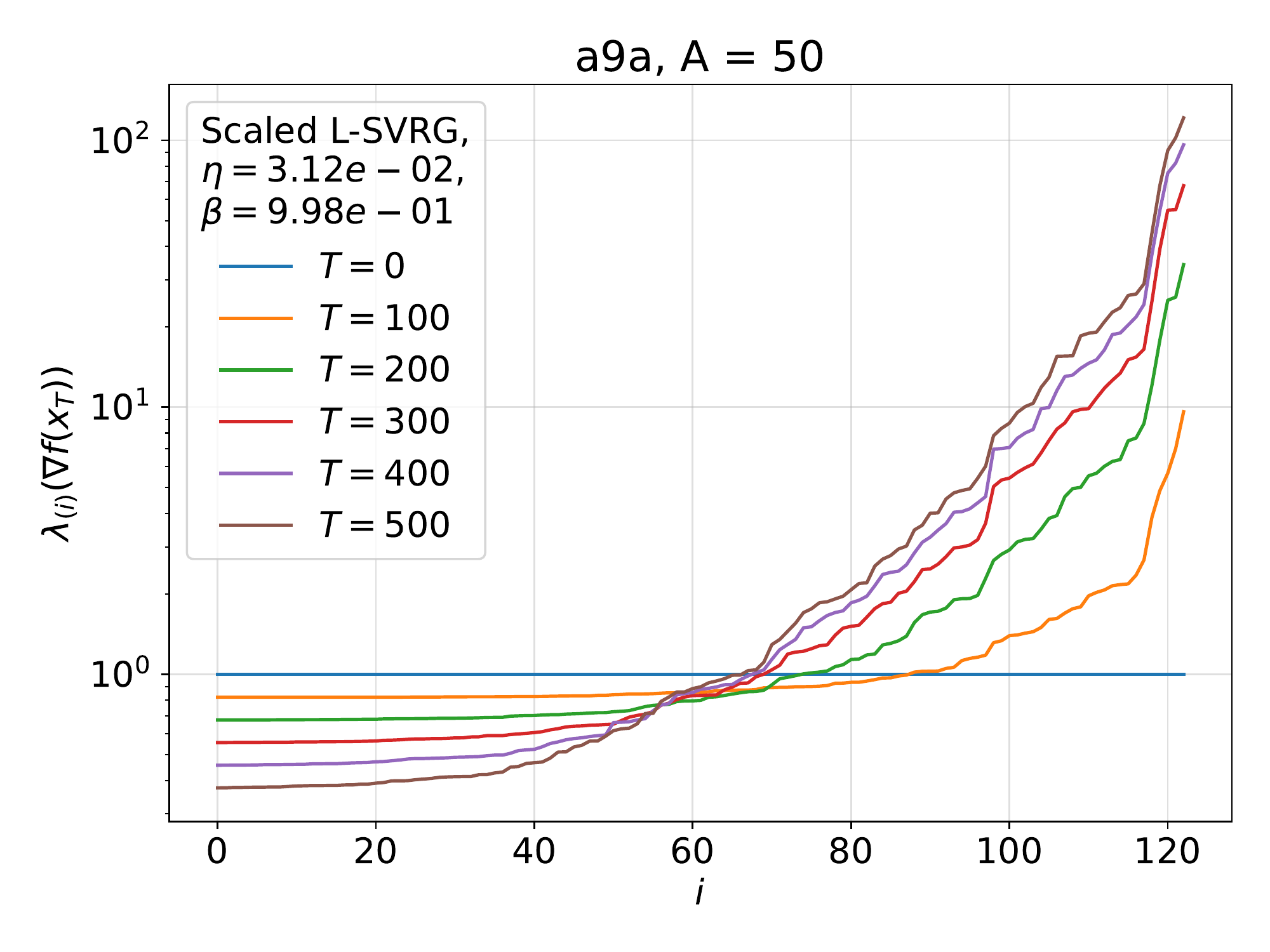}
    \caption{Dependence of spectrum of Hessian approximation on number of iterations.}
    \label{spectrum2}
    \vskip-6pt
\end{figure}

Figure~\ref{spectrum} presents the dynamic of Hessian and scaled Hessian spectrum in two scenarios: diagonal and non-diagonal, for $A \in \{0.01, 100\}$, which is needed to represent both $A < 1$ and $A > 1$ cases. What we see is that the largest eigenvalue of scaled Hessian converges to $1$, which means its increasing in comparison to the largest eigenvalue of Hessian in the case of $A < 1$ and its decreasing~--- in the case of $A > 1$. For the problem we consider, all the eigenvalues of the (scaled) Hessian form two clouds of points of the plot, and whilst the upper cloud is shifted so that the largest eigenvalue tends to $1$, relative position of lower cloud depends on the update type. If we use a diagonal update, lower cloud is shifted in the same direction as the upper one (all the eigenvalues increase or decrease at the same time in the case of $A < 1$ or $A > 1$, correspondingly), which leads to the moving of the smallest eigenvalue of the scaled Hessian away from $1$. It worth to note that diagonal update do not pay enough attention to small eigenvalues that can be seen also from Figure~\ref{spectrum2}. Conversely, if we use non-diagonal update, clouds can be shifted in opposite directions such that both the largest and the smallest eigenvalues of the scaled Hessian converge to $1$.

The following Figure~\ref{optimal_beta2} show the results of the experiments, similar to ones described for Figure~\ref{optimal_beta} before, that is, present the dependencies of optimal $\beta_*$ on smoothness characteristic $A$ and number of iterations $T$, but for non-diagonal preconditioner. In comparison with analogous dependencies for diagonal updates, optimal $\beta_*$ for non-diagonal updates are significantly less sensitive to the change of $A$ and $T$. In particular, dependence of $\beta_*$ on $A$ in non-diagonal case is closer to linear (because it is closer to exponential in logarithmic scale for $A$), so the previous remark on sensitivity of $\beta_*$ to the change of $A \ll 1$ ceases to be relevant. Similarly, the growth of $\beta_*$ with increasing $t$ is significantly slower than in diagonal case and is also less monotonic. Taking into account the range of $\beta$ values on both figures, one can say that $\beta_* = (0.965 \pm 0.005)$ independently on smoothness of the problem and the number of iterations. Thus, our hypothesis is that the faster the smallest eigenvalue (together with the largest one) of the scaled Hessian tends to $1$ with increasing number of iterations, the less dependent optimal $\beta_*$ is on the smoothness and number of iterations, which means in the extreme case that optimal $\beta_*$ is determined by some affine-invariant characteristic of the function. Note that such a $\beta_*$ can be greater than $\beta_*$ obtained for diagonal updates (cf. Figure \ref{optimal_beta}).
\begin{figure}[ht!]
\centering
     \includegraphics[width=0.4\textwidth]{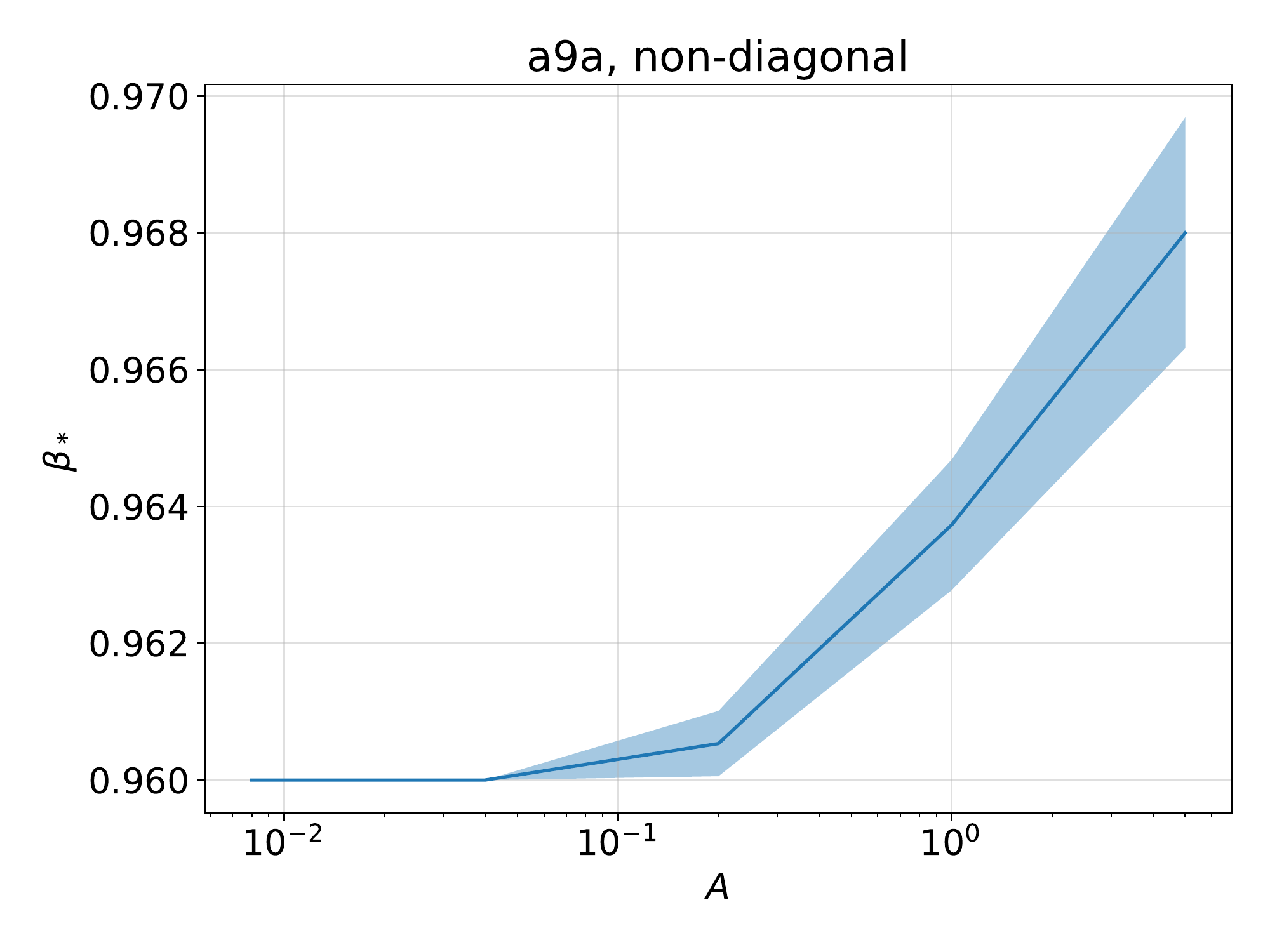}
     \includegraphics[width=0.4\textwidth]{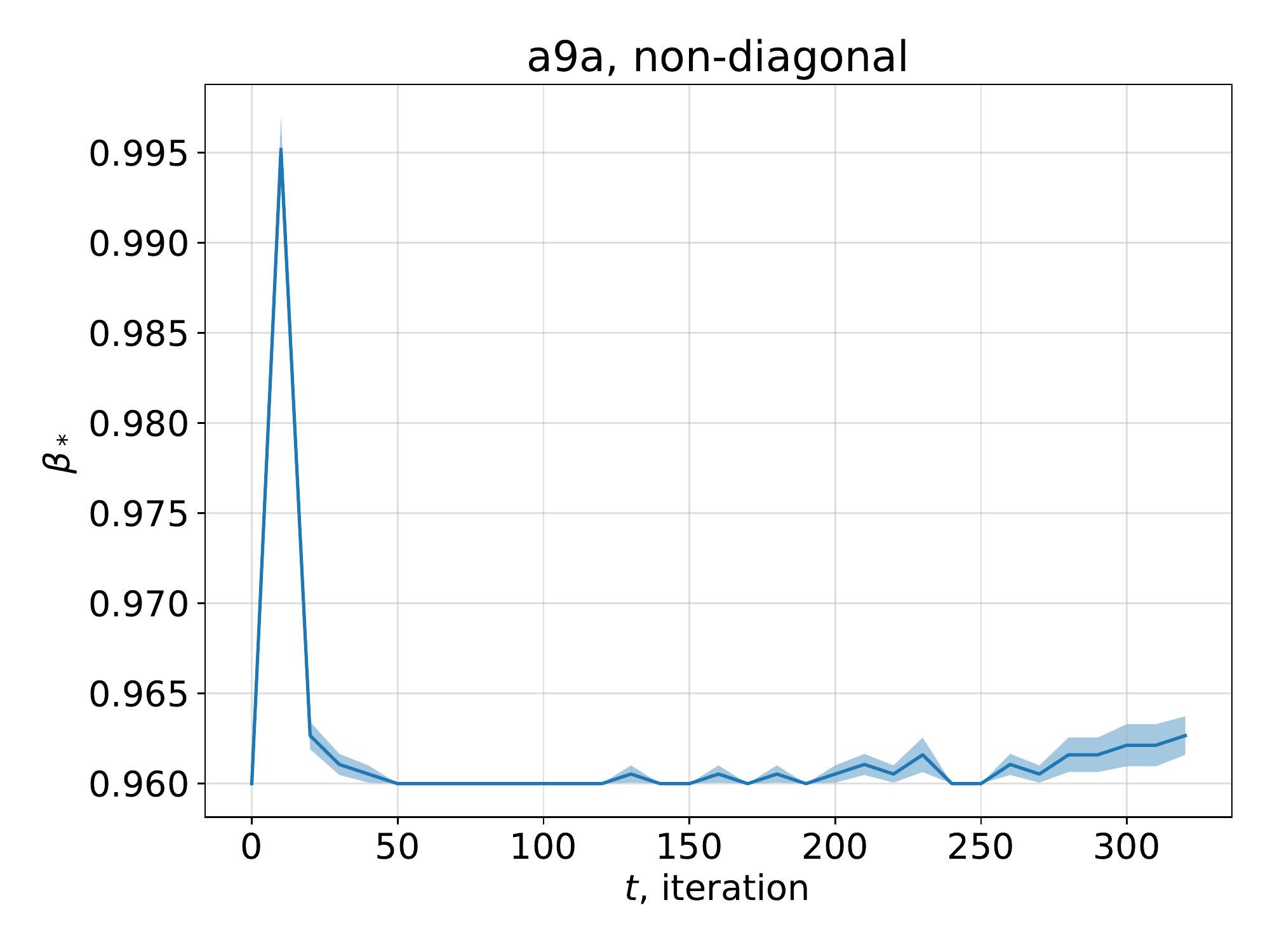}
        \label{optimal_beta_T2}
     \vskip-20pt
    \caption{Dependencies of optimal $\beta_t \equiv \beta$ for non-diagonal preconditioning.}
    \label{optimal_beta2}
\end{figure}

The last group of experiments is related to the version of \texttt{Scaled L-SVRG} with the step sizes chosen with line-search. We use Brent algorithm which searches for $\eta_t \in [0, 1]$ with minimal value of the $f(x_t - \eta_t P_t^{-1} g_t)$ by only function evaluations. This could be non-practical, if calculation of all the objective function's terms is computationally expensive, but in the case when the most expensive operation is evaluation of the gradient it is acceptable. Besides, our interest to \texttt{Scaled L-SVRG} with line-search is more theoretical~--- namely, by this modification we would like to reach the advantage of scaling introduced by averaging of the smoothness constants. Indeed, in previous experiments the step size was fixed, such that Lipschitz constant of the gradient was included in convergence rate as a true constant~--- scaled, but not better than by fixed preconditioner. Here, on the contrary, algorithm exploits scaling as much as possible. So, we compare the performance of algorithms with or without line-search to assess this advantage, and show the dependence of performance on $\beta$.

\begin{figure}[ht!]
\centering
    \vskip-12pt
     \includegraphics[width=0.4\textwidth]{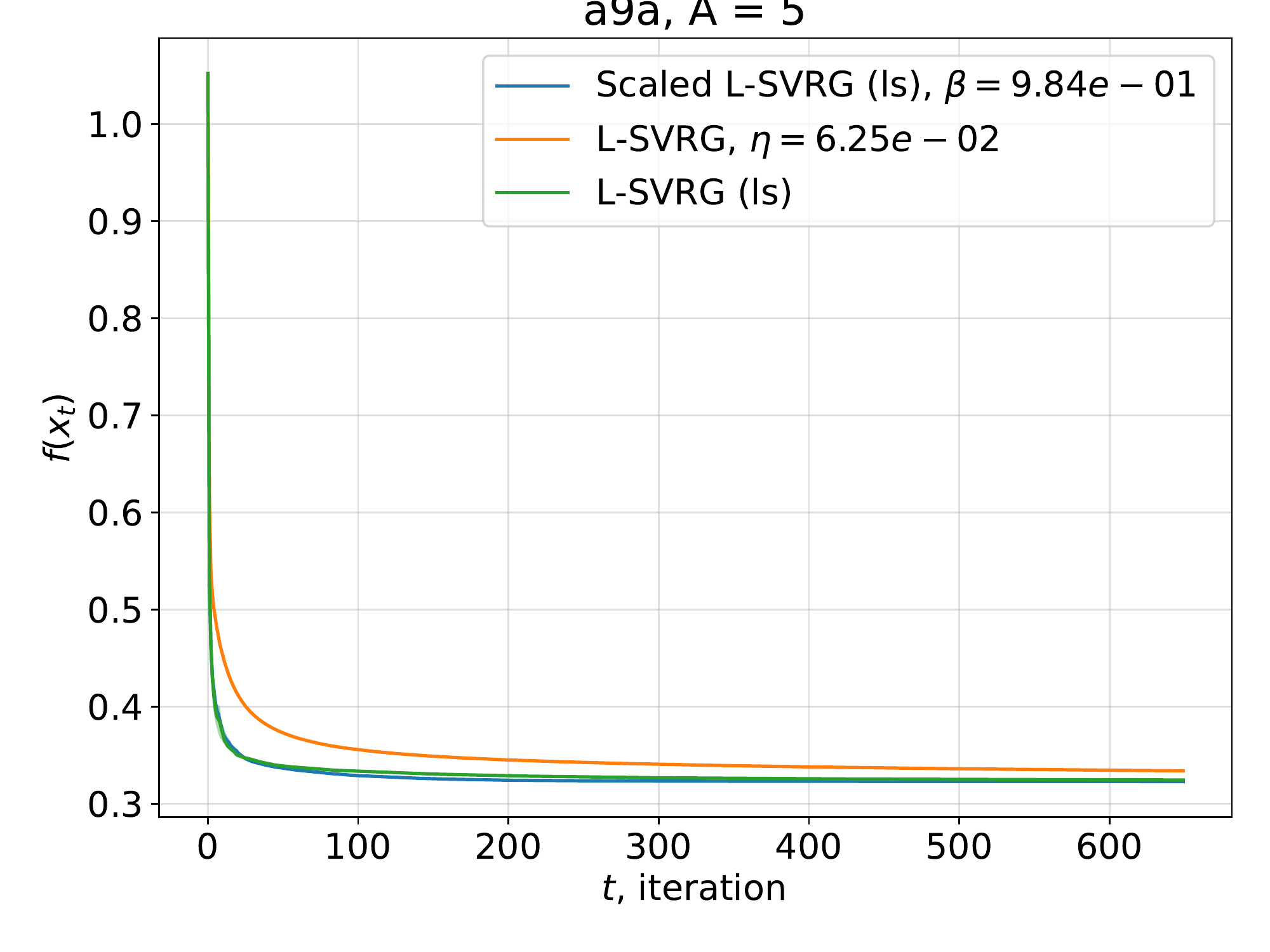}
     \includegraphics[width=0.4\textwidth]{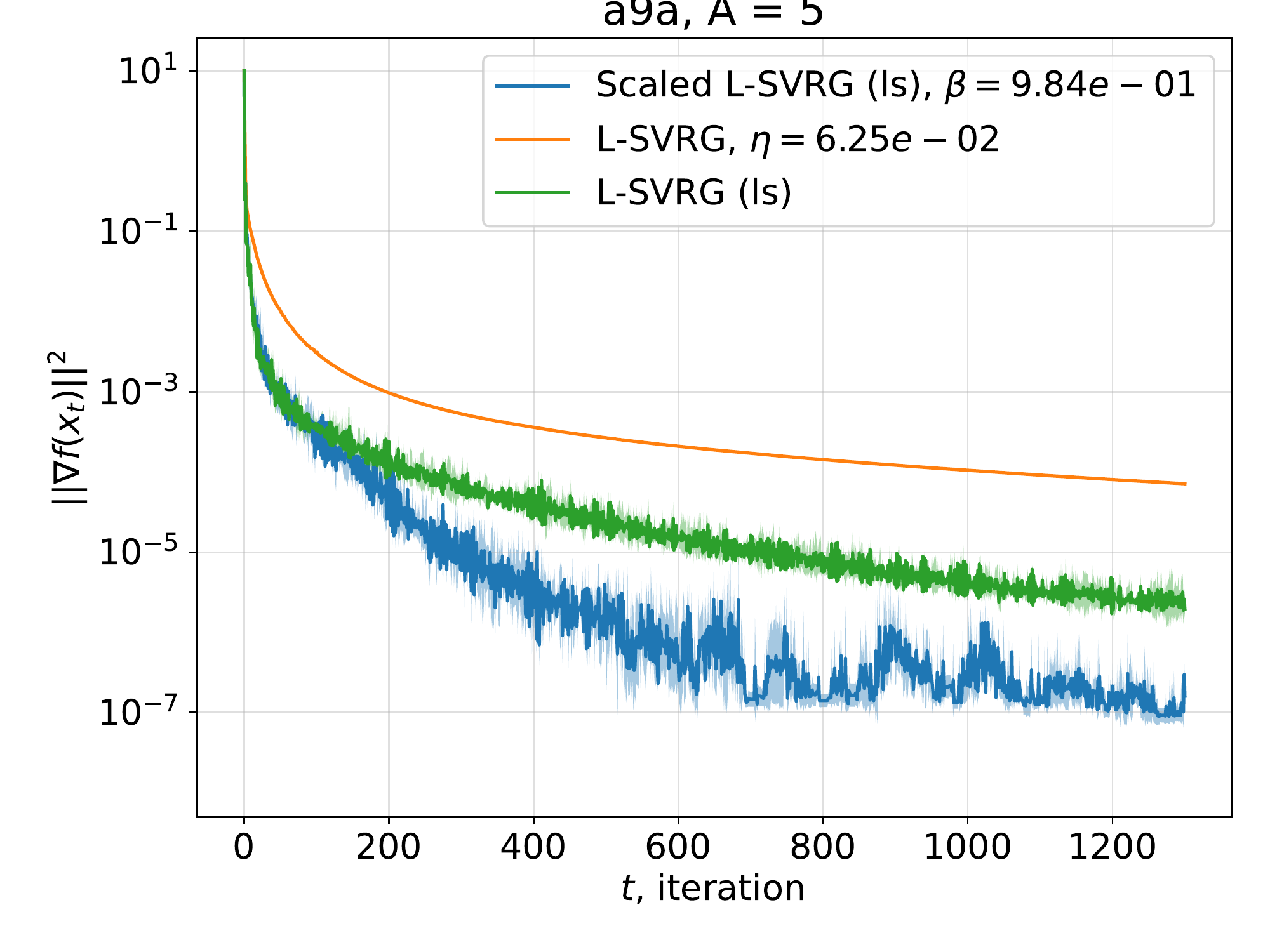}
     \vskip-20pt
    \caption{Convergence curves of \texttt{L-SVRG}, \texttt{L-SVRG} with line-search and \texttt{Scaled L-SVRG} with line-search with optimal choice of $\beta_t \equiv \beta$ and $\eta_t \in [0, 1]$.}
    \label{ls-conv}
    \vskip-12pt
\end{figure}

On the Figure~\ref{ls-conv}, we compare \texttt{L-SVRG}, \texttt{L-SVRG} with line-search (with the aim of fair comparison) and \texttt{Scaled L-SVRG} with line-search. Firstly, precision obtained by the algorithms with line-search is much better. At the same time, performance of \texttt{L-SVRG} with line-search and \texttt{Scaled L-SVRG} with line-search is almost the same until 200 iterations, and advantage of scaling comes clear with increasing of number of iterations and improving the preconditioner (see Figure~\ref{spectrum}). It is natural: the average of smoothness constants decreases with adapting of the preconditioner.
\begin{figure}[ht!]
\centering
     \includegraphics[width=0.4\textwidth]{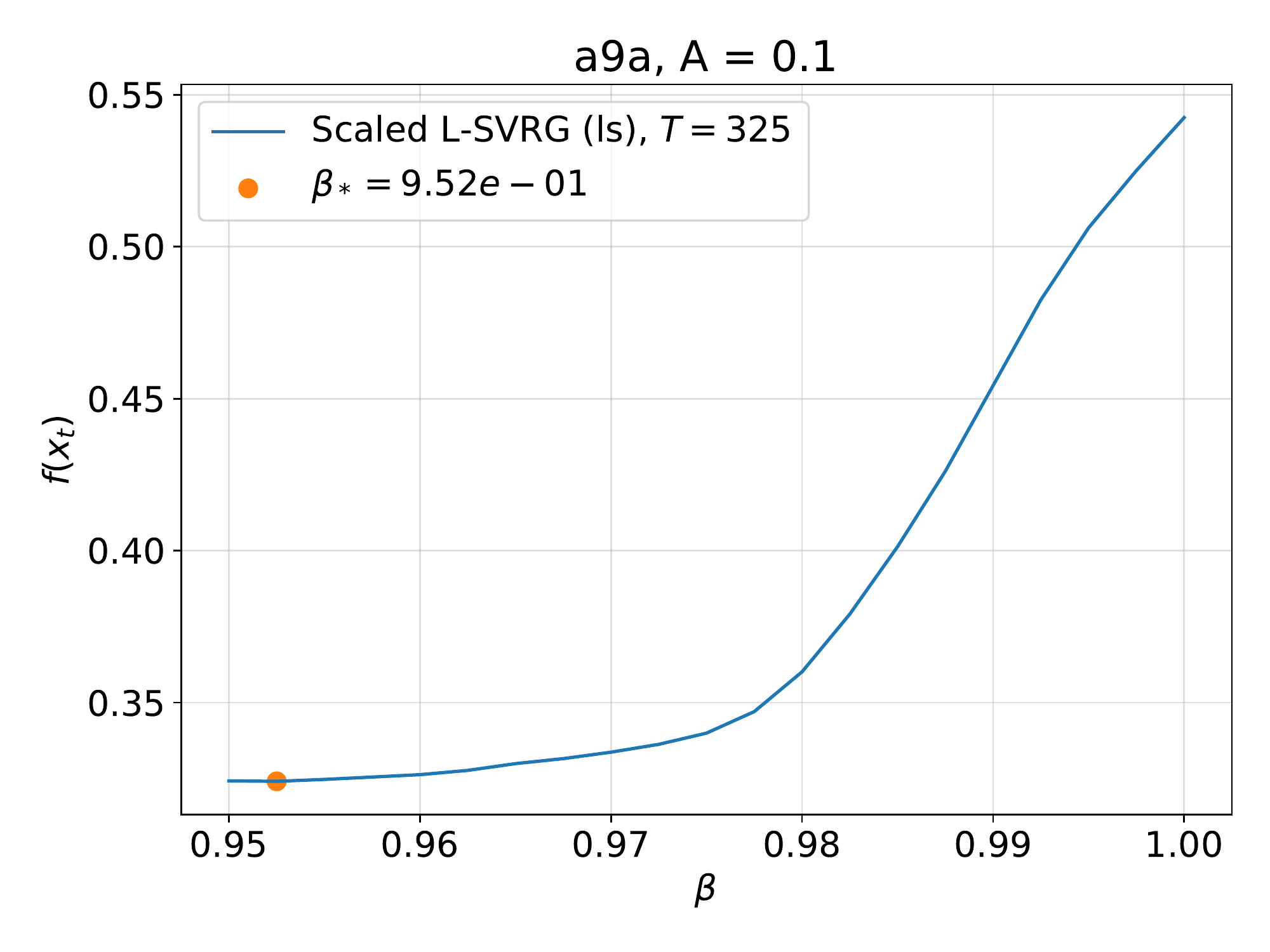}
     \includegraphics[width=0.4\textwidth]{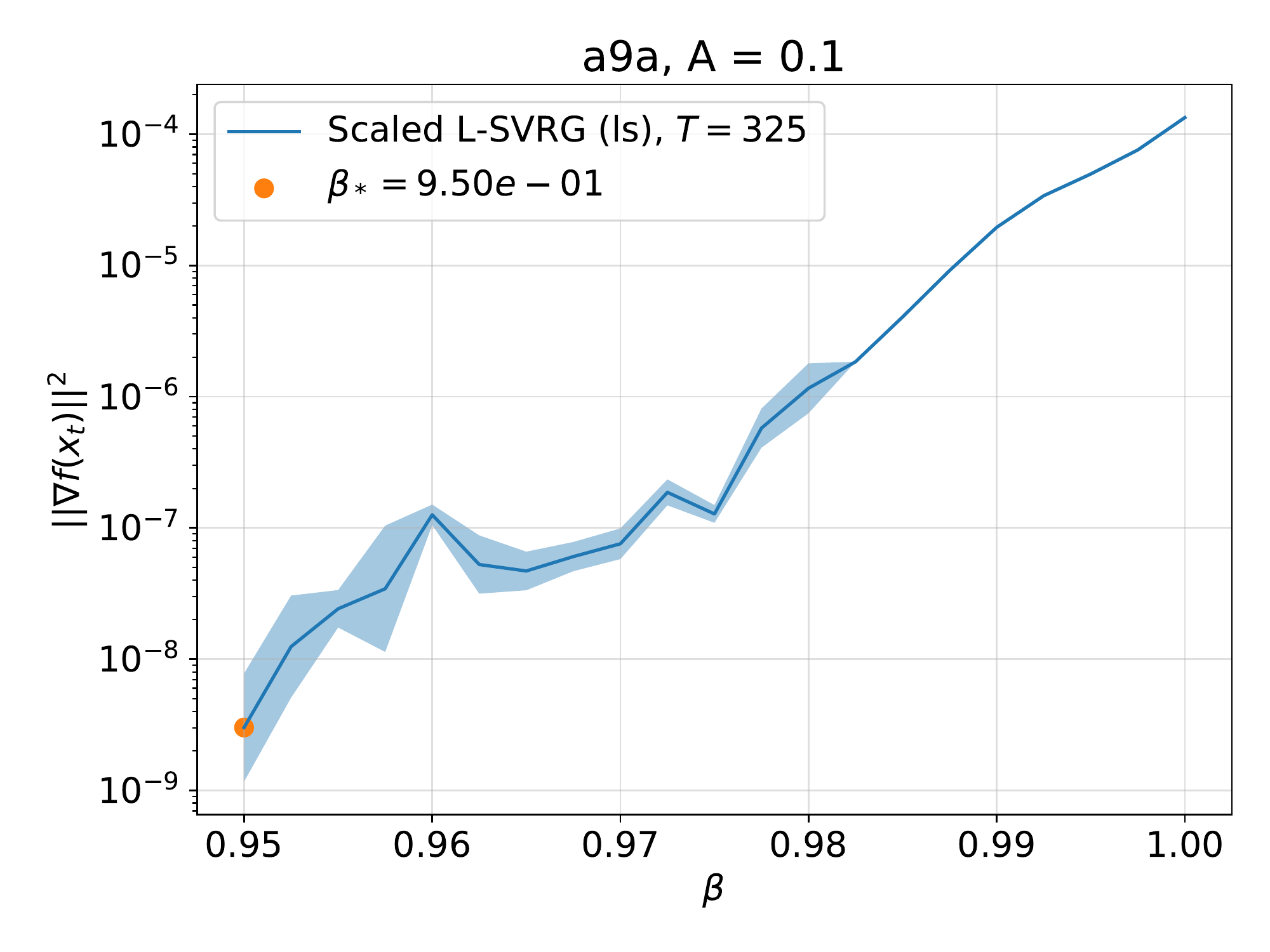}
     \includegraphics[width=0.4\textwidth]{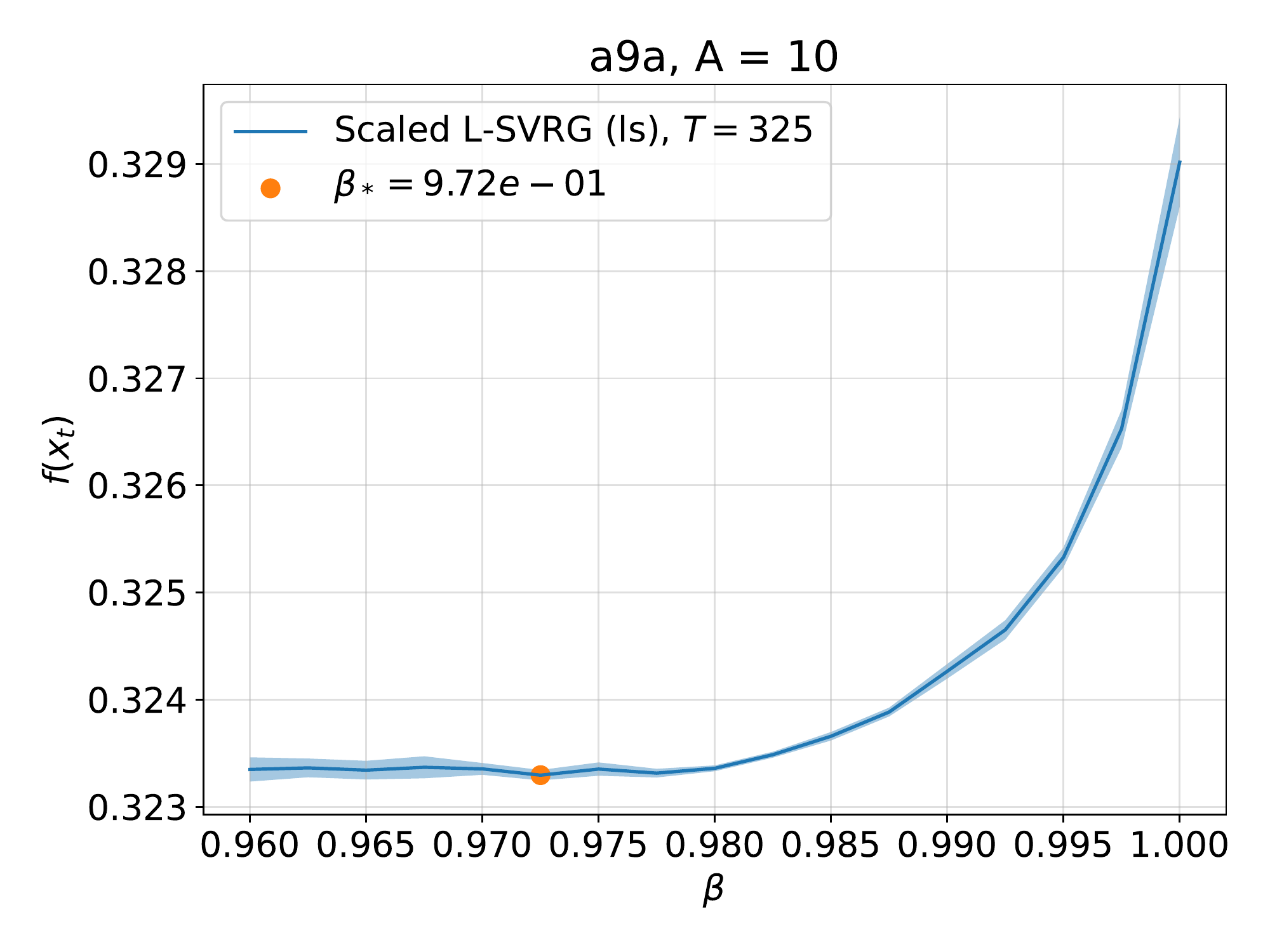}
     \includegraphics[width=0.4\textwidth]{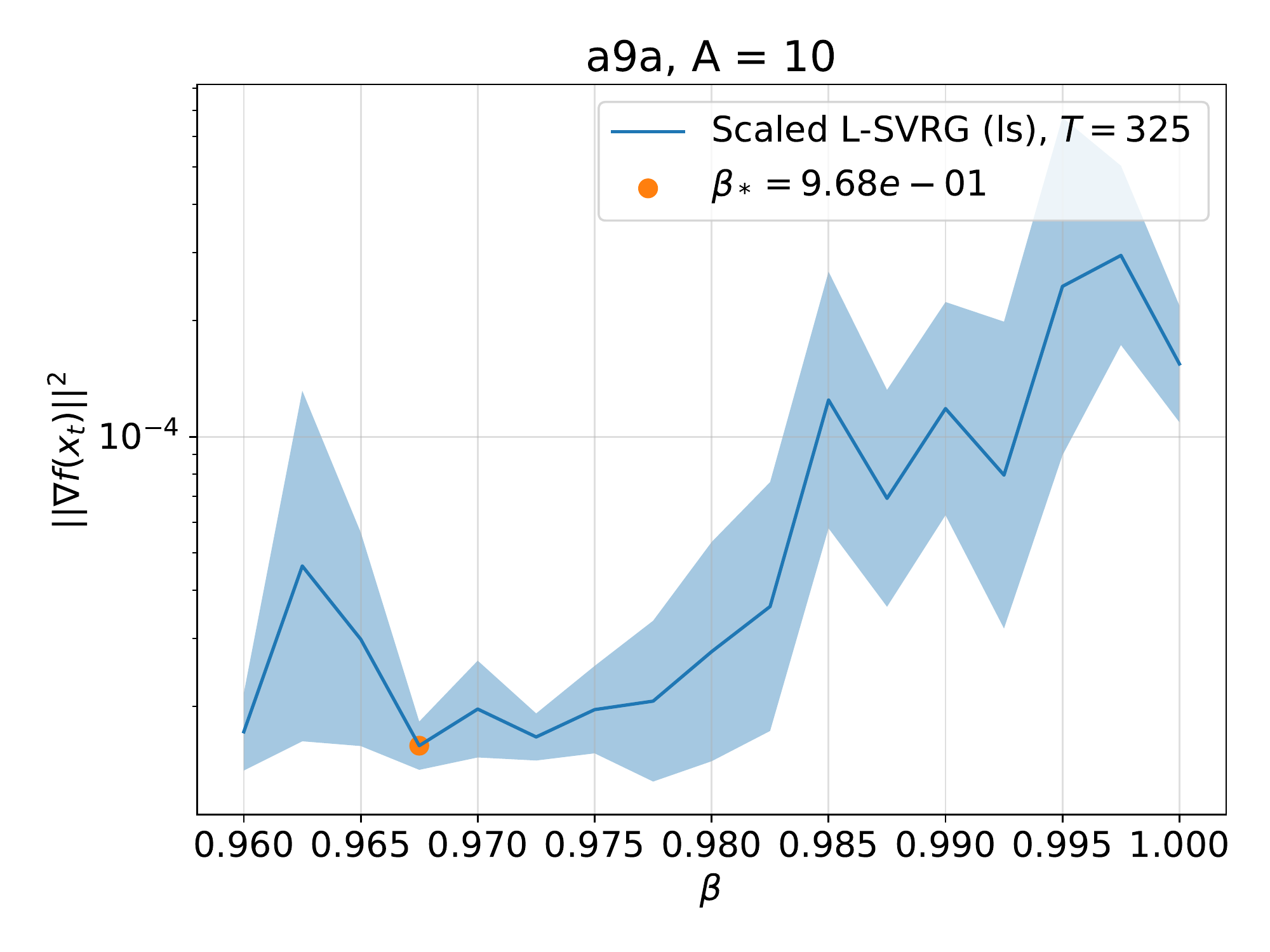}
     \vskip-20pt
    \caption{Dependence of achieved precision on $\beta_t \equiv \beta$ with line-search.}
    \label{betas2}
\end{figure}

Then, we reproduce the comparison of the \texttt{Scaled L-SVRG} operation in dependence on $\beta$. The results are shown on the Figure~\ref{betas2} (cf. Figure~\ref{betas_full}). Summarizing the differences, small $\beta$ values became acceptable even for big values of $A$, so that preconditioner can adapt faster without sacrificing convergence rate. This is a little unexpected, because the convergence slowdown for small $\beta$ values is explained primarily by the variance introduced by changing preconditioner, and the use of line-search does not relieve us of this factor. For now, we cannot explain this effect with certainty.
\newpage
\section{Omitted figures}\label{om_fig}
\subsection{Experiments for LibSVM \texttt{a9a} dataset}

\newpage

\newpage
\subsection{Experiments for LibSVM \texttt{covtype-binary-scaled} dataset}

\begin{figure}[ht!]
\centering
     \includegraphics[width=0.24\textwidth]{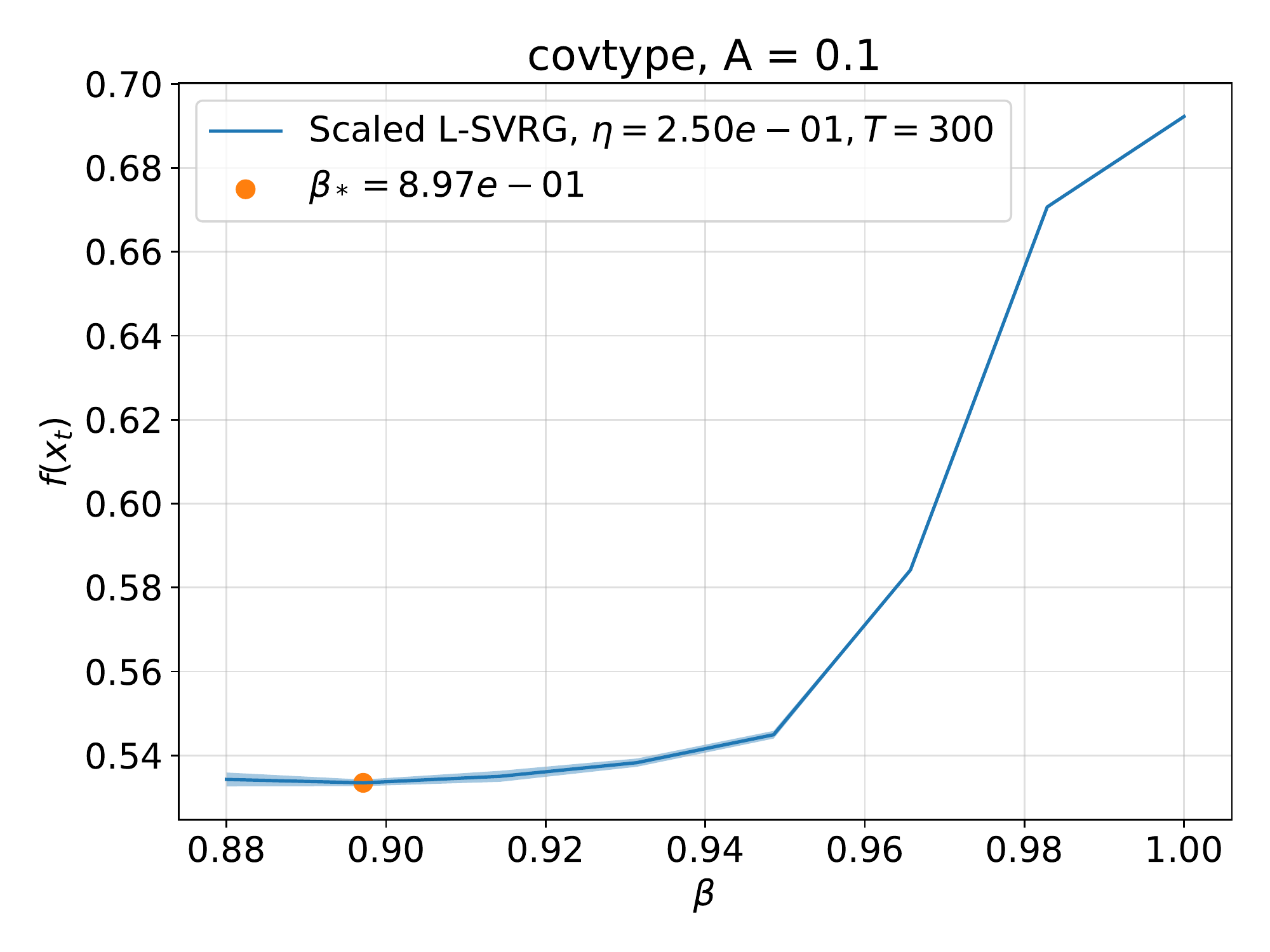}
     \includegraphics[width=0.24\textwidth]{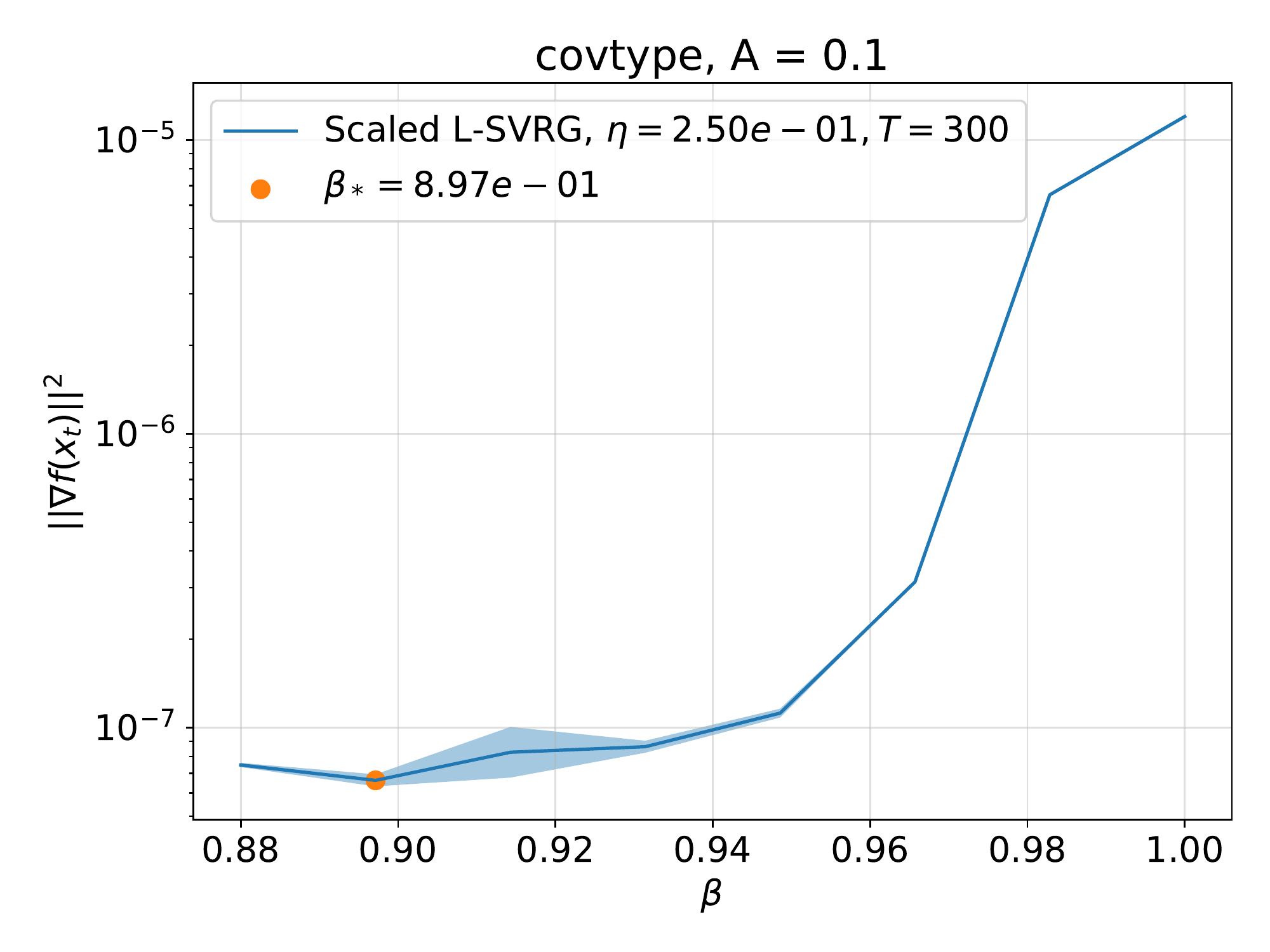}
     \includegraphics[width=0.24\textwidth]{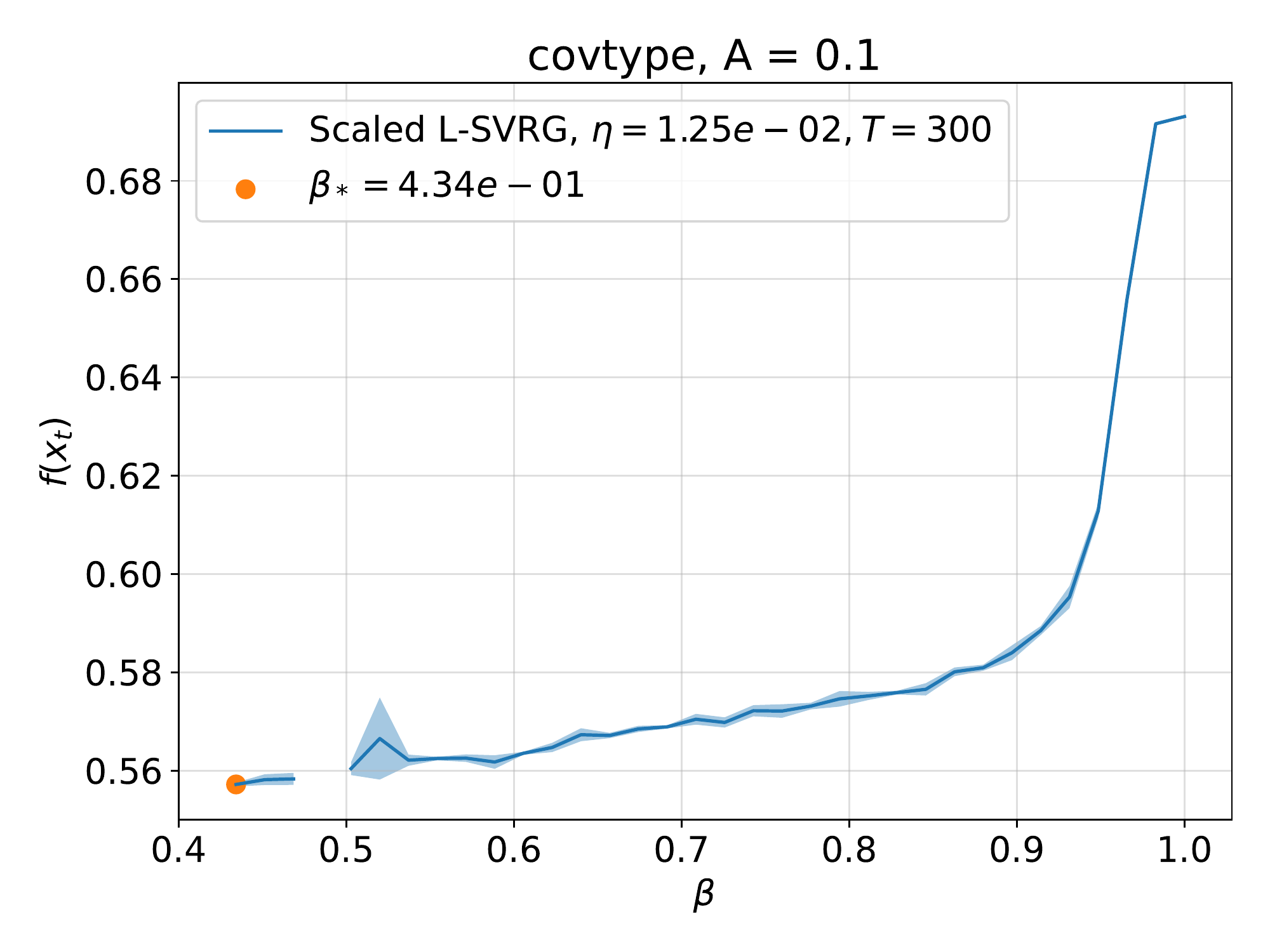}
     \includegraphics[width=0.24\textwidth]{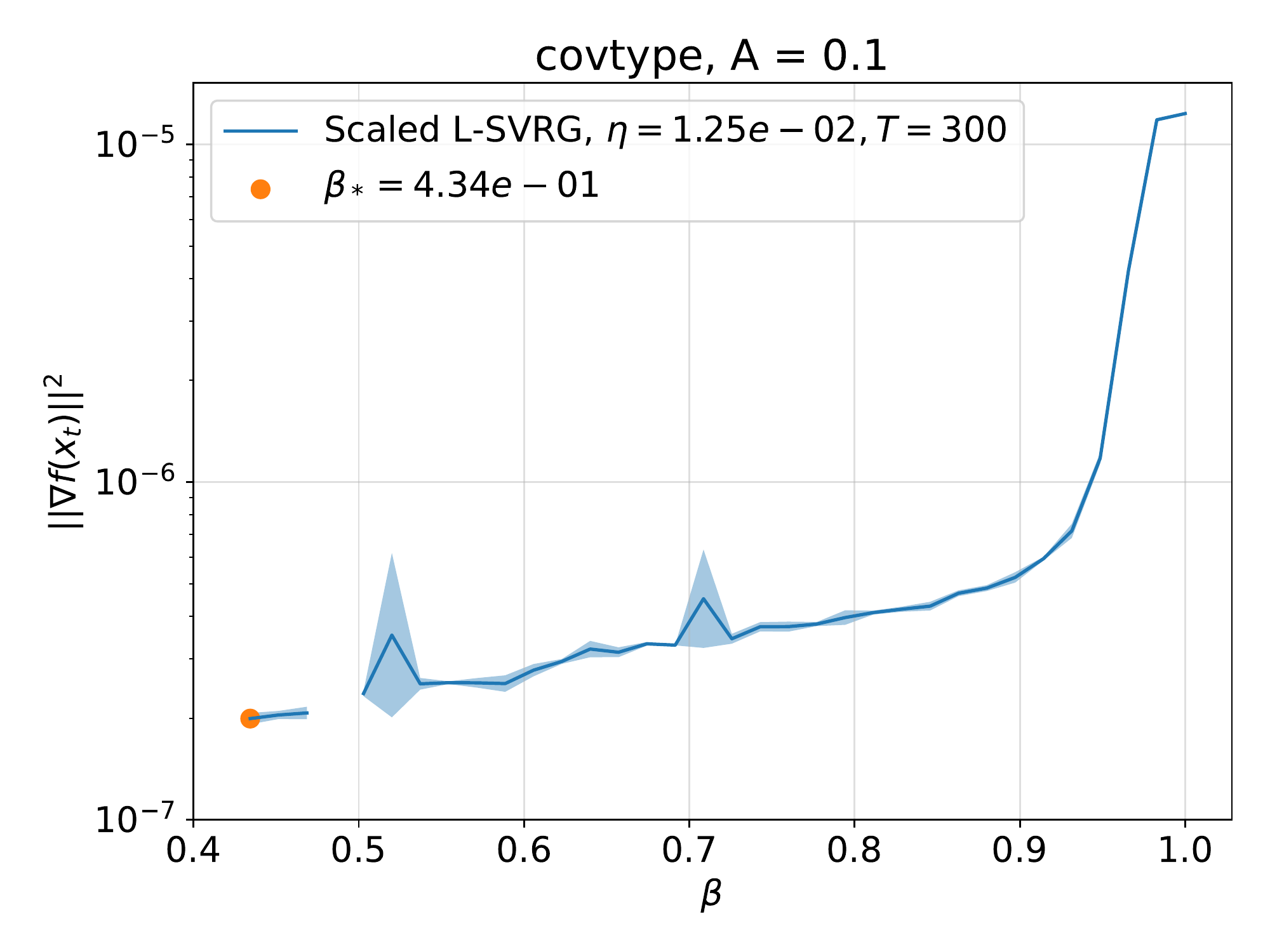}
    \vskip-12pt
    \caption{Dependence of achieved precision on $\beta_t \equiv \beta$, $A = 0.1$. (Note: gaps in curves mean the divergence of the algorithm in at least one of 3 runs)}
    \label{betas-covtype-Lsmall}
\end{figure}

\begin{figure}[ht!]
\centering
    \vskip-6pt
     \includegraphics[width=0.24\textwidth]{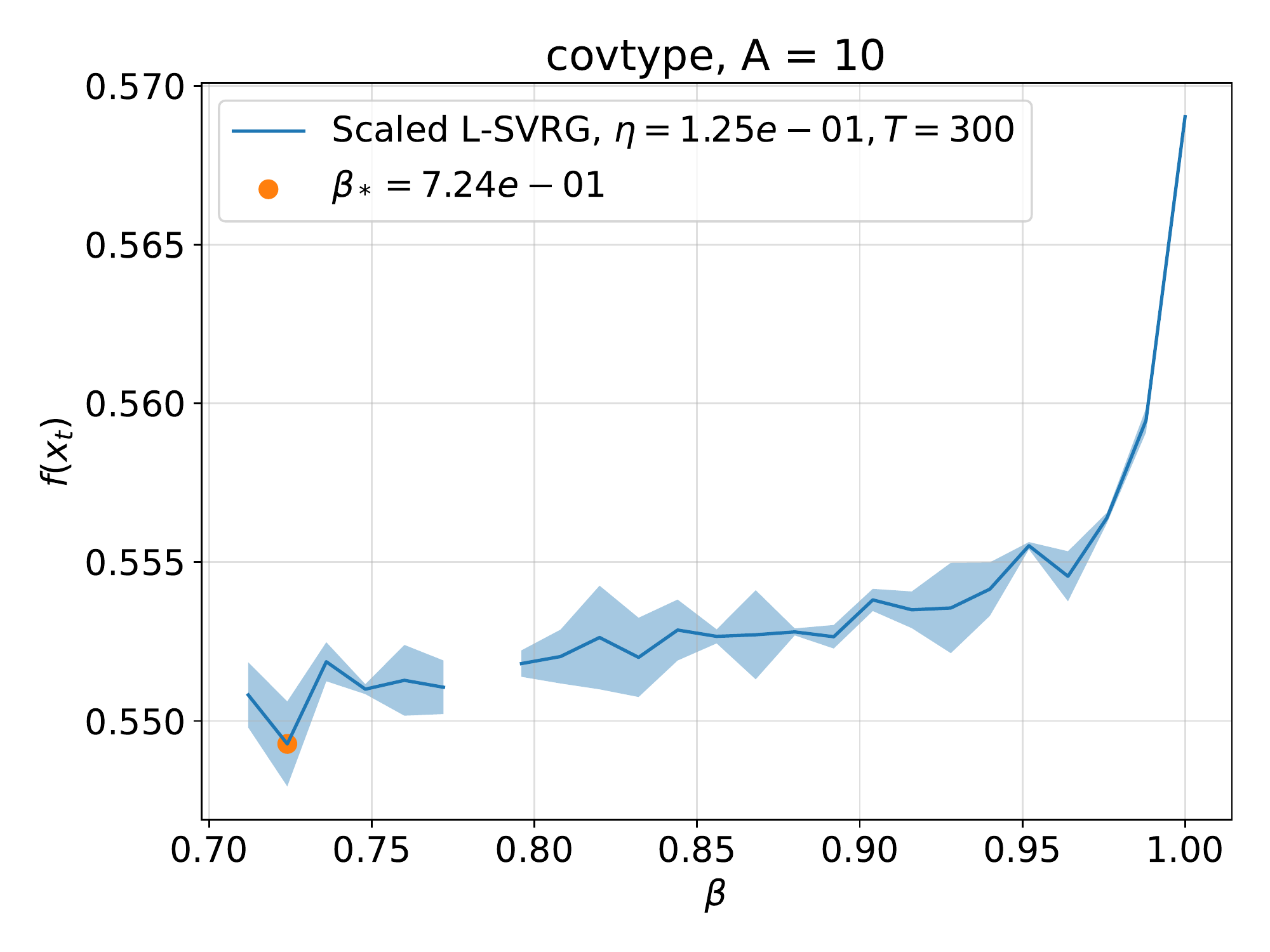}
     \includegraphics[width=0.24\textwidth]{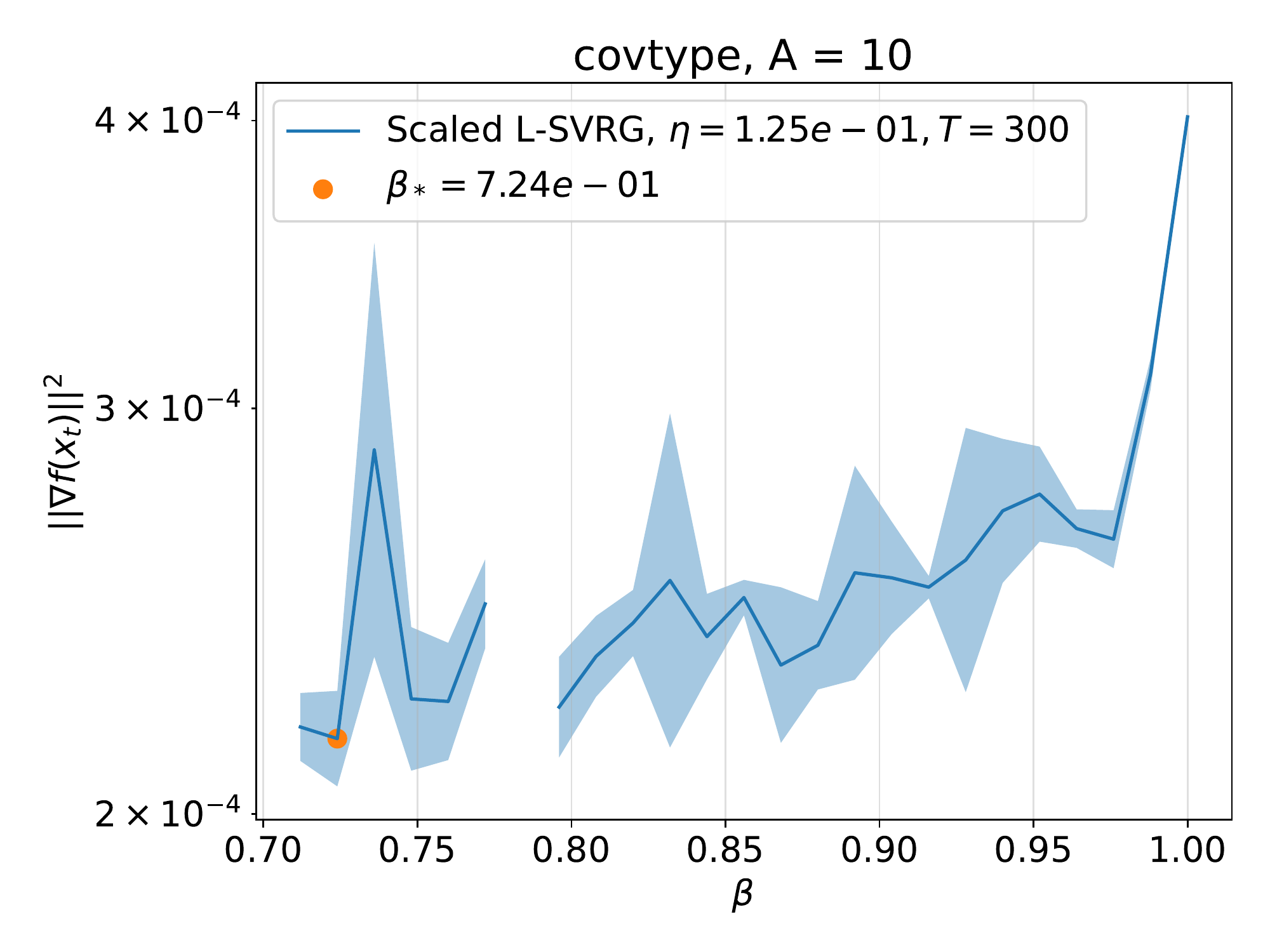}
     \includegraphics[width=0.24\textwidth]{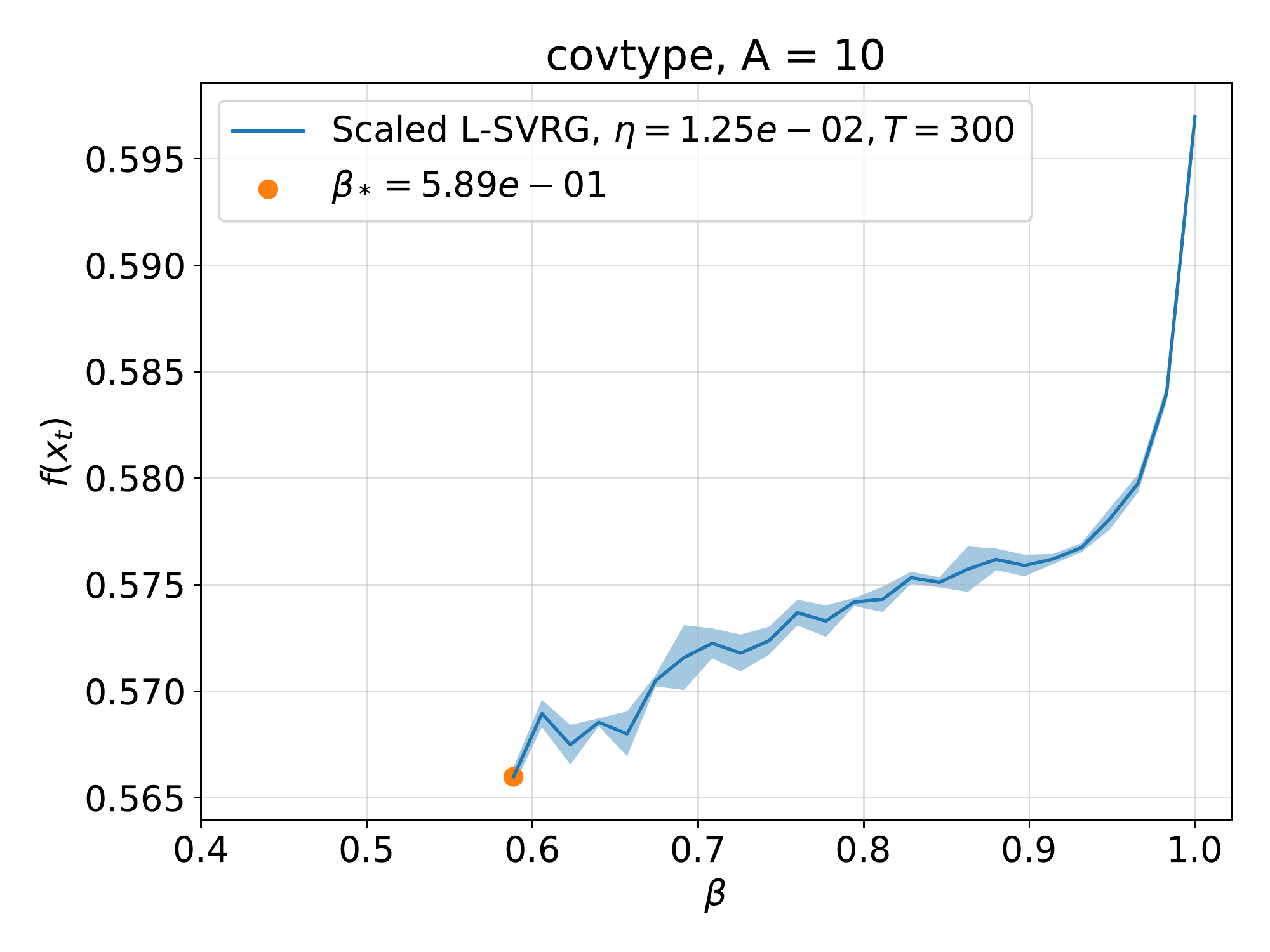}
     \includegraphics[width=0.24\textwidth]{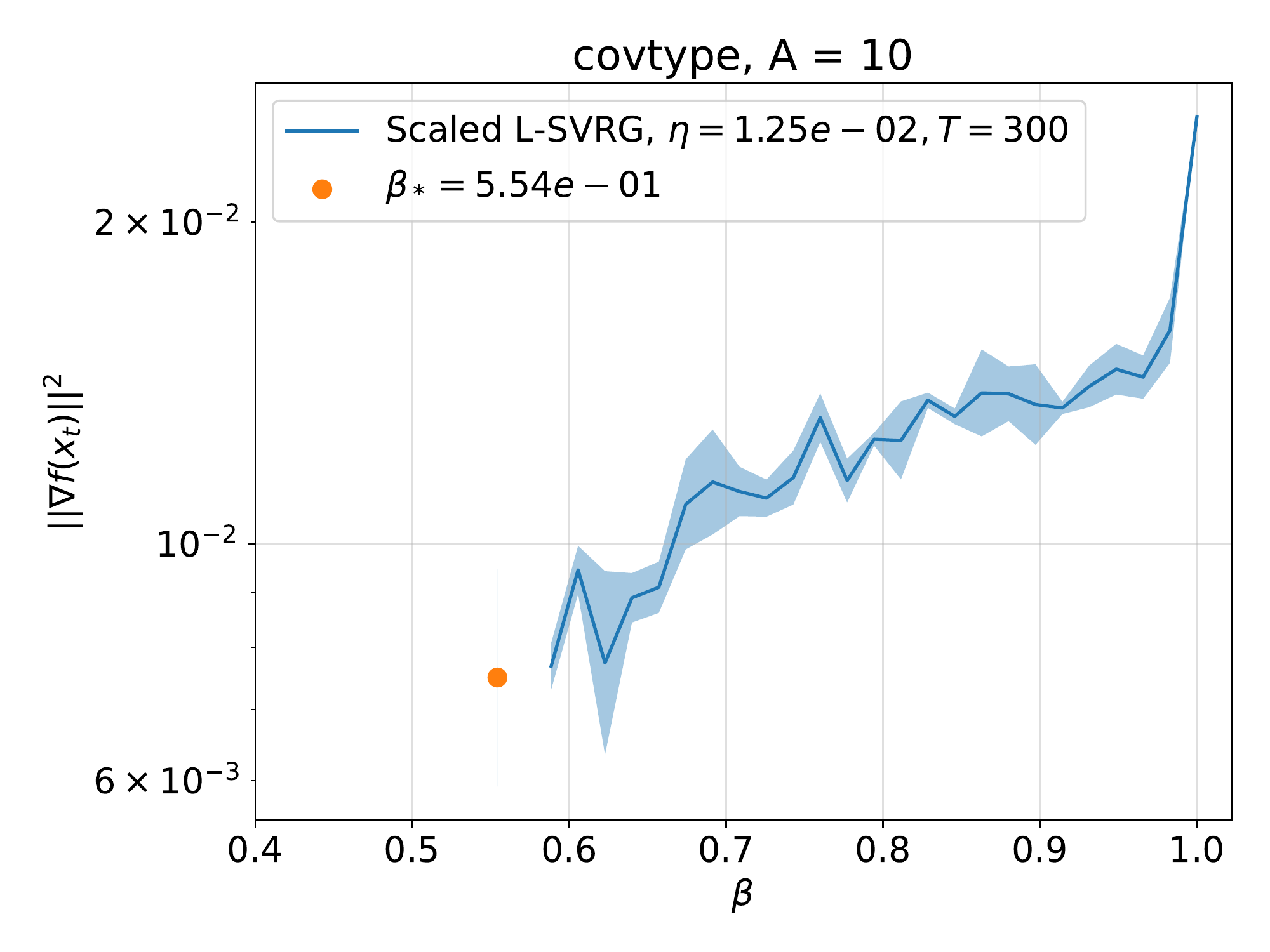}
    \vskip-12pt
    \caption{Dependence of achieved precision on $\beta_t \equiv \beta$, $A = 10$. (Note: gaps in curves mean the divergence of the algorithm in at least one of 3 runs)}
    \label{betas-covtype-Lbig}
\end{figure}

\begin{figure}[ht!]
\centering
    \vskip-6pt
     \includegraphics[width=0.3\textwidth]{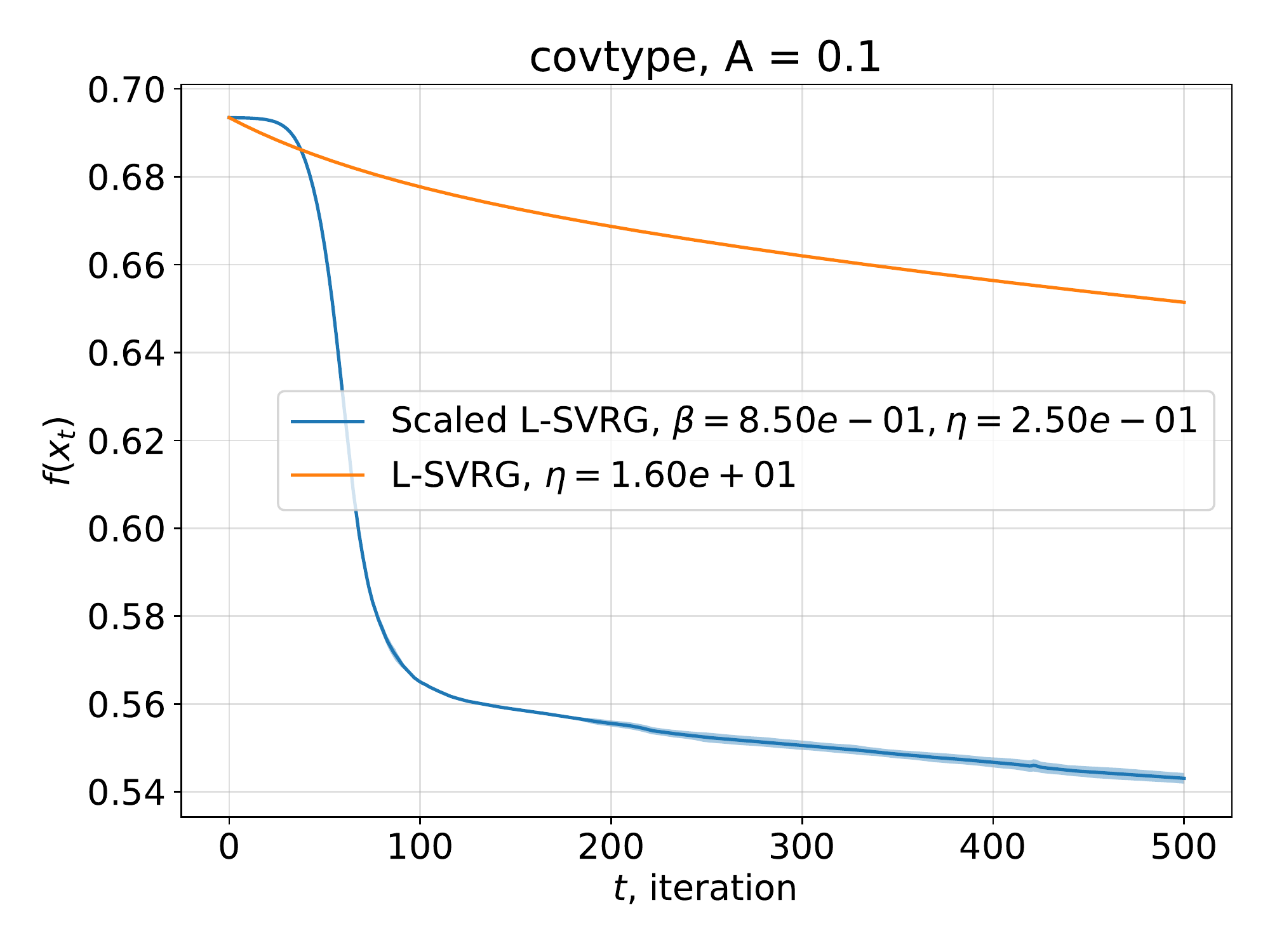}
     \includegraphics[width=0.3\textwidth]{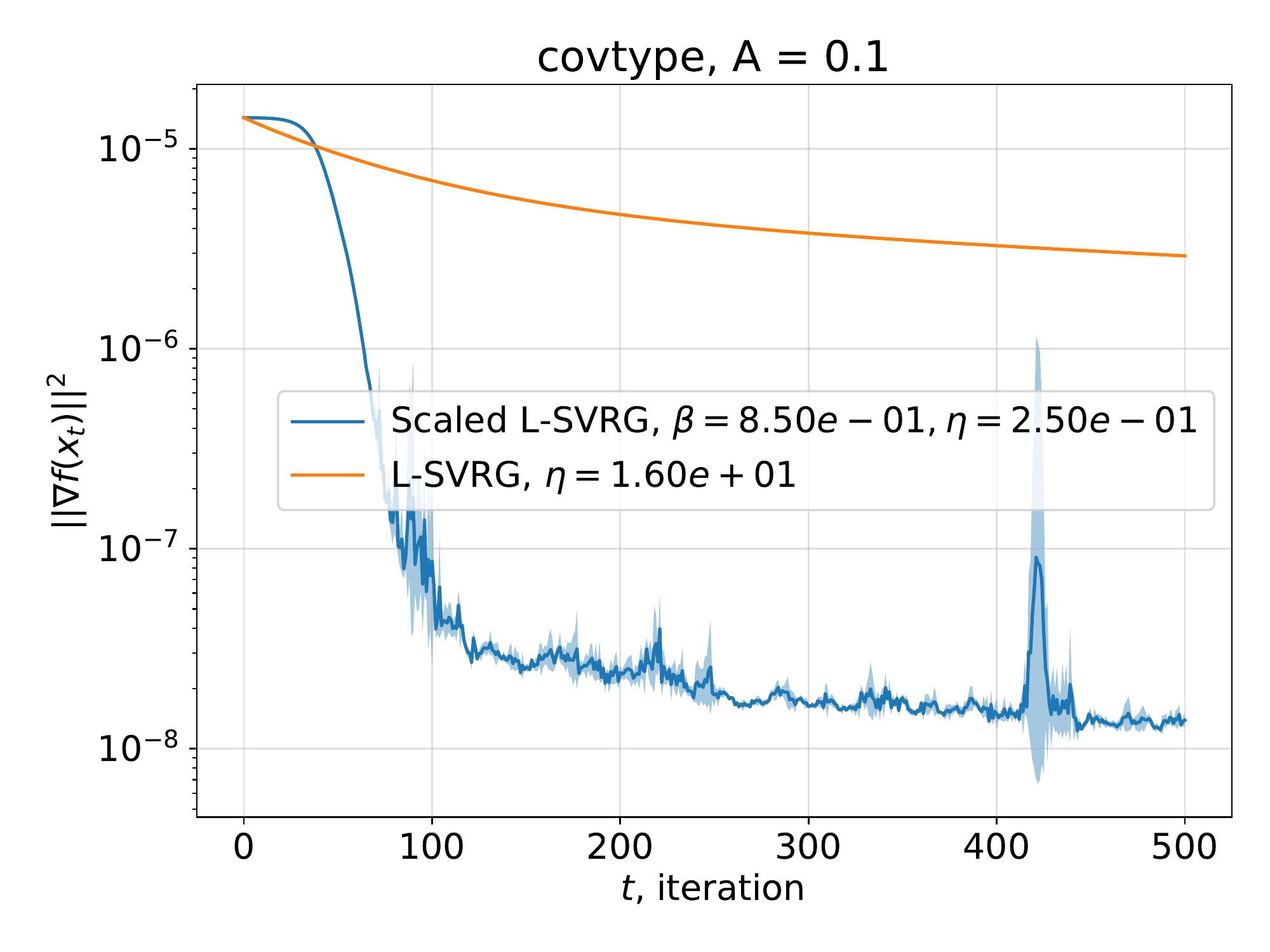}
     \includegraphics[width=0.3\textwidth]{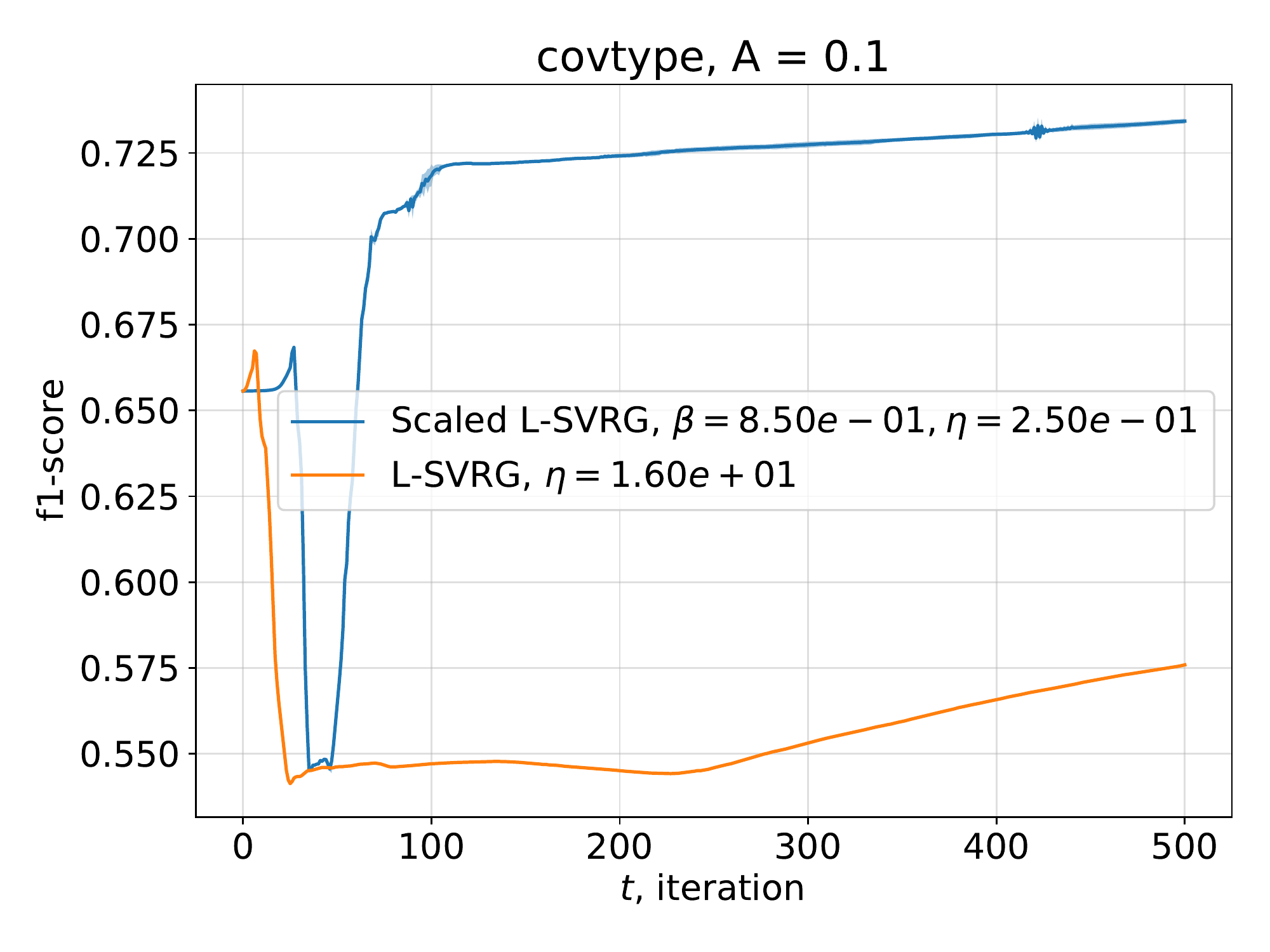}
     \includegraphics[width=0.3\textwidth]{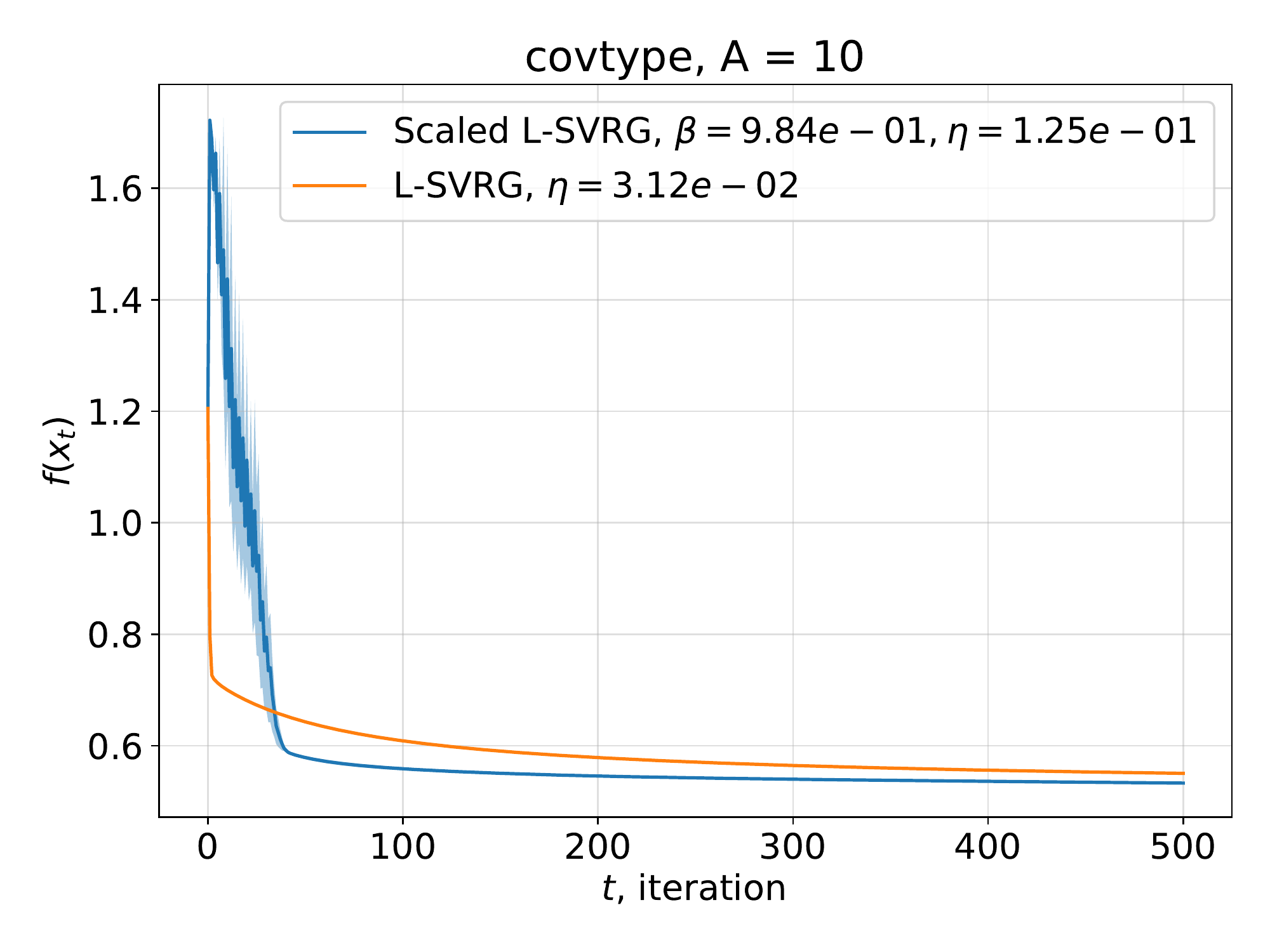}
     \includegraphics[width=0.3\textwidth]{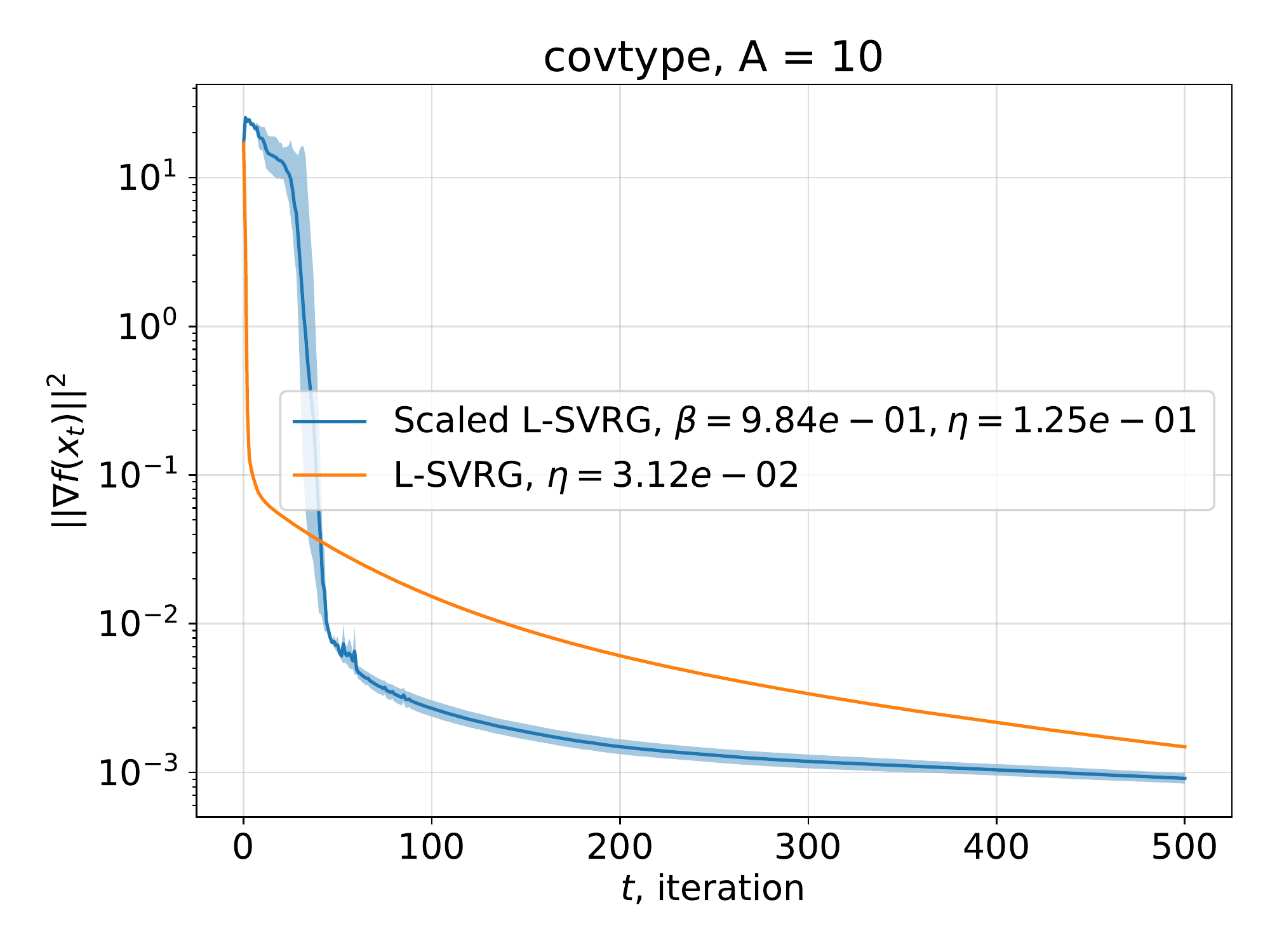}
     \includegraphics[width=0.3\textwidth]{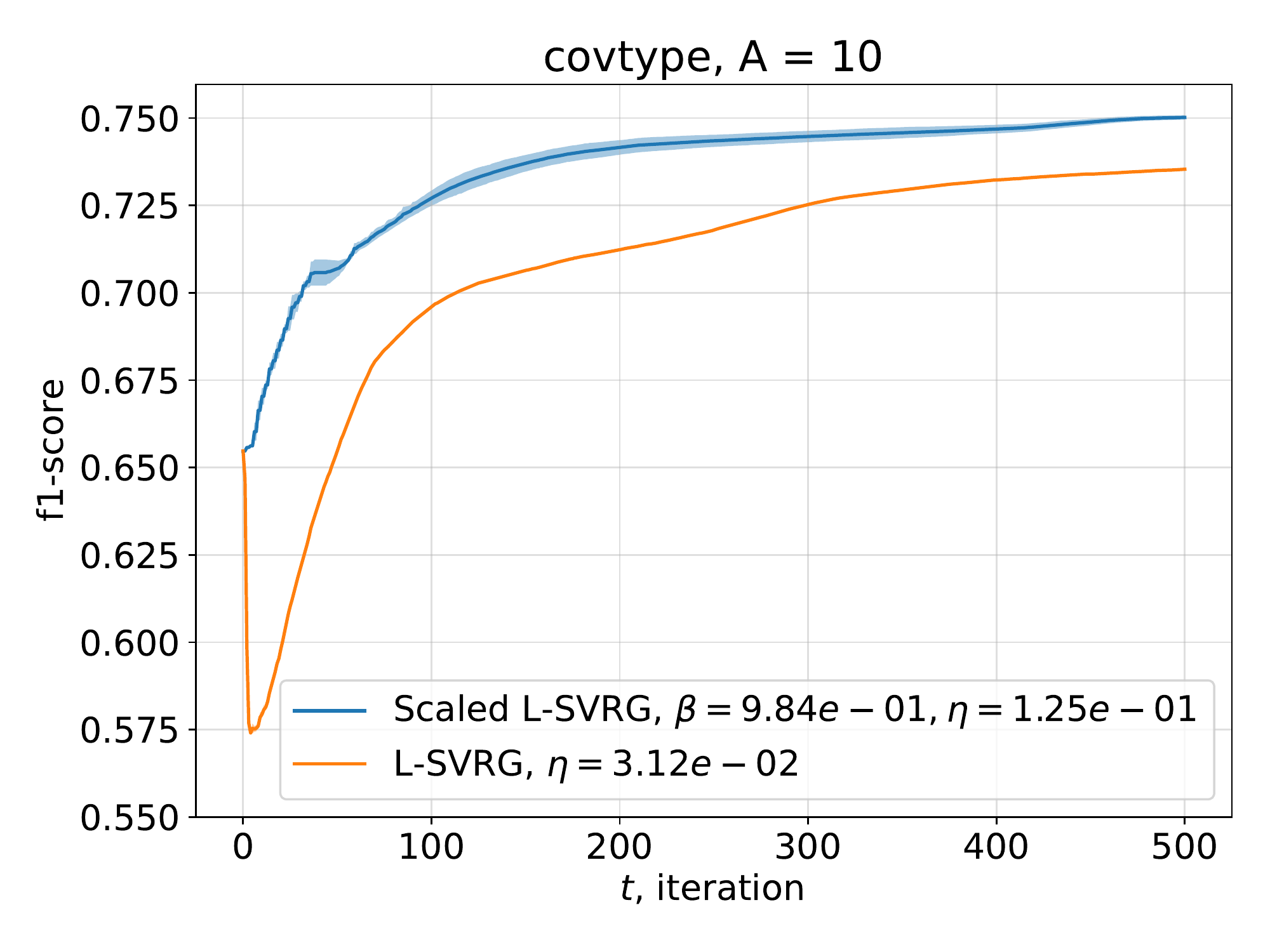}
    \vskip-12pt
    \caption{Convergence curves of \texttt{L-SVRG} and \texttt{Scaled L-SVRG} with optimal choice of $\beta_t \equiv \beta$ and $\eta_t \equiv \eta$.}
    \label{conv-covtype}
\end{figure}

\end{document}